\documentclass{amsart}

\usepackage{amsmath,amsthm,amssymb,amsfonts}
\usepackage[a4paper]{geometry}
\usepackage{a4wide}
\usepackage{graphicx}
\usepackage{mathpazo}
\usepackage[T1]{fontenc}
\usepackage[utf8]{inputenc}
\usepackage{enumitem}
\usepackage{hyperref}
\usepackage{verbatim}
\usepackage{xcolor}
\usepackage{todonotes}

\hypersetup{
    colorlinks,
    linkcolor={red!80!black},
    citecolor={blue!80!black},
    urlcolor={blue!80!black}
}
\usepackage{tikz}

\usetikzlibrary{matrix,arrows, cd, calc}  
\usepackage[margin=10pt,font=normalsize,labelfont=bf, labelsep=endash]{caption}

\makeatletter
\def\thm@space@setup{\thm@preskip=4pt
\thm@postskip=2pt}
\makeatother

\newtheorem{theorem}{Theorem}[section]

\newtheorem{corollary}[theorem]{Corollary}
\newtheorem{fact}[theorem]{Fact}

\newtheorem{proposition}[theorem]{Proposition}

\newtheorem{lemma}[theorem]{Lemma}

\theoremstyle{definition}
\newtheorem{definition}[theorem]{Definition}
\newtheorem{notation}[theorem]{Notation}

\newtheorem{example}[theorem]{Example}

\newtheorem{remark}[theorem]{Remark}

\newcommand{\ZZ}{{\mathbb Z}}
\newcommand{\NN}{{\mathbb N}}

\newcommand{\cO}{{\mathcal O}}

\newcommand{\cM}{{\mathcal M}}

\newcommand{\cH}{{\mathcal H}}

\newcommand{\cB}{{\mathcal B}}

\newcommand{\goodquotient}{\mathbin{
  \mathchoice{\left/\mkern-6mu\right/}
    {/\mkern-5mu/}
    {/\mkern-5mu/}
    {/\mkern-5mu/}}}

\DeclareMathOperator{\Hom}{Hom}
\DeclareMathOperator{\Spec}{Spec}

\DeclareMathOperator{\Proj}{Proj}
\DeclareMathOperator{\Sym}{Sym}
\DeclareMathOperator{\im}{im}
\DeclareMathOperator{\Ext}{Ext}

\DeclareMathOperator{\Mor}{Mor}

\newcommand{\Set}{Set}

\newcommand{\Artk}{\textit{Art}_{\kk}}

\newcommand{\Gmult}{\mathbb{G}_m}%
\newcommand{\kk}{\mathbf{k}}%
\newcommand{\varX}{X}%
\newcommand{\Xplus}{\varX^+}%
\newcommand{\ione}[1]{i_{#1}}%
\newcommand{\ioneX}{\ione{\varX}}%
\newcommand{\iinfty}[1]{\pi_{#1}}%
\newcommand{\iinftyX}{\iinfty{\varX}}%
\newcommand{\isection}[1]{s_{#1}}%
\newcommand{\isectionX}{\isection{\varX}}%
\newcommand{\BBname}{Bia{\l}ynicki-Birula}%
\DeclareMathOperator{\Hilb}{Hilb}%

\newcommand{\Affine}{\mathbb{A}}%
\newcommand{\G}{\mathbb{G}}%
\newcommand{\Pone}{\mathbb{P}^1}%
\newcommand{\Pn}{\mathbb{P}^n}%
\newcommand{\mm}{\mathfrak{m}}%
\newcommand{\Fperp}{F^{\perp}}%
\newcommand{\eps}{\varepsilon}%
\newcommand{\Fperpx}{F_{\Bar{x}}^{\perp}}%
\newcommand{\Fperpy}{F_{\Bar{y}}^{\perp}}%
\newcommand{\xx}{\Bar{x}}%
\newcommand{\yy}{\Bar{y}}%
\newcommand{\aalpha}{\Bar{\alpha}}%
\newcommand{\bbeta}{\Bar{\beta}}%
\newcommand{\wcirc}{\mathbin{\lrcorner}}%
\newcommand{\Mf}{fractal family }%

\newcommand{\symbf}{f}

\newcommand{\symba}{\eta}
\newcommand{\symbb}{\mu}

\newcommand{\symbJo}{I}
\newcommand{\sgamma}{\gamma}
\newcommand{\symbT}{T}
\newcommand{\symbJ}{J}
\newcommand{\nameX}{\mathcal{X}}
\newcommand{\nameY}{\mathcal{Y}}
\newcommand{\nameq}{u}
\newcommand{\nameD}{\mathcal{D}}
\newcommand{\namep}{c}

\newcommand{\symbv}{\alpha}
\newcommand{\symbp}{S}

\newcommand{\symbc}{\delta}

\newcommand{\symb}{\gamma}

\newcommand{\symbtan}{\mathfrak{x}}
\newcommand{\symboly}{\mathfrak{y}}
\newcommand{\symbphi}{\phi}

\begin{document}

\title{Non-reducedness of the Hilbert schemes of few points}
\author{Michał Szachniewicz} 
\keywords{Algebraic geometry, Cycles and subschemes, Parametrization (Chow and Hilbert schemes)}
\address{Institute for Advanced Study, 1 Einstein Drive, Princeton, New Jersey, 08540, USA}
\email{\url{michals@ias.edu}}

\begin{abstract}
    We use generalised \BBname{} decomposition, apolarity and obstruction theories to prove non-reducedness of the Hilbert scheme of $13$ points on $\Affine^6$. Our argument doesn't involve computer calculations and gives an example of a fractal-like structure on this Hilbert scheme.
\end{abstract}
\maketitle

\tableofcontents

\section{Introduction}
Given a space $X$, many questions about its geometry can be reduced to the behaviour of families of subsets of $X$. If $X$ has some geometric structure, for example is a complex manifold, understanding families of subvarieties can lead to inductive arguments reducing problems about $X$ to lower dimensions. The simplest subvarieties of a variety are finite subsets, say of fixed cardinality $d$. If $d=1$ these are parametrized by $X$ itself. For higher $d$, one can consider configuration spaces $C_d(X) =(X^d \setminus \Delta)/S_d$, where $\Delta$ is the thick diagonal (set of tuples where any two variables are equal), and $S_d$ is the symmetric group acting by permutations. However, $C_d(X)$ lacks the information on how tuples of points can degenerate, when they come closer and closer to a single point of $X$ (or to a point of $X^d$ in $\Delta$). To add this information, one has to partially compactify $C_d(X)$, and Hilbert schemes of points form such a partial compactification (and are compact, if $X$ is). 
More precisely, the Hilbert scheme $\Hilb_{\Affine^n/\mathbb{C}}^{d}$ of $d$ points on $\Affine^n$ is a scheme whose $A$-points, where $A$ is a $\mathbb{C}$-algebra, correspond to ideals $I \subseteq A[x_1, \dots, x_n]$ such that the quotient $A[x_1, \dots, x_n]/I$ is a flat $A$-module of rank $d$.
The configuration space $C_d(\Affine^n)$ embeds into $\Hilb_{\Affine^n/\mathbb{C}}^{d}$, by mapping a tuple of $d$ distinct points to the ideal consisting of functions that vanish on all of them. However, there are $\mathbb{C}$-points of $\Hilb_{\Affine^n/\mathbb{C}}^{d}$, like the ideal $(x_1^d, x_2, \dots, x_n)$ in $\mathbb{C}[x_1, \dots, x_n]$, that do not arise in this fashion. In this text we study pathologies (namely non-reducedness) that can occur around some points outside of $C_d(\Affine^n)$, when $d=13$ and $n=6$.

Hilbert schemes of points have many applications in various areas of mathematics. For example, Haiman proved the $n!$ conjecture using their geometry (see \cite{Haiman}). Fogarty \cite{Fogarty} considered these schemes when the ambient variety is a smooth surface and showed that then the Hilbert scheme of points is smooth as well. If we start from a three or higher dimensional scheme this is not the case and many questions arise, such as when the Hilbert scheme of points is reduced, Cohen-Macaulay, etc. Recently, Jelisiejew closed a long-standing open question in this area and showed that the Hilbert scheme of points on an affine space is not reduced \cite{Hilbpath}; this problem was open for more than fifty years. Still the question about the smallest counter-example to reducedness remains. Methods of Jelisiejew were partially implicit and the non-reducedness witness was represented by a scheme of degree $5082$ on $\Affine^{14}$. In \cite{Jel24b}, using computer calculations, a different example was found, namely it was proven that $\Hilb_{\Affine^4/\kk}^{21}$ has generically non-reduced components. Often, people analyse the Hilbert schemes of $d$ points, where $d$ is small, say $d \leq 25$; see for example \cite{example1} or \cite{example2}. Results on the pathologies of low degree Hilbert schemes are especially important, as they can provide test grounds for various conjectures \cite{maulik2006gromov, Ramkumar_2023}. See also \cite{Ricolfi24, jelisiejew2023behrend} for a counterexample to the constancy of the Behrend function of $\Hilb_{\mathbb{A}^3/\kk}^d$ (for $d \geq 24$) and the relationship with the question whether it is generically reduced.

The aim of this paper is to provide a family of examples of non-reduced points on the Hilbert scheme with both the dimension of an affine space and the number of points much smaller than in \cite{Hilbpath} and without using computer calculations. The algebras in the constructed family are very compressed (see e.g. \cite{Jel24a} for a definition and related problems), and this type of algebras are sources of many other pathologies of Hilbert schemes, see \cite{Iar72b, EI78}. The family that we point out shows that the Hilbert scheme of $d$ points on $\Affine^n$ is non-reduced as soon as $d \geq 13, n \geq 6$, which immediately follows from non-reducedness of $\mathrm{Hilb}^{13}_{\Affine^6/\kk}$. In particular we give:
\begin{example}[Corollary \ref{cor:important}]\label{ex:first}
    Let $I$ be the ideal generated by polynomials $f$ in variables $\frac{\partial}{\partial x_1}, \dots, \frac{\partial}{\partial x_6}$ such that $f$ applied to $F = x_1 x_2 x_4 - x_1 x_5^2 + x_2 x_3^2 + x_3 x_5 x_6 + x_4 x_6^2$ is zero or $f$ is of degree greater or equal three. Then $[I] \in \mathrm{Hilb}^{13}_{\Affine^6/\kk}$ is a non-reduced point.
\end{example}

In order to prove non-reducedness of points similar to the one in Example \ref{ex:first} we use the $\G_m$-action on $\cH := \mathrm{Hilb}^{d}_{\Affine^n/\kk}$ coming from the action of $\G_m$ on $\Affine^n$ by multiplication. We define the barycenter scheme $\cB \subseteq \cH$ which consists of points represented by tuples of points in $\Affine^n$ such that their weighted average is $0$, see Definition \ref{def:Bar}. Next we note that $\cH \simeq \Affine^n \times \cB$ and that the $\G_m$-action on $\cH$ restricts to $\cB$. The point from Example \ref{ex:first} is an example of a hedgehog point.
\begin{definition}[Definition \ref{def:hedgehog:precise}]\label{def:hedgehog}
    A $\G_m$-fixed point $[I] \in \cB$ is a \emph{hedgehog point} if it satisfies the following conditions:
    \begin{itemize}
        \item $(T_{[I]} \cB)_{>0} = 0$,
        \item $T_{[I]} \cB \neq (T_{[I]} \cB)_0$,
        \item there is no closed point $[I'] \neq [I]$ in $\cB$ such that $\lim_{t \to 0} t \cdot [I'] = [I]$.
    \end{itemize}
\end{definition}
The last item can be stated in other terms. Let $V$ be the subscheme of $\cB$ consisting of points for which their $\G_m$-limit at zero is $[I]$. We call it the \emph{negative spike at $[I]$}. The last condition is equivalent to $0$-dimensionality of $V$ or to the fact that $V$ is topologically a point. This justifies the name "hedgehog point", as $\cB$ looks like $\cB^{\G_m}$ with infinitesimal spikes, near $[I]$. Equivalent definition of the negative spike and being a hedgehog point is given in Definition \ref{def:hedgehog:precise}.
The following is the first main theorem of this paper.
\begin{theorem}[Hedgehog point Theorem \ref{thm:Hedgehog_point_theorem:precise}]\label{thm:Hedgehog_point_theorem}
    A hedgehog point $[I] \in \cB$ is non-reduced. Hence $(0, [I]) \in \Affine^n \times \cB \simeq \cH$ is non-reduced as well.
\end{theorem}
The geometric idea behind the proof is as follows:
\begin{enumerate}
    \item By a generalisation of \BBname{} decomposition for the $\G_m$-action on the Barycenter scheme we get that there is an open neighbourhood of $[I]$ such that the barycenter looks like an affine $\Affine^1$-cone over $\G_m$-fixed points. Here the fact that $(T_{[I]} \cB)_{>0} = 0$ is crucial.
    \item The fiber over $[I]$ of this cone is $V$ and it is $0$-dimensional as $[I]$ is a hedgehog point. Assuming (for a contradiction) reducedness of $[I]$ in the barycenter scheme, we get that the aforementioned affine cone is trivial i.e. the barycenter scheme is isomorphic to the fix-points subscheme near $[I]$. This follows from semi-continuity of dimension of fibers for projective morphisms.
    \item The tangent space to the fixed-point subscheme of $\cB$ is $(T_{[I]} \cB)_0$ so it is different than $T_{[I]} \cB$, because $[I]$ is a hedgehog point. Hence, we get a contradiction with $[I]$ being reduced.
\end{enumerate}
The key tool is the generalisation of \BBname{} decomposition from \cite{JS19} together with analysis of tangent spaces at the ambient point and of the fiber $V$.

To construct hedgehog points we start with a cubic in six variables $F \in \kk[x_1, \dots, x_6]$ and proceed as in Example \ref{ex:first}. We produce an ideal $I$ which is defined as $\Fperp + S_{\geq 3}$ where $\Fperp \subseteq S = \kk[\frac{\partial}{\partial x_1}, \dots, \frac{\partial}{\partial x_6}]$ is the set of differential polynomials, such that when applied to $F$ they give $0$. The cubic $F$ is \emph{general enough} if it lies in the open (and non-empty) set of cubics from Definition \ref{def:nice}. For those cubics $I = \left\langle \frac{\partial F}{\partial x_1} , \dots, \frac{\partial F}{\partial x_6} \right\rangle^{\perp}$.
We prove the following.
\begin{theorem}[Corollary \ref{cor:Iishedgehog}]\label{thm:cubics_give_hedgehog}
    For general enough cubic $F$ the given point $[I] \in \mathrm{Hilb}^{13}_{\Affine^6/\kk}$ lies in the barycenter subscheme and is a hedgehog point. In particular it is non-reduced.
\end{theorem}
Non-emptiness of the set of general enough cubics is witnessed by $F$ from Example \ref{ex:first}. Due to the paper \cite{Shafarevich} algebras with their Hilbert function $(1, 6, 6)$ form an irreducible component of the barycenter scheme. The point $[I]$ lies only in that component, what follows from the fact that $\cB$ and $\cB^{\G_m}$ are homeomorphic near $[I]$ (by point (2) in the sketch of the proof of Theorem \ref{thm:Hedgehog_point_theorem:precise}, but without assuming reducedness) and from the fact that $\cB^{\G_m}$ consists of algebras of type $(1,6,6)$ near $[I]$. Moreover, a general point of the component containing $[I]$ is smooth, so non-reducedness here is of the "embedded component" type.

After proving theorems \eqref{thm:Hedgehog_point_theorem} and \eqref{thm:cubics_give_hedgehog} we determine the negative spike $V$ over $[I]$ coming from a general enough cubic $F$. We prove the following.
\begin{theorem}[Theorem \ref{thm:completeV}]\label{thm:intro_completeV}
    If $I = \Fperp + S_{\geq 3}$, for a general enough cubic $F$, then the negative spike at $[I]$ is isomorphic to $\Spec(S/\Fperp)$.
\end{theorem}
As a subscheme of the Hilbert scheme $V$ is induced by a deformation. We describe a related deformation and call it the \emph{fractal family}. We chose the name fractal because the fiber over the unique closed point and the base space of this deformation are almost the same, so it is like a deformation over itself. This exhibits an interesting phenomenon: the Hilbert scheme which parametrises finite schemes, naturally contains some of them (or very close to such), as closed subschemes. Overall, the isomorphism between $\Spec(S/\Fperp)$ and $V$ comes from the composition of the morphism to $\cH$ defined by the fractal family and the projection $\cH \simeq \Affine^n \times \cB \to \cB$ from the Hilbert scheme to its barycenter subscheme.
Recently Oszer \cite{PiotrOszer2026} generalised our notion of fractal families. In particular, he shows that the Gorenstein locus of the Hilbert scheme of points is
non-reduced.

\subsection{Organization}

We briefly discuss the contents of the paper. In Section \ref{sec:prelims} we introduce Hilbert schemes of points, describe a few groups acting on Hilbert schemes of points on affine spaces and define \BBname{} decomposition. Next we introduce obstruction theories, describe the barycenter scheme and present a simple form of Macaulay's duality. In Section \ref{sec:main} we start with the notion of being general enough for cubics in six variables and the details of the procedure of making a point $[I]$ in the $\mathrm{Hilb}^{13}_{\Affine^6/\kk}$ from a general enough cubic $F$. After that we determine the tangent space to $V$ at $[I]$ and calculate primary obstructions for $V$. It allows us to prove the $0$-dimensionality of $V$ and conduct the main argument of the paper. At the end we determine $V$ completely by constructing the fractal family.

\section{Preliminaries}\label{sec:prelims}
    We assume familiarity with schemes and their functorial description (e.g. as in \cite[Part 1, Chapter 2]{FGI+05}). We work over an algebraically closed field $\kk$ of characteristic $0$, although algebraic closedness plays a significant role only at a few places. Characteristic zero assumption matters for the apolarity theory.

    Now we list some of the notations appearing later in the paper. The reason for two notations - $S$ and $P$ - denoting the polynomial algebra will become clear after Definition \ref{def:duality}.
    \begin{itemize}
        \item $\delta_{ij}$ denotes the Kronecker delta.
        \item By $\langle a_1, \dots, a_n \rangle$ we denote the $\kk$ vector space spanned by vectors $a_1, \dots, a_n$, if $a_i$'s are elements of a $\kk$-vector space clear from the context.
        \item $X$ denotes a scheme locally of finite type over $\kk$.
        \item For a scheme $X$ and a scheme $T$ over $\kk$, we denote by $X_T$ the base change of $X$ to $T$.
        \item For $x \in X(\kk)$ the tangent space to $X$ at $x$ will be denoted $T_{x} X$. Moreover, $T_{x}^{\vee} X$ will denote the cotangent space to $X$ at $x$. 
        \item $T$ denotes a test scheme for a functor.
        \item $S$ denotes a polynomial algebra $\kk[\alpha_1, \dots, \alpha_n]$ where $n$ is usually clear from the context (in which case we write $\overline{\alpha}$ for the tuple $(\alpha_1, \dots, \alpha_n)$).
        \item $P$ denotes a polynomial algebra $\kk[x_1, \dots, x_n]$ (sometimes written $\kk[\overline{x}]$) with different names for variables.
        \item If $A$ is a $\kk$-algebra, then  $S_A$ denotes $S \otimes_{\kk} A = A[\alpha_1, \dots, \alpha_n]$.
        \item If $A$ is a local algebra, $\mm_A$ denotes its maximal ideal.
        \item If $n$ is fixed, we will write $\xx$ or $\aalpha$ for respectively $x_1, \dots, x_n$ or $\alpha_1, \dots, \alpha_n$.
        \item $\kk[\eps]$ denotes the ring $\frac{\kk[\eps]}{(\eps)^2}$, so whenever something is denoted by $\eps$ (or $\eps'$) we have $\eps^2 = 0$ (although, sometimes we write $\frac{\kk[\eps]}{(\eps)^2}$ or $\frac{\kk[\eps,\eps']}{(\eps, \eps')^2}$ when we compare two infinitesimals $\eps, \eps'$).
        \item We write $S[\eps]$ for $S_{\kk[\eps]}$.
        \item $\Affine^n$/$\Pn$, $\G_m$ denote respectively the affine/projective space of dimension $n$ and the algebraic torus $\Spec(\kk[t^{\pm 1}])$.
        \item If $\G_m$ acts on a vector space $V$, then the induced grading on $V$ is given by $V_n := \{ v \in V : t \in \G_m(k) \textnormal{ acts by multiplying by } t^{-n} \}$.
        \item If $B$ is a non-negatively graded $\kk$-algebra then its \emph{Hilbert function} is the function $\mathbb{N} \to \mathbb{N} \cup \{\infty\}$ defined by $n \mapsto \dim_\kk(B_n)$. 
        \item For a linear space $W$ let $\Sym_2(W)$ be a subspace of $W \otimes W$ consisting of symmetric tensors. The canonical symmetric bilinear $\kk$-vector space homomorphism from $W \times W$ to $\Sym_2(W)$ is defined by the formula
        \begin{equation*}
            (w, v) \mapsto \frac{1}{2}(w \otimes v + v \otimes w).
        \end{equation*}
        This yields a natural duality between $\Sym_2(W)$ and $\Sym_2(W^{\vee})$ induced by the natural duality between $W \otimes W$ and $W^{\vee} \otimes W^{\vee}$. More precisely we define a map
        \begin{equation}\label{eq:dualityw}
            \cdot : \Sym_2(W^{\vee}) \times \Sym_2(W) \to \kk
        \end{equation}
        by the formula $(\eta \mu) \cdot (w v) = (\frac{1}{2}(\eta \otimes \mu + \mu \otimes \eta)) \cdot (\frac{1}{2}(w \otimes v + v \otimes w)) := \frac{1}{4}(\eta(w) \mu(v) + \eta(v) \mu(w) + \mu(w) \eta(v) + \mu(v) \eta(w)) = \frac{1}{2}(\eta(w) \mu(v) + \eta(v) \mu(w))$.
        This gives an identification of $\Sym_2(W)^{\vee}$ with $\Sym_2(W^{\vee})$.
    \end{itemize}
    
    \subsection{Hilbert schemes}\label{sec:Hilbschm}
    We recall the definition of the Hilbert functor of points on a scheme. 
    \begin{definition}\label{def:Hilbert:first}
        Let $X$ be a scheme over $\kk$. Then the functor $\mathrm{Hilb}^ d_{X/\kk} : (\mathit{Sch}/\kk)^{op} \longrightarrow \Set$ given by
        \begin{equation*}
            \mathrm{Hilb}^ d_{X/\kk}(T) = \left\{  \begin{matrix}  Z \subseteq X_ T\text{ closed subscheme such that }
            \\  Z \to T\text{ is flat and finite of degree }d 
            \end{matrix} \right\}
        \end{equation*}
        is called a Hilbert functor of $d$ points on $X$.
    \end{definition}
    \begin{example}\label{ex:functorH}
        If $X = \Affine^n$ then for a $\kk$-algebra $A$ we have
        \begin{equation*}
            \mathrm{Hilb}^ d_{\Affine^n/\kk}(\Spec(A)) = \left\{  \begin{matrix}  I \subseteq S_A \text{ ideal such that }
            \\  S_A/I \text{ is a flat $A$-module of rank }d 
            \end{matrix} \right\}.
        \end{equation*}
        In this section we fix $n$ and $d$ and write $\cH$ for both this functor and the scheme representing it. If we want to emphasize the number of points $d$ we write $\cH^d$. Also, we will denote by $[I]$ an $A$-point corresponding to an ideal $I \subseteq S_A$ such that $S_A/I$ is flat of rank $d$. Note that if $[I] \in \cH^d(\Spec(A))$ and $f:A \to B$ is a homomorphism, then $J := I \otimes_A B \subseteq S_B$ induces a point $[J] \in \cH^d(\Spec(B))$ and $f$ induces a morphism of short exact sequences
        \[\begin{tikzcd}
    	0 & {I} & {S_A} & {\frac{S_A}{I}} & 0 \\
       	0 & J & {S_B} & {\frac{S_B}{J}} & 0.
       	\arrow[from=1-2, to=1-3]
       	\arrow[from=1-3, to=1-4]
       	\arrow[from=1-4, to=1-5]
       	\arrow[from=1-1, to=1-2]
    	\arrow[from=2-1, to=2-2]
    	\arrow[from=2-2, to=2-3]
    	\arrow[from=2-3, to=2-4]
    	\arrow[from=2-4, to=2-5]
    	\arrow[from=1-4, to=2-4]
    	\arrow[from=1-3, to=2-3]
    	\arrow[from=1-2, to=2-2]
    \end{tikzcd}\]
    \end{example}
    
    In the whole paper we will be only concerned with the Hilbert schemes of points on affine spaces. To study local properties of Hilbert functors/schemes (which is not restrictive as every smooth scheme is étale locally around each point an affine space) the following description of the tangent space will be crucial.
    
    \begin{theorem}\label{thm:tangenttohilb}
        Fix a $\kk$-point $[I] \in \cH$. Then the tangent space $T_{[I]} \cH$ is isomorphic to the $\kk$-vector space $\Hom_S(I, S/I)$.
    \end{theorem}
    
    The full proof of this theorem can be found in \cite[Corollary 6.4.10]{FGI+05}. We will use this identification a few times, so we recall how to obtain a homomorphism in $\Hom_S(I, S/I)$ from an element of $T_{[I]} \cH$.
    First, using the fact that $T_{[I]} \cH = \{ \varphi \in \Mor \left( \Spec(\kk[\eps]), \cH \right) \ | \  \varphi((\eps)) =  [I] \}$, it follows by the identification from Example \ref{ex:functorH}, that an element of $T_{[I]} \cH$ is given by a diagram
    \[\begin{tikzcd}
    	0 & {I'} & {S[\eps]} & {\frac{S[\eps]}{I'}} & 0 \\
       	0 & I & {S} & {\frac{S}{I}} & 0,
       	\arrow[from=1-2, to=1-3]
       	\arrow["\beta", from=1-3, to=1-4]
       	\arrow[from=1-4, to=1-5]
       	\arrow[from=1-1, to=1-2]
    	\arrow[from=2-1, to=2-2]
    	\arrow["\alpha", from=2-2, to=2-3]
    	\arrow[from=2-3, to=2-4]
    	\arrow[from=2-4, to=2-5]
    	\arrow["\gamma"', from=1-4, to=2-4]
    	\arrow[from=1-3, to=2-3]
    	\arrow[from=1-2, to=2-2]
    \end{tikzcd}\]
    where $\gamma$ is the induced map between quotients. The corresponding element of $\Hom_S(I, S/I)$ is given by $\delta$ on the following diagram
    \[\begin{tikzcd}
        &&& 0 \\
    	&&& {\frac{S}{I}} \\
    	0 & {I'} & {S[\eps]} & {\frac{S[\eps]}{I'}} & 0 \\
    	0 & I & {S} & {\frac{S}{I}} & 0. \\
    	&&& 0
    	\arrow[from=3-2, to=3-3]
    	\arrow["\beta", from=3-3, to=3-4]
    	\arrow[from=3-4, to=3-5]
    	\arrow[from=3-1, to=3-2]
    	\arrow[from=4-1, to=4-2]
    	\arrow["\alpha", from=4-2, to=4-3]
    	\arrow[from=4-3, to=4-4]
    	\arrow[from=4-4, to=4-5]
    	\arrow["\gamma"', from=3-4, to=4-4]
    	\arrow[from=3-3, to=4-3]
    	\arrow[from=3-2, to=4-2]
    	\arrow[from=2-4, to=3-4]
    	\arrow[from=1-4, to=2-4]
    	\arrow[from=4-4, to=5-4]
    	\arrow["s", bend left, from=4-3, to=3-3]
    	\arrow["\delta", bend left, crossing over, dashed, from=4-2, to=2-4]
    \end{tikzcd}\]
    Here $s$ is the natural embedding of $S$ into $S[\eps]$ and the vertical short exact sequence comes from the short exact sequence of $\kk[\eps]$-modules
    \[ \begin{tikzcd}
        0 & \kk \simeq (\eps) & \kk[\eps] & \kk & 0
        \arrow[from=1-1, to=1-2]
        \arrow[from=1-2, to=1-3]
        \arrow[from=1-3, to=1-4]
        \arrow[from=1-4, to=1-5]
    \end{tikzcd} \]
    tensored with $S[\eps]/I'$ over $k[\eps]$. The homomorphism $\beta \circ s \circ \alpha$ composed with $\gamma$ is $0$ (by diagram chasing) thus we get the existence of $\delta$. 
    \begin{remark}\label{rem:tangentinverse}
        On the other hand, if we are given $\delta \in \Hom_S(I, S/I)$, then the corresponding ideal $I' \subseteq S[\eps]$ is given by
        \begin{equation*}
            I' = ( f - \eps \cdot \delta(f)' : f \in I ),
        \end{equation*}
        where $\delta(f)' \in S$ is any representative of $\delta(f) \in S/I$. A precise explanation can be found in \cite{FGI+05}.
    \end{remark}
    
    \subsection{Groups acting on $\cH$}\label{subsection:groupsactiongonH}
    
    If an algebraic commutative group $G$ acts on a scheme $X$, then there is an induced action of $G$ on $\mathrm{Hilb}^ d_{X/\kk}$. Indeed, take $T$ - a test scheme. Then the $G$-action $G \times X \to X$ yields (by the base change to $T$) a morphism $G_T \times X_T \to X_T$. If we take $t \in G(T) = \Mor(T, G) = \Mor_T(T, G_T)$, then we get the diagram
    \[\begin{tikzcd}
    	{G_T \times X_T} & {X_T}, \\
    	{X_T}
    	\arrow[from=1-1, to=1-2]
    	\arrow["{(t, id)}", from=2-1, to=1-1]
    	\arrow["{t \cdot}"', from=2-1, to=1-2]
    \end{tikzcd}\]
    where the vertical map $X_T \to G_T \times X_T$ is the identity on the second coordinate and factors through the $T$-point $t$ on the first one. Now if $Z \subseteq X_T$ represents a $T$-point of $\mathrm{Hilb}^ d_{X/\kk}$, then the action of $t$ on $Z$ is defined by the image of map $(t\cdot) : X_T \to X_T$. This map is an isomorphism (because $G$ is an algebraic group), so $t \cdot Z := (t \cdot)(Z) \subseteq X_T \to T$ is again locally finite of degree $d$.
    We provide two examples of lifting actions from $\Affine^n$ to $\cH$.
    
    \begin{example}\label{ex:torusaction}
        Consider the $\G_m$-action on $\Affine^n$ given by the standard grading $\deg(\alpha_i)=1$. Let us describe the induced morphism $\G_m \times \cH \to \cH$ constructed on its $T$-points for $T = \Spec(A)$ where $A$ is a $\kk$-algebra. We are interested in the map
        \begin{equation*}
            \G_m(\Spec(A)) \times \cH(\Spec(A)) \to \cH(\Spec(A)).
        \end{equation*}
        Take $a \in \G_m(\Spec(A)) = A^*$ (invertible elements). The element $a$ yields an isomorphism
        \[\begin{tikzcd}
        	{\Affine_A^n} & {\Affine_A^n},
        	\arrow["{a \cdot}"', from=1-1, to=1-2]
        \end{tikzcd}\]
        which on the level of algebras comes from the $A$-algebra homomorphism $(a \cdot)^{\#}:S_A = A[\aalpha] \to A[\aalpha]$ sending $\alpha_i$ to $a \alpha_i$. Now, a point $[Z] \in \cH(\Spec(A))$ is represented by an ideal $I \subseteq S_A$ and the image point is given by the ideal $a \cdot I := \ker(q \circ (a \cdot)^{\#})$ where $q: S_A \to \frac{S_A}{I}$ is the quotient map. The justification of this formula follows from the next two diagrams (the left hand side is the definition of the action and the right hand side is the corresponding diagram on the level of algebras)
        \[\begin{tikzcd}
        	{\Affine_A^n} & {\Affine_A^n} & {S_A} & {S_A} \\
        	Z & {a \cdot Z} & {\frac{S_A}{I}} & {\frac{S_A}{a \cdot I}}.
        	\arrow["{a \cdot}"', from=1-1, to=1-2]
        	\arrow[hook', from=2-1, to=1-1]
        	\arrow[hook', from=2-2, to=1-2]
        	\arrow[from=2-1, to=2-2]
        	\arrow["{(a \cdot)^{\#}}", from=1-4, to=1-3]
        	\arrow["q", two heads, from=1-3, to=2-3]
        	\arrow[two heads, from=1-4, to=2-4]
        	\arrow[from=2-4, to=2-3]
        \end{tikzcd}\]
    \end{example}
    Using the description of $\G_m$-action on $\cH$ one can prove the following result.
    \begin{proposition}\label{prop:torusequivariance}
        Suppose that $A$ is a $\ZZ$-graded $\kk$-algebra and consider a morphism $\phi : \Spec(A) \to \cH$ given by the subscheme $Z \subseteq \Affine_A^n$ cut out by the ideal $I \subseteq S_A$. Consider the $\ZZ$-grading on $S_A = S \otimes_\kk A$ given by $(S_A)_l = \bigoplus_{i+j = l} S_i \otimes_\kk A_j$. If $I$ is homogeneous with respect to this grading, then the morphism $\phi$ is $\G_m$-equivariant.
    \end{proposition}
    \begin{proof}
        We need to check that the following diagram commutes:
        \[\begin{tikzcd}
        	{\Spec(A)} & \cH \\
        	{\G_m \times \Spec(A)} & {\G_m \times \cH}.
        	\arrow["\mu"', from=2-1, to=1-1]
        	\arrow["\phi"', from=1-1, to=1-2]
        	\arrow["\nu", from=2-2, to=1-2]
        	\arrow["{id \times \phi}", from=2-1, to=2-2]
        \end{tikzcd}\]
        Let $\mu' = \phi \circ \mu$ and $\nu' = \nu \circ (id \times \phi)$. These both are morphisms from $\G_m \times \Spec(A) = \Spec(A[t^{\pm 1}])$ to $\cH$ so they correspond uniquely to some ideals in $S_{A[t^{\pm 1}]}$. It is enough to check that these ideals are the same. Consider a homogeneous polynomial $f = \sum_{N} a_{N} \aalpha^{N} \in I$, so that $\deg(a_{N}) + |{N}| = l$ for all $N$ and some $l \in \ZZ$.
        
        We start with the calculation of the ideal corresponding to the map $\mu'$. Since $\mu' = \phi \circ \mu$ the subscheme $Z'$ such that $\mu'$ is given by $[Z'] \in \cH(\G_m \times \Spec(A))$ is a pullback of $Z \subseteq \Affine_A^n$ via the map $(id \times \mu) : \Affine^n \times (\G_m \times \Spec(A)) \to \Affine^n \times \Spec(A)$. By writing $\G_m \times \Spec(A) = \Spec(A[t^{\pm 1}])$ and looking at algebras we get:
        \[\begin{tikzcd}
        	{\Affine_{A[t^{\pm 1}]}^n} & {\Affine_A^n} & {S_{A[t^{\pm 1}]}} & {S_A} \\
        	{Z'} & Z & {\frac{S_{A[t^{\pm 1}]}}{I'}} & {\frac{S_A}{I}}.
        	\arrow[from=1-1, to=1-2]
        	\arrow[hook', from=2-1, to=1-1]
        	\arrow[hook', from=2-2, to=1-2]
        	\arrow[from=2-1, to=2-2]
        	\arrow[from=1-4, to=1-3]
        	\arrow[two heads, from=1-3, to=2-3]
        	\arrow[two heads, from=1-4, to=2-4]
        	\arrow[from=2-4, to=2-3]
            \arrow["\lrcorner"{anchor=center, pos=0.125, rotate=90}, draw=none, from=2-1, to=1-2]
        	\arrow["\lrcorner"{anchor=center, pos=0.125, rotate=90}, draw=none, from=2-3, to=1-4]
        \end{tikzcd}\]
        Since $Z'$ is a pullback we have $\frac{S_{A[t^{\pm 1}]}}{I'} = S_{A[t^{\pm 1}]} \otimes_{S_A} \frac{S_A}{I}$. By right exactness of the tensor product the ideal $I'$ is generated by the image of the map \[(id \otimes i):S_{A[t^{\pm 1}]} \otimes_{S_A} I \to S_{A[t^{\pm 1}]} \otimes_{S_A} S_A = S_{A[t^{\pm 1}]},\]
        where $i:I \to S_A$ is the inclusion. We calculate the image of $1 \otimes f$ via $id \otimes i$:
        \begin{equation*}
            1 \otimes f = 1 \otimes \sum_{N} a_{N} \aalpha^{N} = \sum_{N} \aalpha^N \otimes a_{N} = \sum_{N} \aalpha^N \cdot a_{N} \cdot t^{\deg(a_{N})} \otimes 1 \mapsto \sum_{N} \aalpha^N \cdot a_{N} \cdot t^{\deg(a_{N})} \in S_{A[t^{\pm 1}]}.
        \end{equation*}
        Since the calculation works for any $f$ in $I$, we get that $I'$ is generated by polynomials of the form $f^* = \sum_{N} a_{N} t^{\deg(a_N)} \cdot \aalpha^{N}$, for $f = \sum_{N} a_{N} \aalpha^{N}$ varying over homogeneous polynomials in $I$.
        
        Now we calculate the ideal corresponding to the map $\nu'$. Again, by the fact that $\G_m \times \Spec(A) = \Spec(A[t^{\pm 1}])$ we can see $id \times \phi$ from the diagram at the beginning of this proof as an element of $\G_m (\Spec(A[t^{\pm 1}])) \times \cH(\Spec(A[t^{\pm 1}]))$. By the definition of $id \times \phi$ this element is the pair $(t, [Z''])$, where $Z''$ is a pullback of $Z$ as on the following diagram
        \[\begin{tikzcd}
        	{\Affine_{A[t^{\pm 1}]}^n} & {\Affine_A^n} \\
        	{Z''} & Z.
        	\arrow[hook', from=2-2, to=1-2]
        	\arrow["\pi"', from=1-1, to=1-2]
        	\arrow[hook', from=2-1, to=1-1]
        	\arrow[from=2-1, to=2-2]
        	\arrow["\lrcorner"{anchor=center, pos=0.125, rotate=90}, draw=none, from=2-1, to=1-2]
        \end{tikzcd}\]
        Here $\pi$ comes from the projection $\G_m \times \Spec(A) \to \Spec(A)$. Hence the ideal of $Z''$ in $S_{A[t^{\pm 1}]}$ is generated by $I$ and we denote it by $\tilde{I}$.
        By Example \ref{ex:torusaction} the image of the pair $(t, [Z''])$ by $\nu$ is given by the subscheme $t \cdot Z'' \subseteq \Affine_{A[t^{\pm 1}]}^n$. The ideal $I''$ of $t \cdot Z''$ is given by the kernel of the composition
        \[\begin{tikzcd}
        	{S_{A[t^{\pm 1}]}} & {S_{A[t^{\pm 1}]}} & {\frac{S_{A[t^{\pm 1}]}}{\tilde{I}}},
        	\arrow["{(t \cdot)^{\#}}", from=1-1, to=1-2]
        	\arrow[from=1-2, to=1-3]
        \end{tikzcd}\]
        where $(t \cdot)^{\#} : S_{A[t^{\pm 1}]} \to S_{A[t^{\pm 1}]}$ is a $A[t^{\pm 1}]$-algebras homomorphism sending $\alpha_i$ to $t\alpha_i$. Thus $I''$ is generated by elements of the form $f(t^{-1} \cdot \aalpha) = \sum_{N} a_{N} \aalpha^{N} \cdot t^{-|N|}$, where $f = \sum_{N} a_{N} \aalpha^{N}$ varies over homogeneous polynomials in $I$.
        
        The element $t \in A[t^{\pm 1}]$ is invertible. So, for any homogeneous $f = \sum_{N} a_{N} \aalpha^{N} \in I$ we have
        \begin{equation*}
            f^* = \sum_{N} a_{N} t^{\deg(a_N)} \cdot \aalpha^{N} = t^l \cdot \sum_{N} a_{N} t^{\deg(a_N) - l} \cdot \aalpha^{N} = t^l \cdot \sum_{N} a_{N} t^{-|N|} \cdot \aalpha^{N} = t^l \cdot f(t^{-1} \cdot \aalpha).
        \end{equation*}
        Hence, $I' = I''$ and we are done.
    \end{proof}
    
    \begin{example}\label{ex:additiveaction}
        Now we analyse the induced $\G_a^n$-action on $\cH$ coming from the action on $\Affine^n$ by translations. Again, for a $\kk$-algebra $A$ we get a map
        \begin{equation*}
            \G_a^n (\Spec(A)) \times \cH(\Spec(A)) \to \cH(\Spec(A)),
        \end{equation*}
        defined by $((a_1, \dots, a_n), I) \mapsto \ker(q \circ (\Bar{a}+)^{\#})$ where $(\Bar{a}+)^{\#}:S_A \to S_A$ is a $A$-homomorphism sending $\alpha_i$ to $\alpha_i + a_i$.
    \end{example} 
    
    From Example \ref{ex:torusaction} we can see that a $\kk$-point $[I] \in \cH$ is a $\G_m$-fixed point if and only if $I$ is homogeneous with respect to the standard grading on $S$. Thus, if $I$ is homogeneous, then $\G_m$ acts on its tangent space $T_{[I]} \cH$. This action induces a grading $\Hom_S(I, S/I) = \bigoplus_{k \in \ZZ} \Hom_S(I, S/I)_k$, where
    \begin{equation*}
        \Hom_S(I, S/I)_k = \{ \delta \in \Hom_S(I, S/I) : \delta(I_n) \subseteq (S/I)_{n+k} \textnormal{ for all integers $n$} \}.
    \end{equation*}
    
    There is no analogue of this decomposition for the $\G_a^n$ additive action on $\cH$ (as in Example \ref{ex:additiveaction}), because there is no $\G_a^n$-invariant point in $\cH$. However, one can still ask what happens locally or how does the orbit map (of a point) look like. More precisely, let $[I]$ be a $\kk$-point of $\cH$ and let $\G_a^n \to \cH$ be its orbit map under the additive action of $\G_a^n$. Then the induced differential $\eta : T_{ \{0\} } \G_a^n \to T_{[I]} \cH$ sends the $i$'th basis vector $\partial_i$ to $\eta(\partial_i) \in T_{[I]} \cH = \Hom_S(I, S/I)$ given by the formula
    \begin{equation*}
        \eta(\partial_i) (f) = \frac{\partial f}{\partial \alpha_i}.
    \end{equation*}
    Thus, by abuse of notation we will write $\partial_i$ for $\eta(\partial_i)$.
    
    Suppose additionally that $[I]$ is a $\G_m$-fixed point of $\cH$. From the description of the grading on $T_{[I]} \cH$ we see that the image of $\eta$ lies in $\Hom_S(I, S/I)_{-1} \subseteq \Hom_S(I, S/I)_{< 0}$. Thus it is natural to consider the following definition.
    
    \begin{definition}\label{def:TNT}
        A $\G_m$-fixed $\kk$-point $[I] \in \cH$ has \emph{trivial negative tangents} (abbreviated TNT), if the map $\eta : \langle \partial_1, \dots, \partial_n \rangle \to \Hom_S(I, S/I)_{< 0}$ from Section \ref{sec:Hilbschm} is surjective.
    \end{definition}
    
    This means that the only negative tangents at $[I]$ come from $\Affine^n$ additive action and are given by derivatives.
    
    This definition is taken from \cite{elemcomp}. Note that because the image of $\eta$ is contained in $\Hom_S(I, S/I)_{-1}$ and $\partial_1, \dots, \partial_n \in \Hom_S(I, S/I)_{-1}$ are linearly independent it follows that $[I]$ as in the definition has TNT if and only if $\dim_\kk \Hom_S(I, S/I)_{<0} = n$.
    
    \subsection{\BBname{} decomposition}\label{subsection:BBdecomp}
    A torus action can be useful to study geometry of Hilbert schemes. To exploit it we will need the theory of \BBname{} decompositions which originates from the papers \cite{ABB1} and \cite{ABB2}. Later in the paper, we will work with the \BBname{} decomposition induced by the $\G_m$-action from Example \ref{ex:torusaction}.
    
    In this subsection we present necessary tools from the paper \cite{JS19} which generalises \BBname's result to non-necessarily smooth schemes. In the paper the authors work with linearly reductive groups, so the class that contains $\G_m$, but we will present the results from there only for the $\G_m$-case.
    
    \begin{definition}\label{def:xplusfunctorial}
        Let $X$ be a $\G_m$-scheme. Consider the functor $\Xplus : (\mathit{Sch}/\kk)^{op} \longrightarrow \Set$ defined by
        \begin{equation*}\label{eq:definitionFunctorial:new}
        \Xplus(T) = \left\{ \symbf \colon \mathbb{A}^1 \times T \to \varX\ |\
        \symbf \mbox{ is $\Gmult$-equivariant} \right\},
        \end{equation*}
        where $\G_m$ acts on $\Affine^1 \times T$ by the natural action of $\G_m$ on $\Affine^1$ (i.e. by multiplication). It is called the \emph{\BBname \ decomposition} of $X$. Moreover, \BBname{} decomposition is functorial in $X$, by what we mean that if $X, Y$ are $\G_m$-schemes and $g : X \to Y$ is $\G_m$-equivariant, then there is $g^+ : X^+ \to Y^+$ such that $g^+(T) : X^+(T) \to Y^+(T)$ is given by $g^+(T)(f : \Affine^1 \times T \to X) = (g \circ f : \Affine^1 \times T \to Y) \in Y^+(T)$ for any test scheme $T$. 
    \end{definition}
    
    From the \BBname \ decomposition of a scheme $X$ we get three natural transformations:
    \begin{itemize}
        \item $\Xplus(T) \to X(T)$ given by restricting $\symbf$ to $1 \times T$. This is called \emph{forgetting about the limit} and will be denoted by $i_X$.
        \item $\Xplus(T) \to X^{\G_m}(T)$ given by restricting $\symbf$ to $0 \times T$. This is called \emph{restricting to the limit} and will be denoted by $\pi_X$.
        \item $X^{\G_m}(T) \to \Xplus(T)$ given by $(\symbf : T \to X^{\G_m}) \mapsto ((j \circ \symbf \circ pr_1): T \times \Affine^1 \to X)$, where $j : X^{\G_m} \to X$ is the closed embedding. This morphism is a section of $\pi_X$ and will be denoted by $s_X$.
    \end{itemize}
    
    Under this notation the following holds.
    
    \begin{theorem}[Existence of \BBname{} decompositions~{\cite[Theorem~1.1]{JS19}} ]\label{ref:introRepresentability:thm}
        Let $\varX$ be a $\G_m$-scheme locally of finite type over
        $\kk$. Then the functor $\Xplus$ defined
        by \eqref{eq:definitionFunctorial:new} is represented by a scheme locally
        of finite type over $\kk$, denoted also by $\Xplus$.
        Moreover, the scheme $\Xplus$ has a natural $\Affine^1$-action. The map
        $\iinftyX\colon \Xplus\to \varX^{\G_m}$ is \emph{affine} of finite type and
        equivariant.
    \end{theorem}
    
    Moreover, if $X$ is separated, then $X^+$ can be seen as a subset of $X$, because we have the following fact.
    \begin{fact} \cite[lemma 5.8]{JS19}\label{fact:mono}
        For separated $X$ the morphism $i_X : X^+ \to X$ is a monomorphism of functors.
    \end{fact}
    Hence, for any scheme $T$ we will treat $X^+(T)$ as a subset of $X(T)$, as long as $X$ is separated. In particular, for a $\kk$-point $x \in X^{\G_m}$, identifying $x$ with $s(x) \in X^+$ and with $x \in X$ we can treat the tangent space to $X^+$ at $x$ as a subset of $T_{x} X$ \[T_{x} X^+ = X^+( \Spec(\kk[\eps]) )_{\{\textnormal{pt.}\} \mapsto x} \subseteq X( \Spec(\kk[\eps]) )_{\{\textnormal{pt.}\} \mapsto x} = T_{x} X.\]
    
    This fact also justifies why $X^+$ is called a decomposition of $X$. Geometrically speaking, $\kk$-points of $X^+$ are these points $p$ of $X$ for which the $\G_m$-orbit of $p$ can be extended to a map from $\Affine^1$ to $X$. The image of $0$ of this extended map is $\G_m$-invariant and it is equal $\pi_X(p)$. One can think of $\pi_X(p)$ as the limit $\lim_{t \to 0} t \cdot p$ where dot represents the action of $\G_m$ on $X$. For a proper scheme such limits always exist. A good example to bear in mind is the \BBname \ decomposition of $\G_m$-action on $\Pone$.
    
    \begin{example}\label{example:projectiveline}
        Consider the case when $X = \Pone$ and $\G_m$ fixes $\infty$ and acts on $\Pone \setminus \{ \infty \} = \Affine^1$ by multiplication. Then $\Xplus = \Affine^1 \amalg \{ \infty \}$ and $X^{\G_m} = \{ 0, \infty \}$. The morphism $\pi_X$ takes $\Affine^1$ to $0$ and $\infty$ to $\infty$. On the other hand, $i_X$ embeds $\Affine^1$ onto $\Pone \setminus\{ \infty \}$ and sends $\infty$ to $\infty \in \Pone$.
        
        In fact for any smooth $X$ the map $\pi_X$ is locally an affine bundle and this was proved by \BBname{} in \cite{ABB1}, \cite{ABB2}.
    \end{example}
    
    As can be seen in Example \ref{example:projectiveline}, near $0$ the morphism $i_X$ is even an isomorphism. Thus one can study the local geometry of $\Pone$ near $0$ by looking at $(\Pone)^+$. This is not a coincidence as seen in the following proposition which will be crucial for our purposes.
        
    \begin{proposition}\label{ref:comparison:intro:prop} \cite[Proposition 1.6]{JS19}
        Let $\varX$ be separated and locally of finite type.
        Let $x\in \varX^{\G_m}$ be a fixed $\kk$-point such that the
        cotangent space $T_{x}^{\vee} \varX$ has no negative degrees (with respect to the grading induced by the $\G_m$-action on $\varX$). Then the map $\ioneX\colon \Xplus \to \varX$ is an open embedding near $x' = \isectionX(x)\in \Xplus$.
        More precisely, there exists an affine open
        $\Affine^1$-stable neighbourhood
        $U$ of $x'$ such that $(\ioneX)_{|U}\colon U\to \varX$ is an open
        embedding. In particular, $x\in \varX$ has an affine open $\G_m$-stable
        neighbourhood.
        Conversely, if a fixed point $x$ has an affine open $\Affine^1$-stable neighbourhood
        in $\varX$, then $T_{x}^{\vee} \varX$ has no negative degrees.
    \end{proposition}
    
    \begin{remark}
        The condition of $T_{x}^{\vee} \varX$ having no negative degrees is equivalent to $T_{x} \varX$ having no positive degrees, because the $\G_m$-grading on the dual space is reversed.
    \end{remark}
    
    It is worth mentioning here, that the idea of using \BBname \ decomposition in order to study Hilbert schemes is not new. However, before the tools only applied to smooth varieties, so these applications were concerning mostly Hilbert schemes of points on smooth surfaces, which by \cite{Fogarty} are smooth. An example of such an application is presented in the paper \cite{ES87} by Ellingsrud and Strømme where they calculate the homology groups of the complex manifolds corresponding to Hilbert schemes of points on the projective plane over complex numbers. Furthermore, Iarrobino \cite{Iar72a} and Briancon \cite{Bri77} compute \BBname \ decompositions, for $\dim X = 2$, without explicitly mentioning it.

    \subsection{Obstruction theories}

    We will need obstruction theories so we briefly introduce them here. We follow presentation from \cite{JS21}. For a more detailed explanation see \cite{DefHar}.
    
    We start with some intuition. Let $x \in X$ be a $\kk$-point of a scheme $X$ locally of finite type over $\kk$. Then $x$ is an image of a morphism $\Spec(\kk) \to X$ and one could ask how does this morphism extends to spectras of Artinian rings other then $\kk$ itself. It is natural then to consider the functor $D_{X, x} : \Artk \to \Set$ from the category $\Artk$ of local Artinian $\kk$-algebras with the residue field $\kk$ given by
    \begin{equation*}\label{ref:def:eq}
        D_{X,x}(R) = \{ \symba: \Spec(R) \to X \textnormal{ | $\symba(\mm_R)=x$}\}.
    \end{equation*}
    It is called the \emph{deformation functor} of $X$ at $x$ and it is pro-representable by the completion of the local ring $\hat{\cO}_{\varX, x}$ (i.e. $D_{X,x}(R) \simeq \Hom(\hat{\cO}_{\varX, x},R)$). 
    
    Let us reformulate the question about extensions of $\Spec(\kk) \to X$ in a purely algebraic way. If $(A, \mm)$ is a Noetherian complete local $\kk$-algebra (e.g. $\hat{\cO}_{\varX, x}$), when can we lift a morphism $\symba_0 : A \to B_0$ to a morphism $\symba: A \to B$ for a local surjection of Artinian rings $B \to B_0$ as on the following diagram?
    \[\begin{tikzcd}
    	B & {B_0} \\
    	A
    	\arrow[from=1-1, to=1-2]
    	\arrow["{\symba_0}"', from=2-1, to=1-2]
    	\arrow["\symba", dashed, from=2-1, to=1-1]
    \end{tikzcd}\]
    
    For technical reasons we assume that $\mm_B \cdot \ker(B \to B_0) =0$ - this is called a \emph{small extension}. Every local surjection of Artinian rings is a composition of small extensions so by restricting to such morphisms we don't lose any information.

    Let $D=\Hom(A, -)$.  The Cohen structure theorem \cite{Cohen} (see \cite[Tag 0323]{Stacks} for the proof in a general setup) says that $A$ can be presented as a quotient of $\hat{\symbp} := \kk [\![\symbv_1, \ldots , \symbv_d]\!]$ for $d = \dim_{\kk}(\mm/\mm^2)$, of the form $A = \hat{\symbp} / \symbJo$, where $\symbJo \subseteq \mm_{\hat{\symbp}}^2$. We can always lift the morphism $\kk [\![\symbv_1, \ldots , \symbv_d]\!] \to B_0$ to $B$ (not necessarily uniquely), so we obtain a diagram
    \[\begin{tikzcd}
    	0 & \symbJo & {\hat{\symbp}} & A & 0 \\
    	0 & K & B & {B_0} & 0,
    	\arrow[from=2-3, to=2-4]
    	\arrow[from=2-2, to=2-3]
    	\arrow[from=2-1, to=2-2]
    	\arrow[from=2-4, to=2-5]
    	\arrow[from=1-3, to=1-4]
    	\arrow[from=1-4, to=1-5]
    	\arrow[from=1-2, to=1-3]
    	\arrow[from=1-1, to=1-2]
    	\arrow["{\symbb|_{\symbJo}}", from=1-2, to=2-2]
    	\arrow["\symbb", from=1-3, to=2-3]
    	\arrow["{\symba_0}", from=1-4, to=2-4]
    \end{tikzcd}\]
    where $K = \ker(B \to B_0)$.
    The restriction $\symbb |_{\symbJo}$ doesn't depend on the choice of $\symbb$ and $\mm_B \cdot K = 0$, so we get a homomorphism $\symbJo/\mm \symbJo \to K$. In other words we get an element of $(\symbJo/\mm \symbJo)^{\vee} \otimes_{\kk} K$ which vanishes if and only if the extension $A \to B$ exists. If we denote $Ob = (\symbJo/\mm \symbJo)^{\vee}$ and $D(C) = \Hom(A,C)$ then we get a functorial (under morphisms of small extensions) sequence
    \begin{equation}\label{ref:smallext:eq}
        D(B) \to D(B_0) \to Ob \otimes_{\kk} K.
    \end{equation}
    It is \emph{exact} which in this case means that $\symba_0 \in D(B_0)$ goes to $0$ if and only if there exists $\symba \in D(B)$ which restricts to $\symba_0$. In general obstruction theories provide such exact sequences for every small extension.
    
    \begin{definition}
        Let $D$ be a functor $\Artk \to \Set$ such that $D(\kk) = \{ \textnormal{pt.} \}$ and let $T_D := D(\kk[\varepsilon])$. An \textit{obstruction theory} for $D$ is a finitely dimensional vector space $Ob_D$ over $\kk$ such that for all small extensions $0 \to K \to B \to B_0 \to 0$ there is a map $ob: D(B_0) \to Ob_D \otimes_{\kk} K$ with the following properties:
        \begin{enumerate}
            \item for all $\symba_0 \in D(B_0)$, $ob(\symba_0)=0$ if and only if $\symba_0$ extends to $\symba \in D(B)$,
            \item if $\symba_0 \in D(B_0)$, then the set of extensions of $\symba_0$ to $\symba \in D(B)$ is affine over $T_D \otimes_{\kk} K$ (see \cite{DefHar} for details).
        \end{enumerate}
        We also require that $ob$ is functorial with respect to morphisms of small extensions which means that if we have a commutative diagram where horizontal lines are small extensions
        \[\begin{tikzcd}
        	0 & K & B & {B_0} & 0 \\
        	0 & {K'} & {B'} & {B_0'} & 0,
        	\arrow[from=1-1, to=1-2]
        	\arrow[from=1-2, to=1-3]
        	\arrow[from=1-3, to=1-4]
        	\arrow[from=1-4, to=1-5]
        	\arrow[from=2-1, to=2-2]
        	\arrow[from=2-2, to=2-3]
        	\arrow[from=2-3, to=2-4]
        	\arrow[from=2-4, to=2-5]
        	\arrow["\rho|_K", from=1-2, to=2-2]
        	\arrow["\rho", from=1-3, to=2-3]
        	\arrow["\rho_0", from=1-4, to=2-4]
        \end{tikzcd}\]
        then we get the following commutative diagram
        \[\begin{tikzcd}
        	{D(B)} & {D(B_0)} & {Ob \otimes_{\kk} K} \\
        	{D(B')} & {D(B_0')} & {Ob \otimes_{\kk} K'}.
        	\arrow[from=1-1, to=1-2]
        	\arrow["ob", from=1-2, to=1-3]
        	\arrow["D(\rho)", from=1-1, to=2-1]
        	\arrow[from=2-1, to=2-2]
        	\arrow["id \otimes \rho|_K", from=1-3, to=2-3]
        	\arrow["ob", from=2-2, to=2-3]
        	\arrow["D(\rho_0)", from=1-2, to=2-2]
        \end{tikzcd}\]
    \end{definition}
    
    \begin{example}\label{ref:primary:ex}
        We discuss \emph{primary obstructions}. Let $D \simeq \Hom(A, -)$ be pro-representable by $(A, \mm)$, where $A = \hat{\symbp} / \symbJo$ for some $\symbJo \subseteq \mm_{\hat{\symbp}}^2$. Assume that $D$ has an obstruction theory $Ob$. For $A/\mm^2 \simeq \hat{\symbp}/\mm_{\hat{\symbp}}^2$ we get a homomorphism $\eta : A \to A/\mm^2 \simeq \hat{\symbp}/\mm_{\hat{\symbp}}^2$. We can try to prolong it to a morphism $\mu : A \to \hat{\symbp}/\mm_{\hat{\symbp}}^3$. Using the obstruction theory of $D$ for this homomorphism we get a sequence as in \eqref{ref:smallext:eq} with $B_0 = \hat{\symbp}/\mm_{\hat{\symbp}}^2$, $B = \hat{\symbp}/\mm_{\hat{\symbp}}^3$ and $K = \mm_{\hat{\symbp}}^2/\mm_{\hat{\symbp}}^3$. Hence, an extension $\mu$ exists if and only if $ob(\eta) \in Ob \otimes_{\kk} K$ is $0$. However, in this case an element of $Ob \otimes_{\kk} K$ can be interpreted as a homomorphism $ob_0:K^{\vee} \to Ob$. Consider the identifications
        \begin{equation}\label{eq:KandSym}
            K^{\vee} = \left( \frac{\mm_{\hat{\symbp}}^2}{\mm_{\hat{\symbp}}^3} \right)^{\vee} \simeq \Sym_2 \left( \frac{\mm_{\hat{\symbp}}}{\mm_{\hat{\symbp}}^2}\right)^{\vee} \simeq \Sym_2 \left(\left(\frac{\mm_{\hat{\symbp}}}{\mm_{\hat{\symbp}}^2}\right)^{\vee}\right) = \Sym_2(T_{\mm}A),
        \end{equation}
        where the middle isomorphism comes from duality defined by equation \eqref{eq:dualityw}.
        By identifying $K^{\vee}$ with $\Sym_2(T_{\mm}A)$ we can see $ob_0$ as a map $\Sym_2(T_{\mm}A) \to Ob$, or as a linear map $T_{\mm}A \otimes_{\kk} T_{\mm}A \to Ob$. 
    \end{example}

    Let us note that even though we used a presentation of $A$ as a quotient $A = \hat{\symbp} / \symbJo$, for some $\symbJo \subseteq \mm_{\hat{\symbp}}^2$, the obtained primary obstruction map $ob_0:\Sym_2(T_{\mm}A) \to Ob$ is independent of it. This is because for any other presentation $A = \hat{\symbp}' / \symbJo'$, where $\hat{\symbp}' = \kk [\![\symbv_1', \ldots , \symbv_d']\!]$ and $\symbJo' \subseteq \mm_{\hat{\symbp}'}^2$, it is possible to lift $\hat{\symbp}' \to A$ to a homomorphism $\phi:\hat{\symbp}' \to \hat{\symbp}$. The map $\phi$ induces an isomorphism on tangent spaces, hence it is an isomorphism. Thus by functoriality of obstructions, we get the asserted independence.
    
    \subsection{Calculating primary obstructions} In this subsection we describe some methods of computing primary obstructions. We fix a deformation functor $D \simeq \Hom(A, -)$ pro-representable by $(A, \mm)$ where $A = \hat{\symbp} / \symbJo$ is a quotient of $\hat{\symbp} = \kk [\![\symbv_1, \ldots , \symbv_d]\!]$ by an ideal $\symbJo \subseteq \mm_{\hat{\symbp}}^2$, and $\mm_{\hat{\symbp}}$ is the maximal ideal of $\hat{\symbp}$. We assume that $D$ has an obstruction theory $Ob$ and we use the identification (\ref{eq:KandSym}), namely $K^{\vee} \simeq \Sym_2(T_{\mm}A)$, where $K = \mm_{\hat{\symbp}}^2/\mm_{\hat{\symbp}}^3$. We write $ob_0:\Sym_2(T_{\mm}A) \to Ob$ for the corresponding primary obstruction map.

    \begin{example}\label{ref:primary_compare:ex}
    Let $(A', \mm')$ be a Noetherian complete local $\kk$-algebra with a surjection $\zeta : A \to A'$. Assume that $A' = \hat{\symbT}/\symbJ$ for $\hat{\symbT} := \kk[\![ \sgamma_1, \dots, \sgamma_m]\!]$, $\mm_{\hat{\symbT}} := (\Bar{\sgamma})$, and $\symbJ \subseteq \mm_{\hat{\symbT}}^2$. By the universal property of $\hat{\symbp}$, it is possible to fit $\zeta$ in the following diagram of local $\kk$-algebras 
    \[\begin{tikzcd}
        A & {A'} \\
        {\hat{\symbp}} & {\hat{\symbT}}.
        \arrow[from=2-2, to=1-2]
        \arrow[from=2-1, to=1-1]
        \arrow["\zeta", from=1-1, to=1-2]
        \arrow["\rho", from=2-1, to=2-2]
    \end{tikzcd}\]
    Let $\eta: A \to \hat{\symbp}/\mm_{\hat{\symbp}}^2$  and $\mu : A \to \hat{\symbT}/\mm_{\hat{\symbT}}^2$ be the natural homomorphisms. Consider the morphism of small extensions induced by $\rho$:
    \[\begin{tikzcd}
    	0 & K & {\hat{\symbp}/\mm_{\hat{\symbp}}^3} & {\hat{\symbp}/\mm_{\hat{\symbp}}^2} & 0 \\
    	0 & {K'} & {\hat{\symbT}/\mm_{\hat{\symbT}}^3} & {\hat{\symbT}/\mm_{\hat{\symbT}}^2} & 0,
    	\arrow[from=1-1, to=1-2]
    	\arrow[from=1-2, to=1-3]
    	\arrow["{\rho'|_K}", from=1-2, to=2-2]
    	\arrow[from=1-3, to=1-4]
    	\arrow["\rho'", from=1-3, to=2-3]
    	\arrow[from=1-4, to=1-5]
    	\arrow["{\rho_0'}", from=1-4, to=2-4]
    	\arrow[from=2-1, to=2-2]
    	\arrow[from=2-2, to=2-3]
    	\arrow[from=2-3, to=2-4]
    	\arrow[from=2-4, to=2-5]
    \end{tikzcd}\]
    where $K' = \mm_{\hat{\symbT}}^2/ \mm_{\hat{\symbT}}^3 = \frac{(\sgamma_1, \dots, \sgamma_m)^2}{(\sgamma_1, \dots, \sgamma_m)^3}$. Treating $\eta$ as an element of $D(\hat{\symbp}/\mm_{\hat{\symbp}}^2)$ and $\mu$ as an element of $D(\hat{\symbT}/\mm_{\hat{\symbT}}^2)$ the fixed obstruction theory $Ob$ of $D$ yields elements 
    \[ ob(\eta) \in Ob \otimes_{\kk} K , ob(\mu) \in Ob \otimes_{\kk} K'.\] Note that $(K')^{\vee} \simeq T_{\mm'}A'$ as in (\ref{eq:KandSym}). We identify $T_{\mm'}A'$ as a subspace of $T_{\mm}A$ via $\zeta$. Functioriality of obstructions gives the following commutative diagram
    \[\begin{tikzcd}
        {\Sym_2(T_{\mm}A)} & Ob \\
        {\Sym_2(T_{\mm'}A').}
        \arrow["{ob(\eta)}", from=1-1, to=1-2]
        \arrow["{ob(\mu)}"', from=2-1, to=1-2]
        \arrow["", from=2-1, to=1-1]
    \end{tikzcd}\]
    The map $ob(\eta)$ is the primary obstruction map $ob_0$ by its definition.
    Hence, $ob(\mu)$ is the restriction of primary obstructions to $\Sym_2(T_{\mm'}A')$. 
    \end{example}
    \begin{proposition}\label{prop:primaryobs_tool}
        Fix notation from Example \ref{ref:primary_compare:ex}. Consider the map $ob(\mu)^{\vee} : Ob^{\vee} \to K'$. Then:
        \begin{equation*}
            \im (ob(\mu)^{\vee}) \subseteq \frac{(\sgamma_1, \dots, \sgamma_m)^3 + J}{(\sgamma_1, \dots, \sgamma_m)^3}.
        \end{equation*}
    \end{proposition}
    \begin{proof}
        Consider the commutative diagram
        \[\begin{tikzcd}
            A & {A'} \\
            {\frac{\kk[\![\sgamma_1, \dots, \sgamma_m]\!]}{(\sgamma_1, \dots, \sgamma_m)^3 + J}} & {\frac{\kk[\![\sgamma_1, \dots, \sgamma_m]\!]}{(\sgamma_1, \dots, \sgamma_m)^2}},
            \arrow["\zeta", from=1-1, to=1-2]
            \arrow[from=1-2, to=2-1]
            \arrow[from=1-2, to=2-2]
            \arrow[from=2-1, to=2-2]
            \arrow["\nu"', from=1-1, to=2-1]
            \arrow["\mu"{description}, crossing over, from=1-1, to=2-2]
        \end{tikzcd}\]
        where maps from $A'$ to ${\frac{\kk[\![\sgamma_1, \dots, \sgamma_m]\!]}{(\sgamma_1, \dots, \sgamma_m)^3 + J}}$ and ${\frac{\kk[\![\sgamma_1, \dots, \sgamma_m]\!]}{(\sgamma_1, \dots, \sgamma_m)^2}}$ are projections, and $\nu$ is the composition of $\zeta$ with the projection $A' \to {\frac{\kk[\![\sgamma_1, \dots, \sgamma_m]\!]}{(\sgamma_1, \dots, \sgamma_m)^3 + J}}$. Now look at the map of extensions
        \[\begin{tikzcd}
            0 & {K' = \frac{(\sgamma_1, \dots, \sgamma_m)^2}{(\sgamma_1, \dots, \sgamma_m)^3}} & {B :=\frac{\kk[\![\sgamma_1, \dots, \sgamma_m]\!]}{(\sgamma_1, \dots, \sgamma_m)^3}} & {B_0 := \frac{\kk[\![\sgamma_1, \dots, \sgamma_m]\!]}{(\sgamma_1, \dots, \sgamma_m)^2}} & 0 \\
            0 & {K'' := \frac{(\sgamma_1, \dots, \sgamma_m)^2}{(\sgamma_1, \dots, \sgamma_m)^3 + J}} & {B' := \frac{\kk[\![\sgamma_1, \dots, \sgamma_m]\!]}{(\sgamma_1, \dots, \sgamma_m)^3 + J}} & {B_0' :=\frac{\kk[\![\sgamma_1, \dots, \sgamma_m]\!]}{(\sgamma_1, \dots, \sgamma_m)^2} } & 0.
            \arrow[from=2-2, to=2-3]
            \arrow[from=2-3, to=2-4]
            \arrow[from=2-1, to=2-2]
            \arrow[from=2-4, to=2-5]
            \arrow[from=1-2, to=1-3]
            \arrow[from=1-1, to=1-2]
            \arrow[from=1-4, to=1-5]
            \arrow[from=1-3, to=2-3]
            \arrow["{=}", from=1-4, to=2-4]
            \arrow["u", from=1-2, to=2-2]
            \arrow[from=1-3, to=1-4]
        \end{tikzcd}\]
        By applying $D$, we get the following
        \[\begin{tikzcd}
            {D(B)} & {D(B_0)} & {Ob \otimes_\kk K'} && \mu & {ob(\mu)} \\
            {D(B')} & {D(B_0')} & {Ob \otimes_\kk K''}, & {\nu} & \mu & 0.
            \arrow[from=1-1, to=1-2]
            \arrow[from=1-2, to=1-3]
            \arrow["{id \otimes u}", from=1-3, to=2-3]
            \arrow["{=}", from=1-2, to=2-2]
            \arrow[from=1-1, to=2-1]
            \arrow[from=2-1, to=2-2]
            \arrow[from=2-2, to=2-3]
            \arrow[maps to, from=2-4, to=2-5]
            \arrow[maps to, from=1-5, to=2-5]
            \arrow[maps to, from=2-5, to=2-6]
            \arrow[maps to, from=1-5, to=1-6]
            \arrow[maps to, from=1-6, to=2-6]
        \end{tikzcd}\]
        The fact that $(id \otimes u)(ob'(\mu)) = 0$ translates into
        \[\begin{tikzcd}
            {(K')^{\vee}} & Ob \\
            {(K'')^{\vee}.}
            \arrow["{ob(\mu)}", from=1-1, to=1-2]
            \arrow["{u^{\vee}}", from=2-1, to=1-1]
            \arrow["0"', from=2-1, to=1-2]
        \end{tikzcd}\]
        So $ob(\mu) \circ u^{\vee} = 0$, which means that $\im (ob(\mu)^{\vee}) \subseteq \ker u = \frac{(\sgamma_1, \dots, \sgamma_m)^3 + J}{(\sgamma_1, \dots, \sgamma_m)^3}$.
    \end{proof}
    The following lemma will be useful for calculating primary obstructions. In its statement we use the fact that in the category of rings over $\kk$, the product of $\frac{\kk[\eps]}{(\eps)^2}$ and $\frac{\kk[\eps']}{(\eps')^2}$ is $\frac{\kk[\eps, \eps']}{(\eps, \eps')^2}$ (not to be confused with the tensor product which is $\frac{\kk[\eps, \eps']}{(\eps^2, \eps'^2)}$).
        \begin{lemma}\label{fact:primary_vanish}
            Let $\delta: A \to \frac{\kk[\eps]}{(\eps)^2}$, $\delta' : A \to \frac{\kk[\eps']}{(\eps')^2}$ be homomorphisms of local algebras and denote by $\delta, \delta'$ the corresponding elements in $T_{\mm} A$. Suppose that the induced homomorphism $(\delta \times \delta') : A \to \frac{\kk[\eps, \eps']}{(\eps, \eps')^2}$ can be extended to a homomorphism to $\frac{\kk[\eps, \eps']}{(\eps^2, \eps'^2)}$, i.e., there exists a homomorphism $\gamma$ filling
            \[\begin{tikzcd}
            	A & {\frac{\kk[\eps, \eps']}{(\eps, \eps')^2}} \\
            	& {\frac{\kk[\eps, \eps']}{(\eps^2, \eps'^2)}}.
            	\arrow[from=2-2, to=1-2]
            	\arrow["\delta \times \delta'", from=1-1, to=1-2]
            	\arrow["\gamma", dashed, from=1-1, to=2-2]
            \end{tikzcd}\]
            Then the primary obstruction map $ob_0$ vanishes on $\delta \cdot \delta' \in \Sym_2(T_{\mm} A)$.
        \end{lemma}
        \begin{proof}
            Consider a morphism of small extensions
            \[\begin{tikzcd}
            	0 & {\frac{(\eps, \eps')^2}{(\eps, \eps')^3}} & {\frac{\kk[\eps, \eps']}{(\eps, \eps')^3}} & {\frac{\kk[\eps, \eps']}{(\eps, \eps')^2}} & 0 \\
            	0 & {\frac{(\eps, \eps')^2}{(\eps^2, \eps'^2)}} & {\frac{\kk[\eps, \eps']}{(\eps^2, \eps'^2)}} & {\frac{\kk[\eps, \eps']}{(\eps, \eps')^2}} & 0.
            	\arrow[from=1-2, to=1-3]
            	\arrow[from=1-3, to=1-4]
            	\arrow[from=1-4, to=1-5]
            	\arrow[from=1-1, to=1-2]
            	\arrow[from=2-1, to=2-2]
            	\arrow[from=2-2, to=2-3]
            	\arrow[from=2-3, to=2-4]
            	\arrow["{=}"', from=1-4, to=2-4]
            	\arrow[from=2-4, to=2-5]
            	\arrow[from=1-3, to=2-3]
            	\arrow["s", from=1-2, to=2-2]
            \end{tikzcd}\]
            By applying $D = \Hom(A, -)$ on it we get
            \[\begin{tikzcd}
            	{D\left(\frac{\kk[\eps, \eps']}{(\eps, \eps')^3}\right)} & {D\left(\frac{\kk[\eps, \eps']}{(\eps, \eps')^2}\right)} & {Ob \otimes_\kk \frac{(\eps, \eps')^2}{(\eps, \eps')^3}} && {\delta \times \delta'} & {ob(\delta \times \delta')} \\
            	{D\left(\frac{\kk[\eps, \eps']}{(\eps^2, \eps'^2)}\right)} & {D\left(\frac{\kk[\eps, \eps']}{(\eps, \eps')^2}\right)} & {Ob \otimes_\kk \frac{(\eps, \eps')^2}{(\eps^2, \eps'^2)}}, & \gamma & {\delta \times \delta'} & 0.
            	\arrow[from=1-1, to=1-2]
            	\arrow[from=2-1, to=2-2]
            	\arrow["{=}"', from=1-2, to=2-2]
            	\arrow[from=1-1, to=2-1]
            	\arrow[from=2-2, to=2-3]
            	\arrow["ob", from=1-2, to=1-3]
            	\arrow["{id \otimes s}"', from=1-3, to=2-3]
            	\arrow[maps to, from=2-4, to=2-5]
            	\arrow[maps to, from=1-5, to=2-5]
            	\arrow[maps to, from=2-5, to=2-6]
            	\arrow[maps to, from=1-5, to=1-6]
            	\arrow[maps to, from=1-6, to=2-6]
            \end{tikzcd}\]
            The fact that $(id \otimes s)(ob(\delta \times \delta')) = 0$ translates into a commutative diagram
            \[\begin{tikzcd}
            	{\left(\frac{(\eps, \eps')^2}{(\eps, \eps')^3}\right)^{\vee}} & Ob \\
            	{\left(\frac{(\eps, \eps')^2}{(\eps^2, \eps'^2)}\right)^{\vee}.}
            	\arrow["{ob(\delta \times \delta')}", from=1-1, to=1-2]
            	\arrow["{s^{\vee}}", from=2-1, to=1-1]
            	\arrow["0"', from=2-1, to=1-2]
            \end{tikzcd}\]
            By Example \ref{ref:primary_compare:ex} for $A' = \frac{\kk[\eps, \eps']}{(\eps, \eps')^2}$ we get that $ob(\delta \times \delta')$ is the restriction of the primary obstruction map $ob_0 : \Sym_2(T_{\mm} A) \to Ob$ to the image of $\Sym_2(T_{(\eps, \eps')} \frac{\kk[\eps, \eps']}{(\eps, \eps')^2})$ by the differential $d(\delta \times \delta')$. However, by definition of $\delta \times \delta'$ we have \begin{equation}\label{eq:obs_eq_1}
                d(\delta \times \delta')(\eps^{\vee}) = \delta,
            \end{equation}
            and
            \begin{equation}\label{eq:obs_eq_2}
                d(\delta \times \delta')(\eps'^{\vee}) = \delta'.
            \end{equation}
            The map $s^{\vee}$ is the embedding of $\kk (\eps \cdot \eps')^{\vee}$ into $\Sym_2(\kk \eps^{\vee} \oplus \kk \eps'^{\vee})$. Thus, using the identifications \eqref{eq:obs_eq_1} and \eqref{eq:obs_eq_2} we get that $ob_0(\delta \cdot \delta')=0$. This finishes the proof.
        \end{proof}
    
    \begin{example}\label{ex:hilbert}
        Fix a $\kk$-point $x \in \cH$. Then $x$ corresponds to an ideal $I \subseteq S$ such that $S/I$ of degree $d$. Consider the functor $D_{\cH, x}$ pro-representable by $\hat{\cO}_{\cH, x}$. Its tangent space is $\Hom_{S}(I, S/I)$ by Theorem \ref{thm:tangenttohilb}. It is also known that $D_{\cH, x}$ has an obstruction theory with the obstruction space $Ob = \Ext^1_{S}(I, S/I)$. For a proof see \cite[Corollary 6.4.10]{FGI+05}. Thus in this case the primary obstruction map is a function $ob_0 : \Sym_2(\Hom_{S}(I, S/I)) \to \Ext^1_{S}(I, S/I)$.
    \end{example}
    
    We will need the method of calculating primary obstructions for the functor in Example \ref{ex:hilbert}. It comes from the paper \cite{JS21}.
    
    Consider an ideal $I \subseteq S$ and its free resolution as an $S$-module $(F_i)_{i \geq 1}$. Then the $S$-module $S/I$ has a free resolution which starts with the quotient map $q:F_0 = S \to S/I$ and later is continued by the resolution $(F_i)_{i \geq 1}$ using the composition $F_1 \to I \subseteq S = F_0$. For a given $\symbc \in \Hom_S(I, S/I)$, we can lift it to $S$-module homomorphisms as on the following diagram (not necessarily uniquely)
    
    \[\begin{tikzcd}\label{diagram:hilb_obtruction}
    	& 0 & I & {F_1} & {F_2} & {F_3} & \dots \\
    	0 & {S/I} & {F_0} & {F_1} & {F_2} & \dots
    	\arrow[from=1-3, to=1-2]
    	\arrow["{d_0}"', from=1-4, to=1-3]
    	\arrow["{d_1}"', from=1-5, to=1-4]
    	\arrow["{d_2}"', from=1-6, to=1-5]
    	\arrow[from=2-6, to=2-5]
    	\arrow["{d_1}"', from=2-5, to=2-4]
    	\arrow["{d_0}"', from=2-4, to=2-3]
    	\arrow["q"', from=2-3, to=2-2]
    	\arrow[from=2-2, to=2-1]
    	\arrow[from=1-7, to=1-6]
    	\arrow["\symbc"'{pos=0.7}, from=1-3, to=2-2]
    	\arrow["{s_1(\symbc)}"'{pos=0.8}, from=1-4, to=2-3]
    	\arrow["{s_2(\symbc)}"'{pos=0.8}, from=1-5, to=2-4]
    	\arrow["{s_3(\symbc)}"'{pos=0.8}, from=1-6, to=2-5]
    \end{tikzcd}\]
    
    Under this notation the following theorem holds.
    
    \begin{theorem} \cite[Theorem 4.18]{JS21}\label{thm:hilb_obstruction}
        Consider the point point of the Hilbert scheme $\cH$ corresponding to an ideal $I \subseteq S$ and the associated obstruction theory with obstruction group $Ob = \Ext^1_{S}(I, S/I)$. Then for $\symbc_1, \symbc_2 \in \Hom_{S}(I, S/I)$ the primary obstruction map sends the class of $\symbc_1 \cdot \symbc_2 \in \Sym_2 (\Hom_{S}(I, S/I))$ to the class of
        \begin{equation}\label{eq:hilb_obstruction}
            q \circ (s_1(\symbc_1) \circ s_2(\symbc_2) + s_1(\symbc_2) \circ s_2(\symbc_1))
        \end{equation}
        in $\Ext^1_{S}(I, S/I)$. 
    \end{theorem}
    Note that $q \circ s_1(\symbc_1) \circ s_2(\symbc_2) \circ d_2 = \symbc_1 \circ d_0 \circ d_1 \circ  d_3(\symbc_2) = \symbc_1 \circ 0 \circ d_3(\symbc_2) = 0$, so \eqref{eq:hilb_obstruction} indeed yields an element of $\Ext^1_{S}(I, S/I)$.
    \subsection{Barycenter}\label{subsection:barycenter}
    It will be convenient for us not to work with the Hilbert scheme $\cH$ itself, but with its closed subscheme consisting of tuples which are centered around $0$. In order to construct it we need to know that the support map which takes a point of the Hilbert scheme and associates with it its cycle theoretic support (in the symmetric product of the ambient variety) is a morphism of schemes. More precisely, we have the following theorem from \cite[Theorem 7.1.14]{FGI+05}:
    
    \begin{theorem}\label{thm:HC}
        There is a morphism of schemes called the Hilbert-Chow morphism
        \begin{equation*}
            \theta : \cH \to (\Affine^n)^d \goodquotient \mathbb{S}_d,
        \end{equation*}
        which on the level of $\kk$-points is given by: $[Z] \mapsto \sum_{p} \deg(\cO_{Z,p}) [p]$, where the sum ranges over $p \in \Affine^n$.
    \end{theorem}
    
    Here the quotient $(\Affine^n)^d \goodquotient \mathbb{S}_d$ is by the definition the spectrum of the algebra of global sections of $(\Affine^n)^d$ invariant under the action of the symmetric group $\mathbb{S}_d$ acting by permuting entries of the $d$-tuples.
    
    Now we will define the barycenter scheme. The morphism $(\Affine^n)^d \to \Affine^n$ given by taking the average of $d$ vectors is $\mathbb{S}_d$-equivariant, so it descends to a morphism $Av:(\Affine^n)^d \goodquotient \mathbb{S}_d  \to \Affine^n$. By taking composition with the Hilbert-Chow morphism we get the barycenter morphism
    \begin{equation}\label{equation:barycenter}
        Bar: \cH \to \Affine^n.
    \end{equation}
    \begin{definition}\label{def:Bar}
        The barycenter scheme is $\cB := Bar^{-1}(\{ 0 \})$. If we want to emphasize $d$ we write $\cB^d$.
    \end{definition}
    
    To study the barycenter scheme we give an explicit description of the Hilbert-Chow morphism, or more precisely a formula for the composition $\Spec(A) \to \cH \to (\Affine^n)^d \goodquotient \mathbb{S}_d$ for a ring $A$ such that $\Spec(A) \to \cH$ is given by a free family $R$ over $A$. For $r \in R$ let $Nm_{R/A}(r)$ and $Tr_{R/A}(r)$ be respectively the determinant and the trace of the endomorphism of $R$ (as an $A$-module) given by multiplication by $r$. Then the following theorem from the notes \cite[Theorem 2.16, Corollary 2.18]{notesbertin} holds.
    
    \begin{theorem}\label{thm:HCaffine}
        Suppose a morphism $\Spec(A) \to \cH$ is given by a free $A$ - algebra $R$. Then the composition
        \[\begin{tikzcd}
        	{\Spec(A)} & {\cH} & {(\Affine^n)^d \goodquotient \mathbb{S}_d}
        	\arrow[from=1-1, to=1-2]
        	\arrow[from=1-2, to=1-3]
        \end{tikzcd}\]
        comes from the composite homomorphism
        \[\begin{tikzcd}
        	{(\kk[\aalpha] \otimes_\kk \dots \otimes_\kk \kk[\aalpha])^{\mathbb{S}_d}} & {(R \otimes_A \dots \otimes_A R)^{\mathbb{S}_d}} & A,
        	\arrow[from=1-1, to=1-2]
        	\arrow["LNm", from=1-2, to=1-3]
        \end{tikzcd}\]
        where the first map is the homomorphism $\kk[\aalpha] \to \kk[\aalpha] \otimes_\kk A \to R$ tensored $d$ times, and the second map $LNm$ is the unique $A$-algebra homomorphism that satisfies the equation
        \begin{equation*}
            LNm(r \otimes \dots \otimes r) = Nm_{R/A}(r).
        \end{equation*}
        Moreover, if we consider the composition of the morphism $\Spec(A) \to (\Affine^n)^d \goodquotient \mathbb{S}_d$ with the average morphism $Av : (\Affine^n)^d \goodquotient \mathbb{S}_d \to \Affine^n$ we get a morphism $\Spec(A) \to \Affine^n = \Spec(\kk[\aalpha])$ which on the level of algebras is given by the formula
        \begin{equation*}
            \alpha_j \mapsto \frac{1}{d} Tr_{R/A}(\alpha_j),
        \end{equation*}
        where on the right hand side we identify $\alpha_j$ with its image in $R$.
    \end{theorem}
    
    Every point of the Hilbert scheme can be uniquely shifted so that it is centered at $0$, so there is an isomorphism $\cH \simeq \Affine^n \times \cB$ as one can prove using for example Theorem \ref{thm:HCaffine}.
    \begin{proposition}\label{prop:hilbandbar}
            The barycenter morphism $Bar :\cH \to \Affine^n$ is equivariant under the action of the $n$-dimensional additive group $\G_a^n$. Hence, we have $\cH \simeq \G_a^n \times \cB$.
    \end{proposition}
    \begin{proof}
        First we prove that the barycenter morphism $Bar :\cH \to \Affine^n$ is $\G_a^n$-equivariant. It means that the following diagram commutes:
        \[\begin{tikzcd}
        	{\G_a^n \times \cH} & {\G_a^n \times \Affine^n} \\
        	\cH & {\Affine^n}.
        	\arrow["\mu", from=1-1, to=2-1]
        	\arrow["{+}", from=1-2, to=2-2]
        	\arrow["{id \times Bar}", from=1-1, to=1-2]
        	\arrow["Bar", from=2-1, to=2-2]
        \end{tikzcd}\]
        Here $\mu$ is the $\G_a^n$-action on $\cH$. It is enough to check that for all $\kk$-algebras $A$ the following diagram commutes:
        \[\begin{tikzcd}
        	{\G_a^n(\Spec(A)) \times \cH(\Spec(A))} & {\G_a^n(\Spec(A)) \times \Affine^n(\Spec(A))} \\
        	{\cH(\Spec(A))} & {\Affine^n(\Spec(A))}.
        	\arrow["\mu", from=1-1, to=2-1]
        	\arrow["{+}", from=1-2, to=2-2]
        	\arrow["{id \times Bar}", from=1-1, to=1-2]
        	\arrow["Bar", from=2-1, to=2-2]
        \end{tikzcd}\]
        Fix a pair $(\bar{a}, [Z]) \in \G_a^n(\Spec(A)) \times \cH(\Spec(A))$ where $\bar{a} = (a_1, \dots, a_n) \in A^n$. By passing to a localisation, we can assume that $[Z] \in \cH(\Spec(A))$ is given by a free $A$-algebra $R$ of degree $d$. By Theorem \ref{thm:HCaffine} we calculate the images of $(\bar{a}, [Z])$ by both ways and get:
        \[\begin{tikzcd}
        	{(\bar{a}, [Z])} & {(\bar{a}, (\frac{1}{d} \cdot Tr_{R/A}(\alpha_i))_{i=1}^n)} \\
        	{[\bar{a}+Z]} & {(\frac{1}{d} \cdot Tr_{R/A}(\alpha_i + a_i))_{i=1}^n} & {(a_i + \frac{1}{d} \cdot Tr_{R/A}(\alpha_i))_{i=1}^n}.
        	\arrow[maps to, from=1-1, to=1-2]
        	\arrow[maps to, from=1-1, to=2-1]
        	\arrow[maps to, from=1-2, to=2-3]
        	\arrow[maps to, from=2-1, to=2-2]
        	\arrow["{?}"{marking}, draw=none, from=2-2, to=2-3]
        \end{tikzcd}\]
        The image of $[\bar{a} + Z]$ by $Bar$ is the tuple $(\frac{1}{d} \cdot Tr_{R/A}(\alpha_i + a_i))_{i=1}^n$, because the ideal of $\bar{a} + Z$ in $\Affine_A^n$ is given be the kernel of the composite homomorphism
        \[\begin{tikzcd}
        	{S_A} & {S_A} & R,
        	\arrow[two heads, from=1-2, to=1-3]
        	\arrow["{(\bar{a} +)^{\#}}", from=1-1, to=1-2]
        \end{tikzcd}\]
        as in Example \ref{ex:additiveaction} and under this composition $\alpha_i$ goes to $\alpha_i + a_i \in R$.
        By calculation
        \begin{equation*}
            \frac{1}{d} \cdot Tr_{R/A}(\alpha_i + a_i) = \frac{1}{d} \cdot (Tr_{R/A}(\alpha_i) + Tr_{R/A}(a_i)) = \frac{1}{d} \cdot (Tr_{R/A}(\alpha_i) + d \cdot a_i) = \frac{1}{d} \cdot Tr_{R/A}(\alpha_i) + a_i,
        \end{equation*}
        we get that the question mark in the diagram chase can be replaced with equality. This finishes the part proving that $Bar$ is $\G_m$-equivariant.
        
        Now we pass to the construction of an isomorphism $\cH \simeq \G_a^n \times \cB$.
        Consider the following diagram:
        \begin{equation}\label{diag:proof_of_equiv_of_bar}
        \begin{tikzcd}
        	{\G_a^n \times \cB} & {\G_a^n \times \{0\}} \\
        	{\G_a^n \times \cH} & {\G_a^n \times \Affine^n} \\
        	\cH & {\Affine^n},
        	\arrow["{id \times 0}", from=1-1, to=1-2]
        	\arrow["{id \times Bar}", from=2-1, to=2-2]
        	\arrow["Bar", from=3-1, to=3-2]
        	\arrow["{+}", from=2-2, to=3-2]
        	\arrow["\mu", from=2-1, to=3-1]
        	\arrow["{id \times j}", from=1-1, to=2-1]
        	\arrow["{id \times 0}"', hook, from=1-2, to=2-2]
        	\arrow["\simeq", shift left=5, bend left, from=1-2, to=3-2]
        	\arrow[shift right=2, bend right, from=1-1, to=3-1]
        \end{tikzcd}
        \end{equation}
        where $j: \cB \to \cH$ is the closed embedding. The bottom square is commutative because $Bar$ is  $\G_a^n$-equivariant. The upper square is a pullback diagram, because of the definition of $\cB$. The bottom square is also a pullback, which follows from the existence of group inverse in $\G_a^n$. Thus the left hand side map $\G_a^n \times \cB \to \cH$ is a pullback of an isomorphism, so it is also an isomorphism and we get the conclusion.
    \end{proof}
    The following proposition shows that the isomorphism $\cH \simeq \G_a^n \times \cB$ respects the torus action.
    \begin{proposition}\label{prop:barisgminv}
        The barycenter morphism $Bar : \cH \to \Affine^n$ is $\G_m$-equivariant. Moreover, the isomorphism $\cH \simeq \G_a^n \times \cB$ from Proposition \ref{prop:hilbandbar} is $\G_m$-equivariant, where the torus action on $\cB$ is the one induced from the $\G_m$-action on $\cH$.
    \end{proposition}
    \begin{proof}
        We need to prove commutativity of the diagram
        \[\begin{tikzcd}
        	{\G_m \times \cH} & {\G_m \times \Affine^n} \\
        	\cH & {\Affine^n},
        	\arrow["\nu", from=1-1, to=2-1]
        	\arrow["{\cdot}", from=1-2, to=2-2]
        	\arrow["{id \times Bar}", from=1-1, to=1-2]
        	\arrow["Bar", from=2-1, to=2-2]
        \end{tikzcd}\]
        where $\nu$ is the $\G_m$-action on $\cH$. Exactly as in the proof of Proposition \ref{prop:hilbandbar}, it is enough to prove that for a $\kk$-algebra $A$ the following diagram commutes:
        \[\begin{tikzcd}
        	{\G_m(\Spec(A)) \times \cH(\Spec(A))} & {\G_m(\Spec(A)) \times \Affine^n(\Spec(A))} \\
        	{\cH(\Spec(A))} & {\Affine^n(\Spec(A))}.
        	\arrow["\nu", from=1-1, to=2-1]
        	\arrow["\cdot", from=1-2, to=2-2]
        	\arrow["{id \times Bar}", from=1-1, to=1-2]
        	\arrow["Bar", from=2-1, to=2-2]
        \end{tikzcd}\]
        Again, we assume that $\Spec(A) \to \cH$ is given by a free $A$-module $R$ of degree $d$ and we fix $a \in \G_m(\Spec(A)) = A^*$. By moving the pair $(a, [Z])$ along two paths on the diagram we get:
        \[\begin{tikzcd}
        	{(a, [Z])} & {(a, (\frac{1}{d} \cdot Tr_{R/A}(\alpha_i))_{i=1}^n)} \\
        	{[a \cdot Z]} & {(\frac{1}{d} \cdot Tr_{R/A}(a \cdot \alpha_i))_{i=1}^n} & {(a \cdot \frac{1}{d} \cdot Tr_{R/A}(\alpha_i))_{i=1}^n},
        	\arrow[maps to, from=1-1, to=1-2]
        	\arrow[maps to, from=1-1, to=2-1]
        	\arrow[maps to, from=1-2, to=2-3]
        	\arrow[maps to, from=2-1, to=2-2]
        	\arrow["{=}"{marking}, draw=none, from=2-2, to=2-3]
        \end{tikzcd}\]
        where the right bottom equality follows from the fact that $Tr_{R/A}$ is an $A$-module homomorphism. This finishes the proof of $\G_m$-equivariance of $Bar$.
        
        Now we prove the "moreover" part of the statement. Recall that the isomorphism $\cH \simeq \G_a^n \times \cB$ comes from the diagram (\ref{diag:proof_of_equiv_of_bar}). Since we have proved that $Bar$ is $\G_m$-equivariant it is immediate that all maps on this diagram beside $\mu$ are $\G_m$-equivariant. The $\G_m$-equivariance of $\mu$ is equivalent to commutativity of
        \[\begin{tikzcd}
        	{\G_m \times(\G_a^n \times \cH)} & {\G_m \times \cH} \\
        	{\G_a^n \times \cH} & \cH.
        	\arrow["{id \times \mu}", from=1-1, to=1-2]
        	\arrow[from=1-1, to=2-1]
        	\arrow["\mu"', from=2-1, to=2-2]
        	\arrow["\nu", from=1-2, to=2-2]
        \end{tikzcd}\]
        We check this on $\Spec(A)$-point of $\G_m \times(\G_a^n \times \cH)$ given by the triple $(a, (\bar{b}, [Z])) \in (\G_m \times(\G_a^n \times \cH))(\Spec(A))$. We have
        \[\begin{tikzcd}
        	{(a, (\bar{b}, [Z]))} & {(a, [\bar{b} +Z])} \\
        	{(a \cdot \bar{b}, [a \cdot Z])} & {[a \cdot \bar{b} + a \cdot Z]} & {[a \cdot (\bar{b} + Z)]},
        	\arrow[maps to, from=1-1, to=1-2]
        	\arrow[maps to, from=1-2, to=2-3]
        	\arrow[maps to, from=1-1, to=2-1]
        	\arrow[maps to, from=2-1, to=2-2]
        	\arrow["{=}"{marking}, draw=none, from=2-2, to=2-3]
        \end{tikzcd}\]
        where the equality follows from the distributivity of multiplication over addition. Thus $\mu$ is indeed $\G_m$-equivariant and so also the isomorphism between $\G_a^n \times \cB$ and $\cH$ which is the composition of $\G_m$-equivariant morphisms $id \times j$ and $\mu$ is $\G_m$-equivariant.
    \end{proof}
    
    The $\G_m$-action on $\cB$ inherited from $\cH$ allows us to perform the \BBname{} decomposition on $\cB$ and get the diagram
    \[\begin{tikzcd}
    	{\cB^{+}} & {\cB}. \\
    	{\cB^{\G_m}}
    	\arrow["\pi", from=1-1, to=2-1]
    	\arrow["i"', from=1-1, to=1-2]
    	\arrow["s", bend left, from=2-1, to=1-1]
    \end{tikzcd}\]
    For a $\G_m$-fixed point $[I] \in \cB$ we identify it with $[I] \in \cB^{\G_m}$ and $s([I]) \in \cB^+$ in the spirit of Fact \ref{fact:mono}. Now we can finally give the definition of a hedgehog point.
    \begin{definition}\label{def:hedgehog:precise}
        The \emph{negative spike} at a $\G_m$-fixed point $[I]$ is the fiber $V = \pi^{-1}([I]) \subseteq \cB^+$. The point $[I]$ is said to be a \emph{hedgehog point} if the following conditions are satisfied:
        \begin{enumerate}
            \item $(T_{[I]} \cB)_{>0} = 0$,
            \item $T_{[I]} \cB \neq (T_{[I]} \cB)_0$,
            \item $V$ is $0$-dimensional.
        \end{enumerate}
    \end{definition}
    The equivalence of this formulation and the one in Definition \ref{def:hedgehog} follows from the fact that $V$ has a nontrivial $\Affine^1$ monoid action, so it is $0$-dimensional if and only if it is a point (topologically).
    \subsection{Macaulay duality}
    In order to find hedgehog points on the barycenter scheme we will use the construction of \emph{apolar algebras} following \cite{Jel17}.
    \begin{definition}\label{def:duality}
        For $P = \kk[x_1, \dots, x_n]$ let $\alpha_1 := \frac{\partial}{\partial x_1}, \dots, \alpha_n := \frac{\partial}{\partial x_n}$. Then the polynomial algebra $S = \kk[\overline{\alpha}]$ generated by commuting $n$ derivatives naturally acts on $P$. This action is called evaluation and is denoted by
        \begin{equation*}
            \ \wcirc \ : S \times P \to P.
        \end{equation*}
        For $f \in P$ we also define $ev_f : S \to P$ by the formula $ev_f(h) := h \ \wcirc \ f$.
        If $I \subseteq P$ is a $P$-submodule we define its dual as the ideal $I^{\perp} := Ann(I) = \{ h \in \kk[\overline{\alpha}] : ev_f(h) = 0 \textnormal{ for all } f \in I\}$.
        The algebra $Apolar(I) := \kk[\overline{\alpha}]/I^{\perp}$ is called the \emph{apolar algebra} of $I$. When $I$ is generated by a single polynomial $F$ we simply write $Apolar(F)$ for $Apolar(I)$ and $\Fperp$ for $I^{\perp}$. Moreover, if $F$ is homogeneous, then so is $\Fperp$ hence $Apolar(F)$ is graded.
    \end{definition}
    
    \begin{example}
        If $I$ in Definition \ref{def:duality} is generated by $F$, then the ideal $I^{\perp} = \Fperp$ is given by the kernel of the map $ev_F:\kk[\overline{\alpha}] \to \kk[\overline{x}]$. The quotient $Apolar(F) = \kk[\overline{\alpha}]/\Fperp$ is isomorphic, as a linear space over $\kk$,  to the image \[ ev_F(\kk[\overline{\alpha}]) = \Big\langle F, \frac{\partial F}{\partial x_1}, \dots, \frac{\partial F}{\partial x_n}, \frac{\partial^2 F}{\partial x_i \partial x_j}, \dots, 1 \Big\rangle. \] For example if $F = x^2 + y^2, \alpha = \frac{\partial}{\partial x}, \beta = \frac{\partial}{\partial y}$, then $Apolar(F) = \frac{\kk[\alpha,\beta]}{(\alpha^2 - \beta^2, \alpha \beta)+ (\alpha,\beta)^3} \simeq \langle x^2+y^2, x, y, 1 \rangle$ is four-dimensional, so it gives a point in $\cH^4$.
    \end{example}
    Later, in Subsection \ref{sub:flatness} we will need the relative Macaulay duality. We introduce it there in the case $n=6$.
    
\section{The main example}\label{sec:main}
    The goal of this section is to prove Theorem \ref{thm:Hedgehog_point_theorem} and analyse the geometry around non-reduced hedgehog points coming from some cubics. We fix symbols $i, s, \pi$ for the maps
    \[\begin{tikzcd}
    	{\cB^{+}} & {\cB} \\
    	{\cB^{\G_m}}
    	\arrow["\pi", from=1-1, to=2-1]
    	\arrow["i"', from=1-1, to=1-2]
    	\arrow["s", bend left, from=2-1, to=1-1]
    \end{tikzcd}\]
    on the \BBname{} decomposition for the $\G_m$-action on $\cB$ induced by the $\G_m$-action on $\cH$. This makes sense as $\cB \subseteq \cH$ is $\G_m$-invariant by Proposition \ref{prop:barisgminv}.
    
    We start by proving the Hedgehog point theorem (Theorem \ref{thm:Hedgehog_point_theorem}) which we recall here.
    \begin{theorem}[Hedgehog point theorem]\label{thm:Hedgehog_point_theorem:precise}
        A hedgehog point $[I] \in \cB$ is non-reduced. Hence $(0, [I]) \in \G_a^n \times \cB \simeq \cH$ is non-reduced as well.
    \end{theorem}
    \begin{proof}
        Assume on the contrary that $\cB$ is reduced at $[I]$. 
        By Proposition \ref{ref:comparison:intro:prop} the map $i$ is an isomorphism on some open neighbourhood of $[I]$, because $(T_{[I]} \cB)_{>0} = 0$ as $[I]$ is a hedgehog point. Hence, $\cB^+$ is also reduced at $[I]$.
        From affineness of $\pi$, see Theorem \ref{ref:introRepresentability:thm}, we get that $\cB^+$ corresponds to some sheaf $\cM$ of $\cO_{\cB^{\G_m}}$-algebras and $\pi$ is the natural morphism from the relative $\Spec$ of $\cM$ to $\cB^{\G_m}$. Moreover, the sheaf $\cM$ has an $\mathbb{N}$-grading (as an $\cO_{\cB^{\G_m}}$-module) coming from $\Affine^1$-action. Thus, with the convention that $\dim \emptyset = -1$, for a point $p \in \cB^{\G_m}$ we get
        \begin{equation*}
            \dim \pi^{-1}(p) = \dim \Spec_{\cO_{\cB^{\G_m}}}(\cM)|_{ \{ p \} } = \dim \Proj_{\cO_{\cB^{\G_m}}}(\cM)|_{ \{ p \} } + 1.
        \end{equation*}
        Thanks to the third condition of being a hedgehog point we know that the fiber of $\pi$ over $[I]$ is $0$-dimensional, so $\dim \Proj_{\cO_{\cB^{\G_m}}}(\cM)|_{ \{ [I] \} } = -1$. By closedness of the image of a projective morphism there is an open neighbourhood $U$ of $[I]$ inside $\cB^{\G_m}$ such that for every $p \in U$ we have $\dim \Proj_{\cO_{\cB^{\G_m}}}(\cM)|_{ \{ p \} } = -1$, so $\dim \pi^{-1}(p) = 0$. These fibers all have a unique fixed point of $\Affine^1$-action so they are all topologically points. Thus we get a morphism
        \[\begin{tikzcd}
        	{\pi^{-1}(U)}, \\
        	U
        	\arrow["\pi|_{\pi^{-1}(U)}", from=1-1, to=2-1]
        \end{tikzcd}\]
        which is affine and on the level of closed points is a bijection. Note that the section $s|_U$ is inverse to $\pi|_{\pi^{-1}(U)}$ on the level of topological spaces, so $\pi|_{\pi^{-1}(U)}$ is in fact a homeomorphism. Consider the diagram
        \[\begin{tikzcd}
        	U && {\pi^{-1}(U)} \\
        	& U.
        	\arrow["{=}"', from=1-1, to=2-2]
        	\arrow["s|_U", from=1-1, to=1-3]
        	\arrow["\pi|_{\pi^{-1}(U)}", from=1-3, to=2-2]
        \end{tikzcd}\]
        The morphism $\pi$ is affine, so separated. Thus, by Cancellation theorem \cite[10.1.19]{Vakil} we get that $s|_U$ is a closed embedding. But it is also a homeomorphism, so it is in fact an isomorphism of schemes. Hence, $T_{[I]} \cB^{\G_m} = T_{[I]} U \simeq T_{[I]} \pi^{-1}(U) = T_{[I]} \cB^+ = T_{[I]} \cB$. This gives us a contradiction, because $T_{[I]} \cB^{\G_m} = (T_{[I]} \cB)_0$ and $[I]$ is a hedgehog point so $T_{[I]} \cB \neq (T_{[I]} \cB)_0$.
    \end{proof}
    
    Now we concentrate on the search for hedgehog points. We fix $n=6$ and $d=13$. Therefore $\cH = \mathrm{Hilb}^{13}_{\Affine^6/\kk}$ and $\cB$ is the corresponding barycenter scheme so $\cH \simeq \G_a^6 \times \cB$. Moreover, $\cH^{14}$ denotes the Hilbert scheme of $14$ points on $\Affine^6$. We start with the definition of being general enough.
    \begin{definition}\label{def:nice}
        Let $F$ be a homogeneous cubic in six variables $x_1, \dots, x_6$. We call $F$ \emph{general enough} if the following conditions are satisfied:
        \begin{enumerate}
            \item $(F, \frac{\partial F}{\partial x_1}, \dots, \frac{\partial F}{\partial x_6}, x_1, \dots, x_6, 1)$ are linearly independent, so form a basis of $\im(ev_F)$ (from Definition \ref{def:duality}). In particular, $\Fperp$ gives a point $[\Fperp] \in \cH^{14}$.
            \item The minimal graded free resolution of $\Fperp$ starts with
            \[\begin{tikzcd}
            	0 & \Fperp & {S^{\oplus 15}(-2)} & {S^{\oplus 35}(-3)} & \dots
            	\arrow[from=1-5, to=1-4]
            	\arrow[from=1-4, to=1-3]
            	\arrow[from=1-3, to=1-2]
            	\arrow[from=1-2, to=1-1]
            \end{tikzcd}\]
            \item $[\Fperp] \in \cH^{14}$ has TNT (which is equivalent to $\dim (T_{[\Fperp]} \cH^{14})_{<0} = 6$, see Definition \ref{def:TNT}).
        \end{enumerate}
    \end{definition}
    First, condition (1) assures that $[\Fperp]$ gives us a point in $\cH^{14}$. Secondly, the condition about the minimal graded free resolution encodes two properties: it implies that: (a) $\Fperp$ is generated by its second degree $(\Fperp)_2$ (which is $15$-dimensional, because $\dim_\kk (\Fperp)_2 = \dim_\kk S_2 - \dim_\kk (S/\Fperp)_2 = 21 - 6 = 15$), and (b) the first syzygies are linear. More thorough analysis of this resolution will appear in Subsection \ref{sub:tangent}. The TNT condition assures tame behaviour of the tangent space to $[\Fperp]$.
     
     Notice that all of these conditions are open (on the moduli space of cubics). Indeed, condition (1) is open, because it can be translated into a statement about the rank of some minors of $ev_F:S \to P$. The openness of the condition (2) follows from semi-continuity of Betti numbers (for the definition of Betti numbers see \cite[Chapter 1B]{Eis1}). At last, the openness of the condition (3) is implied by the semi-continuity of the function $[I] \mapsto \dim\left((T_{[I]}\cH^{14})_{-1}\right)$ for $[I] \in (\cH^{14})^{\G_m}$, together with the fact that $6$ is the minimum possible dimension of the negative part of the tangent space. Since the $\G_m$-action on $\Affine^6$ lifts to the restriction of tangent sheaf of $\cH^{14}$ to $(\cH^{14})^{\G_m}$, it decomposes as a sum of sheaves of integer degrees and the result follows from half-continuity of the dimension of a coherent sheaf on a scheme. Moreover, the set of general enough cubics is not only open, but also non-empty, as the following example shows.
    \begin{example}\label{ex:nonempty}
        Let $F = x_1 x_2 x_4 - x_1 x_5^2 + x_2 x_3^2 + x_3 x_5 x_6 + x_4 x_6^2$ as in Example \ref{ex:first}. One can check by hand, or using a computer, that this $F$ satisfies conditions (1) and (2) in Definition \ref{def:nice}. The fact that it satisfies the condition (3) can be deduced from \cite[lemma 3.6]{Jel18}.
    \end{example}
    
    \subsection{The analysis of the point $[I]$}\label{sub:tangent}
    Now we describe a procedure of making a point $[I] \in \cH$ from a general enough cubic $F$. By the definition of being general enough the $\kk$-linear space $\left\langle F, \frac{\partial F}{\partial x_1}, \dots, \frac{\partial F}{\partial x_6}, \frac{\partial^2 F}{\partial x_i \partial x_j}, \dots, 1 \right\rangle$ is $14$-dimensional and $\left(F, \frac{\partial F}{\partial x_1}, \dots, \frac{\partial F}{\partial x_6}, x_1, \dots, x_6, 1\right)$ is its basis. As a consequence, the algebra $S/\Fperp$ has a natural grading with the Hilbert function $(1,6,6,1)$. Hence the third degree of $S/\Fperp$ is generated by a single element $g$ and we can assume that $ev_F(g) = 1$. We fix a representative of the class of $g$ in $S$ and denote it also by $g$.
    Let $I = (\Fperp, g) = \Fperp + S_{\geq 3} \subseteq S$. Then $S/I \simeq Apolar(F)/(g)$ is $13$-dimensonal, so it represents a point $[I] \in \cH$. Note that the ideal $I$ is homogeneous and the only maximal ideal which contains it is the ideal of $0 \in \Affine^6$. Thus, $[I]$ is in fact an element of $\cB^{\G_m} \subseteq \cB \subseteq \cH$ by the description of the Hilbert-Chow morphism in Theorem \ref{thm:HC} and Definition \ref{def:Bar}. Furthermore, $S/I$ is an example of a \textit{very compressed} algebra.  This means that the radical of $I$ is the ideal of the origin, and that its Hilbert function is of the form
    \begin{equation*}
    H(i) = \begin{cases}
    {i+n-1}\choose{n-1} &\text{for $i < s $,}\\
    d - {{n+s-1}\choose{n}} &\text{for $i= s$,}\\
    0 &\text{for $i > s$,}
    \end{cases}
    \end{equation*}
    where $n$ is the number of variables, $d$ is the degree of the algebra, and $s$ is a natural number such that ${{n+s-1}\choose{n}} < d <{{n+s}\choose{n}}$ (in our case the Hilbert function has values $(1,6,6)$). Very compressed algebras can be used to prove non-reducibility of certain Hilbert schemes of points, and lead to examples of non-smoothable subschemes of $\Affine^4$, see \cite{Iar72b, EI78, Sha90}.
    
    From now on we fix a general enough cubic $F$ and $g$ such that $ev_F(g)=1$. Also, $[I]$ denotes a point of $\cH$ constructed from $F$ by taking $I = (\Fperp, g) = \left\langle \frac{\partial F}{\partial x_1}, \dots, \frac{\partial F}{\partial x_6} \right\rangle^{\perp}$.
    For our applications, the analysis of the $\pi$-fiber over $[I]$ will be important, so we fix a notation $V := \pi^{-1}([I])$. Note that $V$ depends on $I$, however we omit it in this notation. Recall that by Theorem \ref{ref:introRepresentability:thm}, the fiber $V$ has a natural $\Affine^1$-action and it is affine as the inverse image of an affine subscheme of $\cB^{\G_m}$. The structure of the rest of this section is as follows. First we  study the tangent space $T_{[I]} \cH$ and check that it satisfies assumptions of Proposition \ref{ref:comparison:intro:prop}. Then we take a closer look at the resolution of the ideal $I$ and using primary obstructions we prove that $V$ is $0$-dimensional. At the end we study a deformation of $S/I$ coming from the fractal family, to determine the structure of $V$ completely.
    
    We start with the following proposition.
    \begin{proposition}\label{prop:nopositiveH}
            The induced action of $\G_m$ on $T_{[I]} \cH$ has no positive degrees.
    \end{proposition}
    \begin{proof}
        This is because $S/I$ is very compressed.
    \end{proof}
    
    \begin{corollary}\label{cor:nopositivedeg}
        The induced action of $\G_m$ on $T_{[I]} \cB$ has no positive degrees.
    \end{corollary}
    \begin{proof}
        This follows from the Proposition \ref{prop:nopositiveH}, because $T_{[I]} \cB$ is a $\G_m$-equivariant linear subspace of $T_{[I]} \cH$.
    \end{proof}
    
    \begin{corollary}\label{cor:bplusisb}
        The map $i: \cB^+ \to \cB$ is an open embedding near $[I]$ i.e. there exist an open $\G_m$-invariant neighbourhood $U \subseteq \cB^+$ of $s([I])$ such that $i_{|U} : U \to \cB$ is an open immersion. In particular, we identify $T_{[I]} \cB^+ = T_{[I]} \cB$.
    \end{corollary}
    \begin{proof}
        This follows from Corollary \ref{cor:nopositivedeg} and Proposition \ref{ref:comparison:intro:prop}.
    \end{proof}
    
    To study a neighbourhood of $[I] \in \cH$ and make use of obstruction theory on $\cH$ we need to analyse resolutions of $\Fperp$ and $I$. First, we take the minimal graded free resolution of the $S$-module $\Fperp$ which is of the following form, because $F$ is general enough (Definition \ref{def:nice})
        \[\begin{tikzcd}
        	0 & \Fperp & {S^{\oplus 15}(-2)} & {S^{\oplus 35}(-3)} & \dots
        	\arrow["d_1"', from=1-4, to=1-3]
        	\arrow["d_0"', from=1-3, to=1-2]
        	\arrow[from=1-2, to=1-1]
        	\arrow[from=1-5, to=1-4]
        \end{tikzcd}\]
        Here $d_0$ sends the natural basis $(E_i)_{1 \leq i \leq 15}$ of $S^{\oplus 15}(-2)$ to fifteen generators of $\Fperp$. Next we readjust this resolution so that it works for $I$
        \[\begin{tikzcd}
        	0 & I & {S^{\oplus 15}(-2) \oplus S(-3)} & {S^{\oplus 35}(-3) \oplus S^{\oplus 6}(-4)} & \dots
        	\arrow["{d_1}"', from=1-4, to=1-3]
        	\arrow["{d_0}"', from=1-3, to=1-2]
        	\arrow[from=1-2, to=1-1]
        	\arrow["{d_2}"', from=1-5, to=1-4]
        \end{tikzcd}\]
        Here for $G$ a generator of $S(-3)$, we define $d_0(G) = g$, and for $(H_k)_{1 \leq k \leq 6}$ - generators of $S^{\oplus 6}(-4)$, we put $d_1(H_k) = \sum_{i = 1}^{15} b_{ki} \cdot E_i + \alpha_k \cdot G$ for quadrics $b_{ki}$'s such that $\alpha_k \cdot g + \sum_{i = 1}^{15} b_{ki} \cdot d_0(E_i) = 0$. The existence of such quadrics follows from the fact that $(\Fperp)_4 = S_4$. The ideal $I$ is generated by $\Fperp$ and $g$, so $d_0$ is surjective. Moreover, for $f \in \ker d_0$ there exists $h$ in $S^{\oplus 6}(-4)$ such that $f - d_1(h) \in \ker d_0 \cap S^{\oplus 15}(-2)$. But $\ker d_0 \cap S^{\oplus 15}(-2) = d_1(S^{\oplus 35}(-3))$, because we have started from a resolution of $\Fperp$. Summarising, we got a partial resolution of $I$. We will use it to study tangent vectors at $[I] \in \cH$. First, let us define a subspace of $T_{[I]} \cH$ that appears in the description of $T_{[I]}V$.
        \begin{definition}\label{def:subspacew}
            Let $\nameX := \langle \symbtan_Q : Q \in (S/I)_2 \rangle \subseteq \Hom_{S}(I, S/I)$ be the subspace consisting of $\symbtan_Q$'s that map $(\Fperp)_2$ to $0$ and $g$ to $Q$.
            The dependence of $\nameX$ on $I$ is omitted in this notation. Note that $\symbtan_Q$ is an $S$-module homomorphism as $Q$ is annihilated by $\alpha_1, \dots, \alpha_6$. Fix $Q_1, \dots, Q_6 \in S$ such that $ev_F(Q_i) = x_i$. Such elements exist because $F$ is general enough. Their classes in $S/I$, which we will denote also by $Q_1, \dots, Q_6$, form a basis of $(S/I)_2$. It follows that $\symbtan_{Q_1}, \dots, \symbtan_{Q_6}$ is a basis of $\nameX$. More generally, if $Q \in (S/I)_2$, we fix its lift $Q \in S$ given by the same linear combination of $Q_1, \dots, Q_6$. We write $\symbtan_i$ for $\symbtan_{Q_i}$. 
        \end{definition}
        \begin{example}\label{ex:tangentsw}
            We will determine the element $[I'] \in  \{ \varphi \in \Mor \left( \Spec( \kk[\eps] ), \cH \right) \ | \ \varphi((\eps))= [I] \} = T_{[I]} \cH$ which corresponds to $\symbtan_Q \in \nameX$.
            The point $[I']$ is given by an ideal $I' \subseteq S[\eps]$. By Remark \ref{rem:tangentinverse} we have
            \[I' = (f - \symbtan_Q (f)' \cdot \eps : f \in I).\]
            Since $I$ is generated by $\Fperp$ and $g$ we get
            \[ I' = \Fperp \cdot S[\eps] + (g - \eps Q) \cdot S[\eps].\]
            The $\kk$-algebra $\kk[\eps]$ is local, so $\frac{S[\eps]}{I'}$ as a flat, finitely generated $\kk[\eps]$-module is free. Consider elements $(1, \alpha_1, \dots, \alpha_6, Q_1(\aalpha), \dots, Q_6(\aalpha))$. Their classes mod $(\eps)$ form a basis of $\frac{S}{I} = \frac{S[\eps]}{(I',\eps)}$, so by the Nakayama's lemma we get that $(1, \alpha_1, \dots, \alpha_6, Q_1(\aalpha), \dots, Q_6(\aalpha))$ is a basis of the $\kk[\eps]$-module $\frac{S[\eps]}{I'}$.
        \end{example}
        \begin{remark}\label{rem:identification_evaluation}
            The evaluation map allows us to define an isomorphism
            \begin{equation*}\label{eq:iso_evaluation}
                \lambda : \nameX \to \kk[x_1, \dots, x_6]_1
            \end{equation*}
            by sending $\symbtan_Q$ to $ev_F(Q) \in \kk[x_1, \dots, x_6]_1$.
        \end{remark}
        
        \begin{definition}\label{def:subspacey}
            Let $\nabla : \kk[x_1, \dots, x_6]_1 \to \kk[\partial_1, \dots, \partial_6]_1$ be the isomorphism sending $x_i$ to $\partial_i$. We define
            \begin{equation*}
                \nameY := \Big\langle \symbtan_Q - \frac{1}{13} \nabla(\lambda(\symbtan_Q)) : Q \in (S/I)_2 \Big\rangle.
            \end{equation*}
            In other words $\nameY = \langle \symbtan_i - \frac{1}{13} \partial_i : i=1, \dots, 6 \rangle$. We will write $\symboly_i$ for $\symbtan_i - \frac{1}{13} \partial_i$ so that $\nameY = \langle \symboly_i : i=1, \dots, 6 \rangle$.
        \end{definition}
        
        Now we are ready to prove the following geometric characterization of the subspace $Y$:
        \begin{proposition}\label{prop:subspacew}
            The subspace $Y$ is the tangent space to the fiber of $\pi : \cB^+ \to \cB^{\G_m}$ over $[I] \in \cB^{\G_m}$. In other words $\nameY = T_{[I]} V$ under the chain of inclusions
            \begin{equation*}
                T_{[I]} V \subseteq T_{[I]} \cB^+ = T_{[I]} \cB \subseteq T_{[I]} \cH = \Hom_{S}(I, S/I).
            \end{equation*}
        \end{proposition}
        \begin{proof}
            From Proposition \ref{prop:hilbandbar} we know that $\cH = \cB \times \G_a^6$, so $T_{[I]} \cH = T_{[I]} \cB \oplus T_{ \{0\} } \G_a^6$ and the additive action of $v \in \G_a^6(\kk)$ on $[I] = ([I], 0) \in \cB \times \G_a^6$ is given by $v + ([I], 0) = ([I], v)$. Thus the subspace $T_{ \{0\} } \G_a^6 \subseteq T_{[I]} \cH$ is the image of the differential of the $\G_a^6$-invariant morphism from $\G_a^6$ to $\cH$ that maps $0$ to $[I]$. 
            By Section \ref{sec:Hilbschm} we know explicit formulas for elements of this subspace, so we get
            \begin{equation}\label{eq:subspacew1}
                T_{[I]} \cH = T_{[I]} \cB \oplus \langle \partial_1, \dots, \partial_6 \rangle.
            \end{equation}
            By Corollary \ref{cor:bplusisb} we have $T_{[I]} \cB^+ = T_{[I]} \cB$. The morphisms $\pi, s$ giving the diagram
            \[\begin{tikzcd}
            	\cB^+ \\
            	{\cB^{\G_m}}
            	\arrow["\pi", from=1-1, to=2-1]
            	\arrow["s", shift left=1, bend left, from=2-1, to=1-1]
            \end{tikzcd}\]
            are $\G_m$-equivariant and induce a $\G_m$-invariant decomposition $T_{[I]} \cB^+ = \ker d\pi \oplus \im ds$. Note that $T_{[I]} \cB^{\G_m} \simeq \im ds = (T_{[I]} \cB)_0$. Thus, by $\G_m$-invariance $\ker d\pi = (T_{[I]} \cB)_{\neq 0} = (T_{[I]} \cB)_{< 0}$ where the last equality holds, because $T_{[I]} \cB$ doesn't have any positive degrees by Corollary \ref{cor:nopositivedeg}. By the definition of $V$ we see that $T_{[I]} V = \ker d\pi \subseteq T_{[I]} \cB$. Hence, overall we get that
            \begin{equation}\label{eq:subspacew2}
                T_{[I]} V = (T_{[I]} \cB)_{< 0} \subseteq T_{[I]} \cH.
            \end{equation}
            In order to determine $(T_{[I]} \cB)_{< 0}$ we calculate $(T_{[I]} \cH)_{< 0}$. Take an arbitrary $\delta \in \Hom_S(I, S/I)$ homogeneous of negative degree. It induces the following maps from the fixed resolution of $I$
            \[\begin{tikzcd}
            	0 & I & {S^{\oplus 15}(-2) \oplus S(-3)} & {S^{\oplus 35}(-3) \oplus S^{\oplus 6}(-4)} & \dots \\
            	&& {S/I},
            	\arrow[from=1-5, to=1-4]
            	\arrow["{d_1}"', from=1-4, to=1-3]
            	\arrow["{d_0}"', from=1-3, to=1-2]
            	\arrow[from=1-2, to=1-1]
            	\arrow["\delta"', from=1-2, to=2-3]
            	\arrow["{\Bar{\delta}}"', from=1-3, to=2-3]
            	\arrow["0", from=1-4, to=2-3]
            \end{tikzcd}\]
            where $\Bar{\delta}$ is the composition $\Bar{\delta} = \delta \circ d_0$. By restricting $\delta$ to $\Fperp$ and $\Bar{\delta}$ to $S^{\oplus 15}(-2)$ we get the diagram
            \[\begin{tikzcd}
            	0 & \Fperp & {S^{\oplus 15}(-2)} & {S^{\oplus 35}(-3)} & \dots \\
            	&& {S/I} \\
            	&& {S/\Fperp}
            	\arrow[from=1-5, to=1-4]
            	\arrow["{d_1}"', from=1-4, to=1-3]
            	\arrow["{d_0}"', from=1-3, to=1-2]
            	\arrow[from=1-2, to=1-1]
            	\arrow["{\Bar{\delta}}"', from=1-3, to=2-3]
            	\arrow["0", from=1-4, to=2-3]
            	\arrow["\mu", from=3-3, to=2-3]
            	\arrow["\delta"', from=1-2, to=2-3]
            	\arrow["{\delta'}"', bend right, dashed, from=1-2, to=3-3]
            	\arrow["0", bend left, dashed, from=1-4, to=3-3]
            \end{tikzcd}\]
            where $\mu$ is the quotient map. The homomorphism $\delta$ is of negative degree so $\Bar{\delta}$ sends generators $E_i$'s of $S^{\oplus 15}(-2)$ to $(S/I)_{\leq 1}$. However, $\mu$ is an isomorphism in degrees $\leq 2$. Hence, in fact $\Bar{\delta}$ factor through $S/\Fperp$ and the composition with $d_1$ is still $0$. Thus, there exists the unique $\delta'$ such that the diagram commutes. But this $\delta'$ belongs to $\Hom_S(\Fperp, S/\Fperp) \simeq T_{ [\Fperp] } \cH^{14}$ so it can only be a linear combination of derivatives $\partial_i$ (because $F$ is general enough so in particular $[\Fperp] \in \cH^{14}$ has TNT).
            
            Now we consider two cases:
            \begin{enumerate}
                \item $\delta$ is of degree $\leq -2$: Then $\delta'$ is also of degree $\leq -2$ so $\delta' = 0$, because partial derivatives are of degree $-1$ and $\delta'$ lies in their span. Thus, $\delta|_{(\Fperp)_2} = 0$ or equivalently $\bar{\delta}|_{S^{\oplus 15}(-2)} = 0$. Let $\bar{\delta}(G) = \delta(g) =: s \in (S/I)_{\leq 1}$. We know that $\bar{\delta} \circ d_1 = 0$, so
                \begin{equation*}
                    0 = \bar{\delta} \circ d_1 ( H_k ) = \bar{\delta} \left(\sum_{i = 1}^{15} b_{ki} \cdot E_i + \alpha_k \cdot G \right) = \bar{\delta} (\alpha_k \cdot G) = \alpha_k \cdot s \in (S/I)_{\leq 2}.
                \end{equation*}
                Note that $(S/I)_{\leq 2} \simeq (S/\Fperp)_{\leq 2}$, so an element from $(S/I)_{\leq 2}$ is $0$ if and only if its evaluation on $F$ is $0$. Hence, we get
                \begin{equation*}
                    ev_F( \alpha_k \cdot s) = \frac{\partial}{\partial x_k} ev_F(s) = 0 \textnormal{ for all $k=1, \dots, 6$}.
                \end{equation*}
                The polynomial $ev_F(s)$ is of degree $\geq 2$, because $s \in (S/I)_{\leq 1}$. But this means that $ev_F(s) = 0$ because it has all partial derivatives equal to $0$. Thus $s=0$ in $S/I$ and so in this case $\delta = 0$.
                
                \item $\delta$ is of degree $-1$: The homomorphism $\delta' \in \Hom_S(\Fperp, S/\Fperp)$ is a linear combination of partial derivatives in this case, because $[\Fperp] \in \cH^{14}$ has TNT. If $\partial := \sum_{i=1}^6 c_i \partial_i$ is defined using the same coefficients, but for derivatives $\partial_i \in \Hom_S(I, S/I)$, then $\delta - \partial$ is $0$ under the restriction to $(\Fperp)_2$. Thus $\delta - \partial$ maps $(\Fperp)_2 \mapsto 0$ and $g \mapsto s$ for some element $s$ of $(S/I)_2$, because both $\delta$ and $\partial$ are of degree $-1$. Hence $\delta - \partial \in \nameX$ and we see that $\Hom_S(I, S/I)_{-1} = \nameX \oplus \langle \partial_1, \dots, \partial_6 \rangle$, as we have presented $\delta$ as a sum of elements from $\langle \partial_1, \dots, \partial_6 \rangle$ and $\nameX$ in a unique way, namely $\delta = \partial + (\delta - \partial)$.
            \end{enumerate}
            Summarising, we have proved that
            \begin{equation*}
                (T_{[I]} \cH)_{< 0} = (T_{[I]} \cH)_{-1} = \nameX \oplus \langle \partial_1, \dots, \partial_6 \rangle.
            \end{equation*}
            By Definition \ref{def:subspacew} and Definition \ref{def:subspacey} we get that $\nameX \oplus \langle \partial_1, \dots, \partial_6 \rangle = \nameY \oplus \langle \partial_1, \dots, \partial_6 \rangle$ so we have
            \begin{equation}\label{eq:subspacew3}
                (T_{[I]} \cH)_{< 0} = (T_{[I]} \cH)_{-1} = \nameY \oplus \langle \partial_1, \dots, \partial_6 \rangle.
            \end{equation}
            We will now show that $\nameY$ lies in $T_{[I]} \cB$. 
            
            By Example \ref{ex:tangentsw} we have an explicit description of morphisms $\symbtan_Q$ as elements of $\cH( \Spec(\kk[\eps]) )$. We apply Theorem \ref{thm:HCaffine}. Fix $\symbtan_i \in \nameX$ and consider the composition starting with the morphism coming from $\symbtan_i$:
            \[\begin{tikzcd}
            	{\Spec(\kk[\eps])} & \cH & {{(\Affine^6)^{13} \goodquotient \mathbb{S}_{13}}} & {\Affine^6}
            	\arrow[from=1-2, to=1-3]
            	\arrow[from=1-1, to=1-2]
            	\arrow[from=1-3, to=1-4]
            \end{tikzcd}\]
            The composition, on the level of rings, gives a homomorphism $\kk[\aalpha] \to \kk[\eps]$ which is the composition of
            \[\begin{tikzcd}
            	{\kk[\aalpha]} & {(\kk[\aalpha] \otimes_\kk \dots \otimes_\kk \kk[\aalpha])^{\mathbb{S}_d}} & {\left(\frac{S[\eps]}{I'} \otimes_{\kk[\eps]} \dots \otimes_{\kk[\eps]} \frac{S[\eps]}{I'}\right)^{\mathbb{S}_d}} & {\kk[\eps]}
            	\arrow[from=1-2, to=1-3]
            	\arrow["LNm", from=1-3, to=1-4]
            	\arrow[from=1-1, to=1-2]
            \end{tikzcd}\]
            where $I' = \Fperp \cdot S[\eps] + (g - \eps Q_i)$ is the ideal from Example \ref{ex:tangentsw}. Recall that $Q_i$ were fixed elements of $S$ such that $ev_F(Q_i)= x_i$ for $i=1, \dots, 6$. By Theorem \ref{thm:HCaffine} this whole composition takes $\alpha_j$, for $j=1, \dots, 6$, to $\frac{1}{13} Tr_{\frac{S[\eps]}{I'}/\kk[\eps]}(\alpha_j)$. By Example \ref{ex:tangentsw} we can take $(1, \alpha_1, \dots, \alpha_6, Q_1, \dots, Q_6)$ as a basis of $\frac{S[\eps]}{I'}$ over $\kk[\eps]$. Let $A = (a_{lm})_{l,m =1}^{13}$ be the matrix representing the multiplication by $\alpha_j$ with respect to this basis. We analyse its diagonal coefficients:
            \begin{itemize}
                \item $\alpha_j \cdot 1 = \alpha_j$ so $a_{11} = 0$.
                \item $\alpha_j \cdot \alpha_k$ is a $\kk$-linear combination of polynomials $Q_i$ modulo $\Fperp$, so for $k=1, \dots 6$ we have $a_{k+1,k+1} = 0$.
                \item For $k \neq j$ the polynomial $\alpha_j \cdot Q_k$ lies is $\Fperp$ because $ev_F(\alpha_j \cdot Q_k) = ev_{x_k}(\alpha_j) = 0$. Hence $\alpha_j \cdot Q_k = 0$ in $\frac{S[\eps]}{I'}$ and here the diagonal coefficients are also zero.
                \item For $k=j$ the polynomial $\alpha_j \cdot Q_j$ is equal to $g$ modulo $\Fperp$, because $ev_F(\alpha_j \cdot Q_j) = ev_{x_j}(\alpha_j) = 1 = ev_F(g)$. But $g \equiv \eps Q_i (mod \ I')$ because $I' = \Fperp \cdot S[\eps] + (g - \eps Q_i)$. Thus, the coefficient near $Q_j$ is $\delta_{ij} \cdot \eps$ in this case.
            \end{itemize}
            Overall, $\frac{1}{13} Tr_{\frac{S[\eps]}{I'}/\kk[\eps]}(\alpha_j) = \frac{1}{13} \delta_{ij} \cdot \eps$, so the induced tangent vector in $T_{ \{0\} } \Affine^6$ is $\frac{1}{13} \partial_i$. Thus the image of $\symboly_i = \symbtan_i - \frac{1}{13} \partial_i$ under the differential $d Bar$ is $0$. This implies that $Y \subseteq \ker d Bar = T_{[I]} \cB$.
            
            Now we bring everything together. First, equations \eqref{eq:subspacew1}, \eqref{eq:subspacew3} imply that
            \begin{equation*}
                 (T_{[I]} \cH)_{<0} = (T_{[I]} \cB)_{<0} \oplus \langle \partial_1, \dots, \partial_6 \rangle = \nameY \oplus \langle \partial_1, \dots, \partial_6 \rangle.
            \end{equation*}
            But $\nameY \subseteq T_{[I]} \cB$ so also $\nameY \subseteq (T_{[I]} \cB)_{<0}$ as elements of $\nameY$ are of degree $-1$. Hence we get
            \begin{equation*}
                \nameY = (T_{[I]} \cB)_{<0} = T_{[I]} V,
            \end{equation*}
            where the last equality is equation \eqref{eq:subspacew2}.
        \end{proof}
        
        Since the tangent space $T_{[I]}V$ is calculated, we now pass to the analysis of the primary obstructions of $[I] \in \cH$. Recall that as in Example \ref{ex:hilbert} the functor $D_{\cH, [I]}$ pro-representable by $\hat{\cO}_{\cH, [I]}$ has an obstruction theory with obstruction space $Ob = \Ext^1_{S}(I, S/I)$.
        Since we want to determine primary obstructions to $V$, the following proposition will make the calculations easier.
        \begin{proposition}\label{prop:partial_vanish}
            Let $ob_0 : \Sym_2(T_{[I]} \cH) \to \Ext^1_{S}(I, S/I)$ be the primary obstruction map. Take $\partial \in \langle \partial_1, \dots, \partial_6 \rangle \subseteq T_{[I]} \cH$ and $\delta \in T_{[I]} \cH$. Then $ob_0(\partial \cdot \delta) = 0$.
        \end{proposition}
        \begin{proof}
            Suppose that $\delta$ corresponds to the morphism $\delta : \Spec\left(\frac{\kk[\eps]}{(\eps)^2}\right) \to \cH$. Denote by $\mu:\G_a^6 \times \cH \to  \cH$ the additive action on the Hilbert scheme from Example \ref{ex:additiveaction}. By composing
            \[\begin{tikzcd}
            	{\G_a^6 \times \Spec\left(\frac{\kk[\eps]}{(\eps)^2}\right)} & {\G_a^6 \times \cH} & \cH
            	\arrow["\mu"', from=1-2, to=1-3]
            	\arrow["{id \times \delta}"', from=1-1, to=1-2]
            \end{tikzcd}\]
            we get a $\G_a^6$-equivariant morphism which extends $\delta$. The tangent vector $\partial$ is in $\langle \partial_1, \dots, \partial_6 \rangle$ so it is induced by the $\G_a^6$-action on $\cH$. Take a morphism $\bar{\partial} : \Spec\left(\frac{\kk[\eps']}{(\eps')^2}\right) \to \G_a^6$ such that its composition with the orbit map of $[I]$ is $\partial : \Spec\left(\frac{\kk[\eps']}{(\eps')^2}\right) \to \cH$. Consider the composition $\eta$
            \[\begin{tikzcd}
            	{\Spec\left(\frac{\kk[\eps']}{(\eps')^2}\right) \times \Spec \left( \frac{\kk[\eps]}{(\eps)^2} \right)} & {\G_a^6 \times \Spec\left(\frac{\kk[\eps]}{(\eps)^2}\right)} & {\G_a^6 \times \cH} & \cH \\
            	{\Spec\left(\frac{\kk[\eps,\eps']}{(\eps^2, \eps'^2)}\right).}
            	\arrow["{id \times \delta}", from=1-2, to=1-3]
            	\arrow["{\bar{\partial} \times id}", from=1-1, to=1-2]
            	\arrow["{=}"{marking}, draw=none, from=2-1, to=1-1]
            	\arrow["\eta", bend right=10, from=2-1, to=1-4]
            	\arrow["\mu", from=1-3, to=1-4]
            \end{tikzcd}\]
            By its definition $\eta$ restricted to $\Spec\left(\frac{\kk[\eps']}{(\eps')^2}\right) \times \{0\} \cong \Spec\left(\frac{\kk[\eps']}{(\eps')^2}\right)$ is $\partial$ and $\eta$ restricted to $\{0\} \times \Spec\left(\frac{\kk[\eps]}{(\eps)^2}\right) \cong \Spec\left(\frac{\kk[\eps]}{(\eps)^2}\right)$ is $\delta$. Thus the following diagram commutes
            \[\begin{tikzcd}
            	{\Spec\left(\frac{\kk[\eps,\eps']}{(\eps^2, \eps'^2)}\right)} & \cH \\
            	{\Spec\left(\frac{\kk[\eps']}{(\eps')^2}\right) \amalg_{\Spec(\kk)} \Spec\left(\frac{\kk[\eps]}{(\eps)^2}\right)} \\
            	{\Spec\left(\frac{\kk[\eps,\eps']}{(\eps, \eps')^2}\right).}
            	\arrow["\eta", dashed, from=1-1, to=1-2]
            	\arrow[hook', from=2-1, to=1-1]
            	\arrow["{\partial \amalg \delta}"', from=2-1, to=1-2]
            	\arrow["{=}"{marking}, draw=none, from=3-1, to=2-1]
            \end{tikzcd}\]
            On the level of algebras this diagram is dual to the diagram in the statement of Lemma \ref{fact:primary_vanish}. Thus the existence of $\eta$ implies that $ob_0(\partial \cdot \delta)=0$ by Lemma \ref{fact:primary_vanish}.
        \end{proof}
        Fix $Q \in (S/I)_2$. For $\symbtan_Q : I \to S/I$ we will construct a lift of $\symbtan_Q$ to a homogeneous chain complex map as in the diagram from Theorem \ref{thm:hilb_obstruction}. The liftings are not unique but we choose them so that they are amenable to calculations.
        We have
        \[\begin{tikzcd}
    	& 0 & I & {S^{\oplus 15}(-2) \oplus S(-3)} & {S^{\oplus 35}(-3) \oplus S^{\oplus 6}(-4)} & F_3 \\
    	0 & {S/I} & {S} & {S^{\oplus 15}(-2) \oplus S(-3)} & {S^{\oplus 35}(-3) \oplus S^{\oplus 6}(-4)} & \dots
    	\arrow[from=1-3, to=1-2]
    	\arrow["{d_0}"', from=1-4, to=1-3]
    	\arrow["{d_1}"', from=1-5, to=1-4]
    	\arrow["{d_2}"', from=1-6, to=1-5]
    	\arrow["{d_2}"', from=2-6, to=2-5]
    	\arrow["{d_1}"', from=2-5, to=2-4]
    	\arrow["{d_0}"', from=2-4, to=2-3]
    	\arrow["q"', from=2-3, to=2-2]
    	\arrow[from=2-2, to=2-1]
    	\arrow["\symbtan_Q"'{pos=0.8}, from=1-3, to=2-2]
    	\arrow["{s_1(\symbtan_Q)}"'{pos=0.7}, shift right=1, shorten <=3pt, from=1-4, to=2-3]
    	\arrow["{s_2(\symbtan_Q)}"'{pos=0.9}, from=1-5, to=2-4]
    	\arrow["{s_3(\symbtan_Q)}"'{pos=0.9}, from=1-6, to=2-5]
        \end{tikzcd}\]
        The maps here are defined in the following way:
        \begin{itemize}
            \item We define $s_1(\symbtan_Q) : S^{\oplus 15}(-2) \oplus S(-3) \to S$ which takes $S^{\oplus 15}(-2)$ to $0$ and the generator $G$ of $S(-3)$ to $Q \in S$. Then it is obvious that $\symbtan_Q \circ d_0 = q \circ s_1(\symbtan_Q)$.
            
            \item The second map to define is $s_2(\symbtan_Q)$. We declare it to be $0$ on $S^{\oplus 35}(-3)$. Now it is enough to define it on the basis $(H_k)_{1 \leq k \leq 6}$ of $S^{\oplus 6}(-4)$. Recall that $(E_i)_{1 \leq i \leq 15}$ is the natural basis of $S^{\oplus 15}(-2)$. By the definition $s_1(\symbtan_Q)|_{S^{\oplus 15}(-2)} = 0$, so \[s_1(\symbtan_Q) \circ d_1 (H_k) = s_1(\symbtan_Q) \left(\sum_{1\leq i \leq 15} b_{ki} \cdot E_{i} + \alpha_k \cdot G\right) = s_1(\symbtan_Q) (\alpha_k \cdot G) = \alpha_k \cdot Q \in S\] for $b_{ki}$'s as in the description of $d_1$ and $k=1, \dots, 6$. Note that $a_k := ev_F(\alpha_k \cdot Q) \in \kk$ because $\alpha_k \cdot Q \in S_3$. Hence $ev_F(\alpha_k \cdot Q - a_k \cdot g) = 0$, so $\alpha_k \cdot Q - a_k \cdot g \in (\Fperp)_3$ and there exist linear forms $a_{ki}$ such that $\alpha_k \cdot Q - a_k \cdot g$ is the image of $\sum_{1 \leq i \leq 15} a_{ki} E_{i}$ by $d_0$. This is because $\Fperp$ is generated by $(\Fperp)_2$, as $F$ is general enough. We write $a_{ki} = a_{ki} (Q)$, since they depend on $Q$. Now we define $s_2(\symbtan_Q) (H_k) := \sum_{1 \leq i \leq 15} a_{ki} E_{i} + a_k G$ for $k=1, \dots, 6$. The constant $a_k$ also depends on $Q$, so we denote it $a_k(Q) = ev_F(\alpha_k \cdot Q)$.
        \end{itemize}
        
        Note that $\Ext^1_{S}(I, S/I) = \frac{\ker d_2^{*}}{\im d_1^{*}}$ where \[ d_1^{*} : \Hom_{S}(S^{\oplus 15}(-2) \oplus S(-3), S/I) \to \Hom_{S}(S^{\oplus 35}(-3) \oplus S^{\oplus 6}(-4), S/I) \] is composing with $d_1$ from the right and \[ d_2^{*} : \Hom_{S}(S^{\oplus 35}(-3) \oplus S^{\oplus 6}(-4), S/I) \to \Hom_{S}(F_3, S/I) \] is composing with $d_2$ from the right.
        
        It will turn out that for our application only a part of the obstruction space will be needed. More precisely, the $S$-module $\Hom_{S}(S^{\oplus 35}(-3) \oplus S^{\oplus 6}(-4), S/I)$ splits as \[ \Hom_{S}(S^{\oplus 35}(-3), S/I) \oplus \Hom_{S}(S^{\oplus 6}(-4), S/I). \] Let $i_2 : S^{\oplus 6}(-4) \to S^{\oplus 35}(-3) \oplus S^{\oplus 6}(-4)$ be the natural embedding. It induces \[ i_2^* : \Hom_{S}(S^{\oplus 35}(-3) \oplus S^{\oplus 6}(-4), S/I) \to \Hom_{S}(S^{\oplus 6}(-4), S/I). \] Consider the projection
        \begin{equation}\label{eq:projection}
            \Ext^1_{S}(I, S/I) = \frac{\ker d_2^{*}}{\im d_1^{*}} \to \frac{i_2^*(\ker d_2^*)}{i_2^*(\im d_1^{*})}.
        \end{equation}
        We have the following proposition, which will be crucial for proving $0$-dimensionality of the scheme $V$.
        \begin{proposition}\label{prop:obscalc}
                The primary obstruction map $ob_0 : \Sym_2(T_{[I]} \cH) \to \Ext^1_{S}(I, S/I)$ restricted to $\Sym_2(T_{[I]} V) \subseteq \Sym_2(T_{[I]} \cH)$ and composed with the projection \eqref{eq:projection} yields a $\kk$-linear map
                \begin{equation*}
                    \Omega : \Sym_2(T_{[I]} V) \to \frac{i_2^*(\ker d_2^*)}{i_2^*(\im d_1^{*})},
                \end{equation*}
                whose kernel is $\big\langle \frac{\partial F}{\partial x_1}(\symboly_1, \dots, \symboly_6), \dots, \frac{\partial F}{\partial x_6}(\symboly_1, \dots, \symboly_6) \big\rangle \subseteq \Sym_2(T_{[I]} V)$, where $\symboly_p = \symbtan_p - \frac{1}{13} \partial_p$ for $p=1, \dots, 6$.
        \end{proposition}
        \begin{proof}
            Fix $\nameD \in \Sym_2(T_{[I]} V)$ and write it in the form $\nameD = \sum_{1 \leq l, m \leq 6} d_{lm} \symboly_l \symboly_m$, $d_{lm} = d_{ml}$ (we can assume symmetry of $d_{lm}$ since the characteristic of the base field is not two). Using this equation we expand $\nameD$ and get
            \begin{equation*}
                \nameD = \sum_{1 \leq l, m \leq 6} d_{lm} \symboly_l \symboly_m = \sum_{1 \leq l, m \leq 6} d_{lm} \symbtan_l \symbtan_m + \sum_{1 \leq p \leq 6} \partial_p \cdot \delta_p
            \end{equation*}
            for some $\delta_p \in T_{[I]} \cH$. By $\kk$-linearity of $\Omega$ and Proposition \ref{prop:partial_vanish} we get that
            \begin{equation}\label{eq:omegacalc1}
                \Omega(\nameD) = \Omega \left( \sum_{1 \leq l, m \leq 6} d_{lm} \symbtan_l \symbtan_m \right).
            \end{equation}
            We will use the formula \eqref{eq:hilb_obstruction} from Theorem \ref{thm:hilb_obstruction} to determine the right hand side of this equation. First we take any $Q, Q' \in (S/I)_2$ and calculate
            \begin{equation*}
                q \circ s_1(\symbtan_Q) \circ s_2(\symbtan_{Q'}) (H_k) = q \circ s_1(\symbtan_Q) \left(\sum_{1 \leq i \leq 15} a_{ki}(Q') E_{i} + a_k(Q') G \right) = q (a_k(Q') \cdot Q) = a_k(Q') \cdot Q,
            \end{equation*}
            so:
            \begin{equation*}
                q \circ (s_1(\symbtan_Q) \circ s_2(\symbtan_{Q'}) + s_1(\symbtan_{Q'}) \circ s_2(\symbtan_{Q})) (H_k) = a_k(Q') \cdot Q + a_k(Q) \cdot Q',
            \end{equation*}
            for $k=1, \dots, 6$. If $Q = Q_j$, then by the definition $a_k(Q) = a_k(Q_j) = ev_F(\alpha_k \cdot Q_j) = ev_{x_j}(\alpha_k) = \delta_{kj}$. The third equality follows from the fact that $ev_F(Q_j) = x_j$ and from Definition \ref{def:duality}. Hence for $Q = Q_j, Q' = Q_l$ we have
            \begin{equation*}
                q \circ (s_1(\symbtan_{Q_j}) \circ s_2(\symbtan_{Q_l}) + s_1(\symbtan_{Q_l}) \circ s_2(\symbtan_{Q_j})) (H_k) = \delta_{kl} \cdot Q_j + \delta_{kj} \cdot Q_l,
            \end{equation*}
            for all $1 \leq j,k,l \leq 6$.
            
            By equation \eqref{eq:omegacalc1}, $\kk$-linearity of $\Omega$ and Theorem \ref{thm:hilb_obstruction} we get that $\Omega (\nameD) = \sum_{{1 \leq l, m \leq 6}} d_{lm} \Omega(\symbtan_l \symbtan_m)$ is a class (modulo $i_2^*(\im d_1^{*})$) of a homomorphism
            \begin{equation*}
                r : H_k \mapsto \sum_{1 \leq l, m \leq 6} d_{lm} (a_k(Q_l) \cdot Q_m + a_k(Q_m) \cdot Q_l),
            \end{equation*}
            but
            \begin{gather*}
                \sum_{1 \leq l, m \leq 6} d_{lm} (a_k(Q_l) \cdot Q_m + a_k(Q_m) \cdot Q_l) = \sum_{1 \leq l, m \leq 6} d_{lm} (\delta_{kl} \cdot Q_m + \delta_{km} \cdot Q_l) = \\
                \sum_{1 \leq m \leq 6} d_{km} Q_m + \sum_{1 \leq l \leq 6} d_{lk} Q_l = 2 \sum_{1 \leq l \leq 6} d_{lk} Q_l,
            \end{gather*}
            so $\nameD \in \ker \Omega$ if and only if the homomorphism $r : H_k \mapsto 2 \sum_{1 \leq l \leq 6} d_{lk} Q_l$ is in $i_2^*(\im d_1^{*})$.
            
            The condition $r \in i_2^*(\im d_1^*)$ is equivalent to the existence of a map $\nameq$ making the following diagram commute (with $r|_{S^{\oplus 35}(-3)} = 0$)
            \[\begin{tikzcd}
            	{S^{\oplus 15}(-2) \oplus S(-3)} & {S^{\oplus 35}(-3) \oplus S^{\oplus 6}(-4)} & {\sum_{i} b_{ki} \cdot E_{i} + \alpha_k \cdot G} & {H_k} \\
            	{S/I} && {2 \sum_l d_{lk} Q_l}
            	\arrow["{d_1}"', from=1-2, to=1-1]
            	\arrow["\nameq", dashed, from=1-1, to=2-1]
            	\arrow["r", from=1-2, to=2-1]
            	\arrow[maps to, from=1-4, to=1-3]
            	\arrow[dashed, maps to, from=1-3, to=2-3]
            	\arrow[maps to, from=1-4, to=2-3]
            \end{tikzcd}\]
            In other words
            \begin{equation*}
                \nameD \in \ker \Omega \iff r \in i_2^*(\im d_1^*) \iff \textnormal{there exists } \nameq \textnormal{ making the diagram commute.}
            \end{equation*}
            Such a map $\nameq$ is determined by images of $(E_{i})_{1 \leq i \leq 15}$ and $G$. Hence it exists if and only if there exist $(s_{i})_{1 \leq i \leq 15}$ and $s$ in $S/I$ such that $\nameq$ defined by $\nameq(E_{i}) = s_{i}$ and $\nameq(G) = s$ makes the diagram commute. This imposes the following equations on $(s_{i})_{1 \leq i \leq 15}$ and $s$
            \begin{equation*}\label{eq:proof_obs_1}
                \sum_{1 \leq i \leq 15} b_{ki} \cdot s_{i} + \alpha_k \cdot s =  2 \sum_{1 \leq l \leq 6} d_{lk} Q_l \textnormal{ for $k=1, \dots, 6$} \\
                \textnormal{ and } \nameq|_{d_1(S^{\oplus 35}(-3))} = 0.
            \end{equation*}
            Note that $2 \sum_{1 \leq l \leq 6} d_{lk} Q_l$ is of degree $2$ in $S/I$, so we can look only at $s_{i}$'s of degree $0$ (because by the definition $b_{ki}$'s are of degree two) and $s$ of degree 1. In this case $\nameq|_{S^{\oplus 15}(-2)}$ factors through $S/\Fperp$.
            Now, using the fact that $0 \leftarrow \Fperp \leftarrow S^{\oplus 15}(-2) \leftarrow S^{\oplus 35}(-3)$ is exact (because $F$ is general enough), we see that the condition $\nameq|_{d_1(S^{\oplus 35}(-3))} = 0$ implies that $\nameq|_{S^{\oplus 15}(-2)}$ factors through $\Fperp$. We get the following commutative diagram
            \[\begin{tikzcd}
            	0 & \Fperp & {S^{\oplus 15}(-2)} & S^{\oplus 35}(-3) \\
            	&& {S/I} \\
            	&& {S/\Fperp}
            	\arrow[from=1-4, to=1-3]
            	\arrow[from=1-3, to=1-2]
            	\arrow["\nameq", from=1-3, to=2-3]
            	\arrow["{r|_{S^{\oplus 35}(-3)} = 0}", from=1-4, to=2-3]
            	\arrow[from=1-2, to=2-3]
            	\arrow[from=1-2, to=1-1]
            	\arrow["{\hat{\nameq}}"', from=1-2, to=3-3]
            	\arrow[from=3-3, to=2-3]
            \end{tikzcd}\]
            where $\hat{\nameq}$ is the unique homomorphism respecting the gradings that lifts $F^{\perp} \to S/I$. Again, by the fact that $F$ is general enough we know that such $\hat{\nameq}$ must be $0$. Indeed, that is because $\Hom_S(\Fperp, S/\Fperp)_{-2} = 0$ as $[\Fperp] \in \cH^{14}$ has TNT (see Definition \ref{def:TNT}). Hence, the constants $s_i$ must be all $0$.
            
            This means that the existence of $\nameq$ satisfying our assumptions is equivalent to the existence of a linear polynomial $s$ such that for $k = 1, \dots, 6$
            \begin{equation}\label{eq:proof_obs_2}
                \alpha_k \cdot s =  2 \sum_{1 \leq l \leq 6} d_{lk} Q_l.
            \end{equation}
            Note that the equality \eqref{eq:proof_obs_2} takes place in $(S/I)_2 = (S/ \Fperp)_2$, so we can apply the isomorphism $ev_F : S/\Fperp \to ev_F(S)$ and get
            \begin{equation}\label{eq:proof_obs_3}
                ev_F(\alpha_k \cdot s) =  2 \sum_{1 \leq l \leq 6} d_{lk} x_l \textnormal{ for $k=1, \dots, 6$.}
            \end{equation}
            Moreover, if we find $s$ such that equation \eqref{eq:proof_obs_3} holds, then also equation \eqref{eq:proof_obs_2} is satisfied by it, as the evaluation function is an isomorphism on the second degree of $S/ \Fperp$.
            Equation \eqref{eq:proof_obs_3} can be rewritten as
            \begin{equation*}\label{eq:proof_obs_4}
                \frac{\partial}{\partial x_k} (ev_F(s)) = \frac{\partial}{\partial x_k} \left( \sum_{1 \leq l,m \leq 6} d_{lm} x_l x_m \right) \textnormal{ for $k=1, \dots, 6$}.
            \end{equation*}
            and rearranging the terms, we get
            \begin{equation*}\label{eq:proof_obs_5}
                \frac{\partial}{\partial x_k} \left( \sum_{1 \leq l,m \leq 6} d_{lm} x_l x_m - ev_F(s) \right) = 0 \textnormal{ for $k=1, \dots, 6$}.
            \end{equation*}
            The quadric $\sum_{1 \leq l,m \leq 6} d_{lm} x_l x_m - ev_F(s)$ has all partial derivatives $0$ if and only if it is $0$ itself. Thus $\nameq$ exists if and only if $\sum_{1 \leq l,m \leq 6} d_{lm} x_l x_m = ev_F(s)$ for some $s \in S_1$ if and only if $\sum_{1 \leq l,m \leq 6} d_{lm} x_l x_m \in ev_F(S_1) = \langle \frac{\partial F}{\partial x_1}, \dots, \frac{\partial F}{\partial x_6} \rangle$. 
            
            The existence of $\nameq$ is equivalent to the fact that $\nameD \in \ker \Omega$, so we get that
            \begin{equation*}\label{eq:proof_obs_6}
                \ker \Omega = \bigg\langle \frac{\partial F}{\partial x_1}(\symboly_1, \dots, \symboly_6), \dots, \frac{\partial F}{\partial x_6}(\symboly_1, \dots, \symboly_6) \bigg\rangle,
            \end{equation*}
            which concludes the proof.
        \end{proof}

        \subsection{Zero-dimensionality of $V$}\label{sub:descriptionofv}
        We will use Proposition \ref{prop:obscalc} to prove that $V$ is $0$-dimensional.
        Let us start with some analysis of $V$. Recall that it is defined as $V=\pi^{-1}([I])$ for $\pi : \cB^+ \to \cB^{\G_m}$ and that it is affine with $\Affine^1$-action (by multiplication) where the unique fixed point is $[I]$, see Theorem \ref{ref:introRepresentability:thm}. Since it is affine let $V = \Spec B$. The action of $\Affine^1$ on $\Spec B$ translates to $\mathbb{N}$-grading on $B$ so $B = \bigoplus_{n \geq 0} B_n$ with $B_0 = \kk$ as we have only one fixed point. Moreover, we know that $T_{[I]}V$ is $6$-dimensional (Proposition \ref{prop:subspacew} and Definition \ref{def:subspacew}), so for $\mm := \bigoplus_{n \geq 1} B_n$ which is the fix-point ideal of $B$, we get that $\dim_{\kk} \frac{\mm}{\mm^2}=6$. Fix six homogeneous representatives of generators of $\frac{\mm}{\mm^2}$ and call them $\symboly_1^{\vee}, \dots, \symboly_6^{\vee} \in \mm$. Assume they are dual to generators $\symboly_1, \dots, \symboly_6$ of the tangent space at $[I]$. Define a graded homomorphism $\kk[\symb_1, \dots, \symb_6] \to B$ sending $\symb_i$ to $\symboly_i^{\vee}$. By \cite[Lemma 0EKB]{Stacks} it is onto.
        If we denote by $J$ the kernel of this homomorphism, then we get an identification of graded rings $B \simeq \kk[\symb_1, \dots, \symb_6]/J$ with $J$ - homogeneous. Note that also $J \subseteq (\symb_1, \dots, \symb_6)^2$. So to prove that $V$ is $0$-dimensional it is enough to prove that $\kk[\symb_1, \dots, \symb_6]/J$ is $0$-dimensional. We will do so by obtaining a lower bound on the ideal $J$.
        
        First consider the following composition of inclusions and morphisms
        \begin{equation*}
            V \subseteq \cB^+ \to \cB \subseteq \cH.
        \end{equation*}
        By looking near $[I]$ we get a homomorphism $\cO_{\cH, [I]} \to \cO_{V, [I]}$. It is surjective as $V \subseteq \cB^+$ is a closed inclusion, $\cB^+ \to \cB$ is an isomorphism near $[I]$ and $\cB \subseteq \cH$ is a closed inclusion.
        We complete it and get a map $\zeta : \hat{\cO}_{\cH, [I]} \to \hat{\cO}_{V, [I]}$ which is also surjective.
        Using the map $\zeta$ we will prove the following.
        \begin{proposition}\label{prop:boundj}
            The vector space \[(\Fperp)_2(\Bar{\symb}) := \{ g(\symb_1, \dots, \symb_6) \in \kk[\symb_1, \dots, \symb_6] : g(\alpha_1, \dots, \alpha_6) \in (\Fperp)_2 \}\] is contained in $J_2$.
        \end{proposition}
        Before giving the proof, we give an immediate corollary.
        \begin{corollary}\label{cor:zerodim}
            $V$ is $0$-dimensional.
        \end{corollary}
        \begin{proof}
            The cubic $F$ is general enough, so $\Fperp$ is generated by its second degree $(\Fperp)_2$. By Proposition \ref{prop:boundj} $(\Fperp)_2(\Bar{\symb}) \subseteq J$, so also $\Fperp(\Bar{\symb}) \subseteq J$. Hence we get a surjection
            \[\begin{tikzcd}
            	{\frac{\kk[\Bar{\symb}]}{\Fperp(\Bar{\symb})}} & {\frac{\kk[\Bar{\symb}]}{J}.}
            	\arrow[two heads, from=1-1, to=1-2]
            \end{tikzcd}\]
            The left hand side ring is $0$-dimensional, so the other one also is. Since $V \simeq \Spec(\frac{\kk[\Bar{\symb}]}{J})$ we get the conclusion.
        \end{proof}
        \begin{proof}[Proof of the Proposition \ref{prop:boundj}]
            Denote $Ob = \Ext^1_{S}(I, S/I)$ and let $ob_0 : \Sym_2(T_{[I]} \cH) \to Ob$ be the primary obstruction map for $\hat{\cO}_{\cH, [I]}$.
            Let $ob_V : \Sym_2(T_{[I]} V) \to Ob$ be the restriction of $ob_0$ to $\Sym_2(T_{[I]} V)$.
            By Proposition \ref{prop:primaryobs_tool} with $\zeta : \hat{\cO}_{\cH, [I]} \to \hat{\cO}_{V, [I]}$ we get
            \[ \im ob_V^{\vee} \subseteq \frac{(\Bar{\symb})^3 + J}{(\Bar{\symb})^3} = \frac{(\Bar{\symb})^3 + J_2}{(\Bar{\symb})^3}.\] 
            Hence, it is enough to prove that $(\Fperp)_2(\Bar{\symb}) \subseteq \im ob_V^{\vee}$.
            Note that if we compose $ob_V$ from the right with $Ob \to \frac{i_2^*(\ker d_2^*)}{i_2^*(\im d_1^{*})}$ we get the map $\Omega$ from Proposition \ref{prop:obscalc}. Thus we end up with the sequence
            \[\begin{tikzcd}
            	{\ker \Omega} & {\Sym_2(T_{[I]} V)} & Ob & {\frac{i_2^*(\ker d_2^*)}{i_2^*(\im d_1^{*})},}
            	\arrow["j"', from=1-1, to=1-2]
            	\arrow["{ob_V}"', from=1-2, to=1-3]
            	\arrow[from=1-3, to=1-4]
            	\arrow["\Omega"', bend left, from=1-2, to=1-4]
            \end{tikzcd}\]
            so after dualizing we get
            \[\begin{tikzcd}
            	{(\ker \Omega)^{\vee}} & {{\Sym_2(T_{[I]} V)}^{\vee}} & {Ob^{\vee}} & {\left(\frac{i_2^*(\ker d_2^*)}{i_2^*(\im d_1^{*})}\right)^{\vee}.}
            	\arrow["j^{\vee}", from=1-2, to=1-1]
            	\arrow["{ob_V^{\vee}}", from=1-3, to=1-2]
            	\arrow[from=1-4, to=1-3]
            	\arrow["{\Omega^{\vee}}", bend right, from=1-4, to=1-2]
            \end{tikzcd}\]
            Hence $\im \Omega^{\vee} \subseteq \im ob_V^{\vee} \subseteq J_2$. Note that $j$ is the kernel of $\Omega$, so $\im \Omega^{\vee} = \ker j^{\vee}$. 
            
            We will prove that $\ker j^{\vee} = (\Fperp)_2(\Bar{\symb}) \subseteq \Sym_2(T_{[I]}^{\vee} V)$ under the identification $\Sym_2(T_{[I]}^{\vee} V) \simeq \Sym_2(T_{[I]} V)^{\vee}$ coming from the natural duality $\cdot : \Sym_2(W^{\vee}) \times \Sym_2(W) \to \kk$ for $W = T_{[I]} V$ (see Equation \eqref{eq:dualityw}). However, for the proof it will be easier to start with a different duality. First let us write $\Sym_2(T_{[I]} V) = \kk[\symboly_1, \dots, \symboly_6]_2$. Consider the evaluation map from Definition \ref{def:duality} with variables $x_i$ changed to $\symboly_i$
            \begin{equation*}
                \ \wcirc \ : \kk[\alpha_1, \dots, \alpha_6] \times \kk[\symboly_1, \dots, \symboly_6] \to \kk[\symboly_1, \dots, \symboly_6].
            \end{equation*}
            In other words $\alpha_i$ acts as $\frac{\partial}{\partial \symboly_i}$. By restricting this map to degree two we get a perfect pairing
            \begin{equation*}
                \ \wcirc \ : \kk[\alpha_1, \dots, \alpha_6]_2 \times \kk[\symboly_1, \dots, \symboly_6]_2 \to \kk.
            \end{equation*}
            This gives us an identification $\kk[\alpha_1, \dots, \alpha_6]_2 \simeq (\kk[\symboly_1, \dots, \symboly_6]_2)^{\vee} \simeq \Sym_2(T_{[I]} V)^{\vee}$.
            We will now calculate $\ker j^{\vee} \subseteq \Sym_2(T_{[I]} V)^{\vee} \simeq \kk[\alpha_1, \dots, \alpha_6]_2$. Take $\nameD \in \kk[\alpha_1, \dots, \alpha_6]_2$. We have
            \begin{equation*}
                \nameD \in \ker j^{\vee} \textnormal{ if and only if } \nameD|_{\ker \Omega} = 0.
            \end{equation*}
            From Proposition \ref{prop:obscalc} we know that $\ker \Omega = \big\langle \frac{\partial F}{\partial x_1}(\symboly_1, \dots, \symboly_6), \dots, \frac{\partial F}{\partial x_6}(\symboly_1, \dots, \symboly_6) \big\rangle$, so
            \begin{equation*}
                \nameD|_{\ker \Omega} = 0 \iff  \nameD \ \wcirc \ \left(\frac{\partial F}{\partial x_k}(\symboly_1, \dots, \symboly_6)\right)=0 \textnormal{ for all $k=1, \dots, 6$.}
            \end{equation*}
            We transform this condition using the evaluation action
            \[\nameD \ \wcirc \ \left(\frac{\partial F}{\partial x_k}(\symboly_1, \dots, \symboly_6)\right) = (\nameD \cdot \alpha_k) \ \wcirc \ (F(\symboly_1, \dots, \symboly_6)) = (\alpha_k \cdot \nameD) \ \wcirc \ (F(\symboly_1, \dots, \symboly_6)) = \frac{\partial}{\partial \symboly_k}(\nameD \ \wcirc \ F).\]
            Since $\nameD \ \wcirc \ F$ is a linear form, the vanishing of $\frac{\partial}{\partial \symboly_k}(\nameD \ \wcirc \ F)$ for all $k=1, \dots 6$ is equivalent to $\nameD \ \wcirc \ F = 0$. Thus we get
            \begin{equation*}
                \nameD \in \ker j^{\vee} \iff \nameD|_{\ker \Omega} = 0 \iff \nameD \ \wcirc \ F = 0.
            \end{equation*}
            Hence $\ker j^{\vee} = (\Fperp)_2 \subseteq \kk[\alpha_1, \dots, \alpha_6]_2$. The only thing left is to check what happens under the isomorphisms
            \begin{equation}\label{eq:isomorphism}
                \kk[\alpha_1, \dots, \alpha_6]_2 \simeq \Sym_2(T_{[I]} V)^{\vee} \simeq \Sym_2(T_{[I]}^{\vee} V) = \Sym_2(\kk \symb_1 \oplus \dots \oplus \kk \symb_6).
            \end{equation}
            Take $\alpha_i \alpha_j \in \kk[\alpha_1, \dots, \alpha_6]_2$ and $\symboly_i^{\vee} \symboly_j^{\vee} \in \Sym_2(T_{[I]}^{\vee} V)$ where $(\symboly_i^{\vee})_{i=1}^6$ is the dual basis to $(\symboly_i)_{i=1}^6 \subseteq T_{[I]} V$.
            Consider the action of these elements on $\Sym_2(T_{[I]} V)$:
            \begin{itemize}
                \item $\alpha_i \alpha_j \cdot \symboly_i \symboly_j$ is $1$ if $i \neq j$ and $2$ if $i=j$, 
                \item $\symboly_i^{\vee} \symboly_j^{\vee} \cdot \symboly_i \symboly_j = \frac{1}{2}(\symboly_i^{\vee}(\symboly_i)\symboly_j^{\vee}(\symboly_j) +  \symboly_i^{\vee}(\symboly_j)\symboly_j^{\vee}(\symboly_i))$ so it is $\frac{1}{2}$ if $i \neq j$ and $1$ if $i=j$ (in both cases we use Formula \eqref{eq:dualityw}).
             \end{itemize}
            From a calculation it follows that the action on $\symboly_a \symboly_b$ for $\{a, b\} \neq \{i, j\}$ is zero in both cases. Thus under the isomorphism (\ref{eq:isomorphism}), $\alpha_i \alpha_j \simeq 2 \symboly_i^{\vee} \symboly_j^{\vee}$.
            Hence, after two isomorphisms everything is multiplied by $2$, when considered in bases $(\alpha_i \alpha_j)$ and $(\symboly_i^{\vee} \symboly_j^{\vee})$. Now, by the definition of $\symb_i$'s (they were coming exactly from the dual basis to $(\symboly_i)_i$) and homogeneity of $(\Fperp)_2$, we get the conclusion.
        \end{proof}
        
        \begin{corollary}\label{cor:Iishedgehog}
            The ideal $I = (\Fperp, g) = \Fperp + S_{\geq 3}$ yields a hedgehog point $[I] \in \cB$. In particular $[I]$ is a non-reduced point of $\cB$ and $\cH$.
        \end{corollary}
        \begin{proof}
            The point $[I]$ lies in $\cB$ as it is supported at $0 \in \Affine^6$. It is a $\G_m$-fixed point, because $I$ is homogeneous.
            Now we check three properties from Definition \ref{def:hedgehog:precise}:
            \begin{enumerate}
                \item $(T_{[I]} \cB)_{>0} = 0$ by Corollary \ref{cor:nopositivedeg},
                \item $T_{[I]} \cB \neq (T_{[I]} \cB)_0$ by Proposition \ref{prop:subspacew} and Definition \ref{def:subspacew},
                \item the negative spike $V$ at $[I]$ is $0$-dimensional by Corollary \ref{cor:zerodim}.
            \end{enumerate}
            The non-reducedness follows from Theorem \ref{thm:Hedgehog_point_theorem:precise}.
        \end{proof}
        Hence, by Example \ref{ex:nonempty} we get the following.
        \begin{corollary}\label{cor:important}
            The ideal $I = (\Fperp, g) = \Fperp + S_{\geq 3}$ for $F = x_1 x_2 x_4 - x_1 x_5^2 + x_2 x_3^2 + x_3 x_5 x_6 + x_4 x_6^2$ yields a non-reduced point $[I] \in \cH$.
        \end{corollary}
        
    \subsection{Fractal family and its flatness}\label{sub:flatness}
    Since we proved that $V$ is $0$-dimensional in Corollary \ref{cor:zerodim}, we know that the composition of morphisms
    \begin{equation}\label{eq:compo}
        V \subseteq \cB^+ \to \cB \subseteq \cH
    \end{equation}
    is a closed embedding. Indeed, that is because $V$ is concentrated on one point, the first and the last morphisms are closed embeddings, and $i : \cB^+ \to \cB$ is an isomorphism on some open neighbourhood of $[I] \in \cB^+$. The subscheme $V = \Spec(B) \hookrightarrow \cH$ yields a deformation of $[I]$ over $\Spec(B)$. In Subsection \ref{sub:descriptionofv} we presented $B$ as a quotient $\frac{\kk[\symb_1, \dots, \symb_6]}{J}$, where $J$ was a homogeneous ideal. Moreover, we proved that $\Fperp(\Bar{\symb}) \subseteq J$ in Proposition \ref{prop:boundj} and this was the heart of our argument for $0$-dimensionality of $V$ in Corollary \ref{cor:zerodim}. Thus, a natural suspicion is that in fact $J = \Fperp(\Bar{\symb})$, or more generally $V \simeq \Spec(S/\Fperp)$. If this was true, then there would exist a morphism $\Spec(S/\Fperp) \to \cH$ onto $V$. However, such a morphism would induce a deformation of $S/I$ over $S/\Fperp$ - almost over itself! Moreover, since it would factor through $\cB^+$ we should have been able to prolong it to an element of $\cH(\Affine^1 \times \Spec(S/\Fperp))$. Suspecting the existence of such a structure, we go the other way around: we find a deformation over $\Affine^1 \times \Spec(S/\Fperp)$ which yields a morphism to $\cH$ such that after composing it with the projection $\cH = \G_a^6 \times \cB \to \cB$ the induced morphism to $\cB^+$ is an isomorphism onto $V$. In this subsection we deal with the commutative algebra part of the problem, namely we construct a particular family, called the fractal family (that depends on a cubic and on some other choices), and prove that it is flat in Corollary \ref{cor:flatness_of_fractal_family}.
    
    In order to define the fractal family we need a relative version of the Macaulay duality. Although what we do in the following part of this subsection could be done in more generality, we decided to work in six variables for concreteness, because anyhow we are interested in general enough cubics (see Definition \ref{def:nice}), and because the ultimate goal is to prove Theorem \ref{thm:intro_completeV} (which is Theorem \ref{thm:completeV} later in the text). 
    
\begin{notation}\label{notation:relative_Macaulay}
        We define rings of new variables $\kk[\beta_1, \dots, \beta_6]$ and $\kk[y_1, \dots, y_6]$ so that $\beta_i := \frac{\partial}{\partial y_i}$. Additionally, we define a new free variable $t$ which will correspond to the torus action.
        In this case we get the relative evaluation map
        \begin{equation*}
            ev : \kk[t][\alpha_1, \dots, \alpha_6, \beta_1, \dots, \beta_6] \times \kk[t][x_1, \dots, x_6, y_1, \dots, y_6] \to \kk[t][x_1, \dots, x_6, y_1, \dots, y_6].
        \end{equation*}
        The action of $\kk[t, \Bar{\alpha}, \Bar{\beta}]$ on $\kk[t, \Bar{x}, \Bar{y}]$ will be also denoted by $\wcirc$ so that $ev((a, b)) = a \ \wcirc \ b$.  Moreover, if $b \in \kk[t, \Bar{x}, \Bar{y}]$ then $ev_{b} : \kk[t, \Bar{\alpha}, \Bar{\beta}] \to \kk[t, \Bar{x}, \Bar{y}]$ is the map $ev_b(a) := a \ \wcirc \ b$.
        
        Fix a cubic $F$ in six variables $x_1, \dots, x_6$. We assume that there are polynomials 
        $g, Q_1, \dots, Q_6 \in \kk[\aalpha]$ such that $g \wcirc F = 1$ and $Q_i \wcirc F = x_i$ for $i=1, \dots, 6$, with $g$ of degree three, and $Q_i$'s of degree two. We write $F_{\Bar{x}} := F(\Bar{x}), F_{\Bar{y}} := F(\Bar{y})$ (so these are the same cubics, in different variables) and more generally if $G$ is a polynomial in six variables, we write $G(\xx)$ (or $G(\yy), G(\aalpha), G(\xx+\yy)$) for the same polynomial, but with different names of variables. In particular, we use polynomials $F_{\Bar{x}} \cdot F_{\Bar{y}}$ or 
        \[F(\Bar{x} + \Bar{y}) :=F(x_1+y_1, \dots, x_6+y_6), \ F(t \cdot \Bar{x} + \Bar{y}) := F(t \cdot x_1+y_1, \dots, t \cdot x_6+y_6)\] of mixed variables, and in such cases we will usually omit in which ring we are calculating those, unless it is unclear. In particular, if we write for example $\frac{\kk[t, \aalpha]}{F_{\Bar{x}}^{\perp}}$, we mean the algebra $\kk[t, \aalpha]$ divided by the kernel of the restricted evaluation map 
        $ev_{F_{\Bar{x}}}:\kk[t, \aalpha] \to \kk[t,\xx]$. Note that if we treat $F_{\Bar{x}}^{\perp}$ as an ideal in $\kk[\aalpha]$, then the ideal it generates in $\kk[t,\aalpha]$ (or $\kk[t, \aalpha, \bbeta]$) coincides with the kernel of $ev_{F_{\Bar{x}}}:\kk[t, \aalpha] \to \kk[t,\xx]$ (resp. the kernel of $ev_{F_{\Bar{x}}}:\kk[t, \aalpha, \bbeta] \to \kk[t,\xx, \yy]$), which justifies this notation.
    \end{notation}
    
    \begin{definition}\label{def:gammat}
        Suppose $g, Q_1, \dots, Q_6 \in \kk[\aalpha]$ are as in Notation \ref{notation:relative_Macaulay}.
        Let $\Gamma(t) \in \kk[t, \aalpha, \bbeta]$ be defined by
        \begin{equation*}
            \Gamma(t) := g(\aalpha) + t \cdot \sum_{1 \leq i \leq 6} \beta_i \cdot Q_i(\aalpha) + t^2 \cdot \sum_{1 \leq i \leq 6} \alpha_i \cdot Q_i(\bbeta) + t^3 \cdot g(\bbeta).
        \end{equation*}
        Note that $\Gamma(t)$ depends on $g, F$ and $Q_i$'s even though it is not indicated in the notation.
    \end{definition}
    
    Now we can define the fractal family.
    
    \begin{definition}
        Consider a $\frac{\kk[t, \bbeta]}{\Fperpy}$ - algebra $M$ of the form
        \[\begin{tikzcd}
        	{\frac{\kk[t, \bbeta]}{\Fperpy}} & {M := \frac{\kk[t, \aalpha, \bbeta]}{(\Fperpx, \Fperpy, \Gamma(t))},}
        	\arrow[from=1-1, to=1-2]
        \end{tikzcd}\]
        where we take $\Gamma(t)$ from Definition \ref{def:gammat} with $g$ and $Q_i$'s as in Notation \ref{notation:relative_Macaulay}.
        We call it the \emph{fractal family}. The name refers to the self-similarity of the base space and the special fiber. The rest of this subsection aims to prove flatness of the fractal family.
    \end{definition}
    
    \begin{lemma}\label{fact:Gamma}
        The following identity holds
        \begin{equation*}
            \Gamma(t) \ \wcirc \ (F_{\Bar{x}} \cdot F_{\Bar{y}}) = F(t \cdot \xx + \yy).
        \end{equation*}
    \end{lemma}
    \begin{proof}
        The idea is to use Taylor's formula. First note that
        \begin{equation*}
            F(\xx + \yy) = F(\xx) + \sum_{1 \leq i \leq 6} \frac{\partial F}{\partial x_i}(\xx) \cdot y_i + \sum_{1 \leq i,j \leq 6} \frac{\partial^2 F}{\partial x_i \partial x_j}(\xx) \cdot y_i y_j + \dots
        \end{equation*}
        and also
        \begin{equation*}
            F(\xx + \yy) = F(\yy) + \sum_{1 \leq i \leq 6} \frac{\partial F}{\partial y_i}(\yy) \cdot x_i + \sum_{1 \leq i,j \leq 6} \frac{\partial^2 F}{\partial y_i \partial y_j}(\yy) \cdot x_i x_j + \dots
        \end{equation*}
        On the other hand, $F$ is a cubic, so $F(\xx + \yy)$ can be thought of as a cubic with variable $\xx$ and coefficients in $\kk[\yy]$. From the second Taylor expansion we see that the part of $\overline{x}$-degree zero is $F(\yy)$ and the part of $\overline{x}$-degree one is $\sum_{1 \leq i \leq 6} \frac{\partial F}{\partial y_i}(\yy) \cdot x_i$. Moreover, using the first Taylor expansion we get that the part of $F(\xx + \yy)$ of $\overline{x}$-degree at least two is $F(\xx) + \sum_{1 \leq i \leq 6} \frac{\partial F}{\partial x_i}(\xx) \cdot y_i$. Thus in fact
        \begin{equation*}
            F(\xx + \yy) = F(\xx) + \sum_{1 \leq i \leq 6} \frac{\partial F}{\partial x_i}(\xx) \cdot y_i + \sum_{1 \leq i \leq 6} \frac{\partial F}{\partial y_i}(\yy) \cdot x_i + F(\yy),
        \end{equation*}
        so also
        \begin{equation*}
            F(t \cdot \xx + \yy) = t^3 \cdot F(\xx) + t^2 \cdot \sum_{1 \leq i \leq 6} \frac{\partial F}{\partial x_i}(\xx) \cdot y_i + t \cdot \sum_{1 \leq i \leq 6} \frac{\partial F}{\partial y_i}(\yy) \cdot x_i + F(\yy).
        \end{equation*}
        By the assumptions about $g$'s and $Q_i$'s action on $F$ from Definition \ref{def:gammat} we get the conclusion.
    \end{proof}
    The following observation will be useful:
        \begin{lemma}\label{fact:fracfamilyform}
        Consider $(F_{\xx} \cdot F_{\yy})^{\perp} \subseteq \kk[t, \aalpha, \bbeta]$ and $\Fperpx \subseteq \kk[\aalpha], \Fperpy \subseteq \kk[\bbeta]$. Then
        \begin{equation*}
        (F_{\xx} \cdot F_{\yy})^{\perp} = (\Fperpx, \Fperpy)
        \end{equation*}
        as ideals in $\kk[t, \aalpha, \bbeta]$.
    \end{lemma}
    \begin{proof}
        Consider $P = \sum_{I,J,k} a_{I,J,k} t^k \aalpha^I \bbeta^J \in \kk[t, \aalpha, \bbeta]$.
        Then 
        \begin{equation}\label{equation:applying derivatives on both}
            ev_{F_{\xx} \cdot F_{\yy}}(P) = \sum_{I,J,k} a_{I,J,k} t^k \frac{\partial^{I+J}(F_{\xx} \cdot F_{\yy})}{\partial \xx^I \partial \yy^J} = \sum_{I,J,k} a_{I,J,k} t^k \frac{\partial^I F_{\xx}}{\partial \xx^I} \cdot \frac{\partial^J F_{\yy}}{\partial \yy^J}.
        \end{equation}
        Assume that $P \in (F_{\xx} \cdot F_{\yy})^{\perp}$.
        Then, for each natural $k_0 \in \NN$ appearing in the sum defining $P$, we have
        \begin{equation}\label{equation:applying P is zero}
            \sum_{I,J} a_{I,J,k_0} \frac{\partial^I F_{\xx}}{\partial \xx^I} \cdot \frac{\partial^J F_{\yy}}{\partial \yy^J} = 0.
        \end{equation}
        We will now show that for any $k_0$, the polynomial $P_0 := \sum_{I,J} a_{I,J,k_0} \aalpha^I \bbeta^J$ belongs to $(\Fperpx, \Fperpy)$, and thus also $P$ belongs to $(\Fperpx, \Fperpy)$.
        Since $\kk[\xx,\yy] \cong \kk[\xx] \otimes_{\kk} \kk[\yy]$, Equation (\ref{equation:applying P is zero}) holds if and only if there are finitely many elements $b_{p,I}, b_J', c_{q,J}, c_I'$ in $\kk$, where $p$ and $q$ vary among some finite sets of indices, such that for all $I, J$
        \[ a_{I,J,k_0} = \sum_{p} b_{p,I} \cdot b_J' + \sum_q c_{q,J} \cdot c_I' \]
        and for all $p,q$
        \[ \sum_I b_{p,I} \frac{\partial^I F_{\xx}}{\partial \xx^I}  = \sum_J c_{q,J} \frac{\partial^J F_{\yy}}{\partial \yy^J} = 0. \]
        If we put 
        \[ B_p = \sum_{I} b_{p,I} \aalpha^I, \ B' = \sum_{J} b_{J}' \bbeta^J, \ C_q = \sum_{J} c_{q,J} \bbeta^J, \ C' = \sum_{I} c_{I}' \aalpha^I, \]
        then it follows that $B_p \in \Fperpx \subseteq \kk[\aalpha], C_q \in \Fperpy \subseteq \kk[\bbeta]$ for all $p,q$, and 
        \[ P_0 = \sum_{p} B_p \cdot B' + \sum_{q} C_q \cdot C'. \]
        This finishes the proof of the inclusion $(F_{\xx} \cdot F_{\yy})^{\perp} \subseteq (\Fperpx, \Fperpy)$.
        The other inclusion is obvious by Equation (\ref{equation:applying derivatives on both}).
    \end{proof}
    \begin{lemma}\label{fact:free}
        The $\kk[t, \bbeta]$-algebras $\frac{\kk[t, \aalpha, \bbeta]}{F(t \cdot \xx + \yy)^{\perp}}$ and $\frac{\kk[t, \bbeta]}{\Fperpy}$ are isomorphic.
    \end{lemma}
    \begin{proof}
        Consider the natural $\kk$-algebra homomorphism $\mu : \kk[t, \bbeta] \to \frac{\kk[t, \aalpha, \bbeta]}{F(t \cdot \xx + \yy)^{\perp}}$.
        First we show that this map is onto. Note that
        \begin{equation*}
            (\alpha_i - t \cdot \beta_i) \ \wcirc \ F(t \cdot \xx + \yy) = \left(\frac{\partial}{\partial x_i} - t \cdot \frac{\partial}{\partial y_i}\right)(F(t \cdot \xx + \yy)) = \left(\frac{\partial F}{\partial x_i}\right) (t \cdot \xx + \yy) \cdot t - t \cdot \left(\frac{\partial F}{\partial x_i}\right) (t \cdot \xx + \yy) = 0,
        \end{equation*}
        so $\alpha_i \equiv t \cdot \beta_i \ (mod \ F(t \cdot \xx + \yy)^{\perp})$ for all $i=1, \dots, 6$, and we get surjectivity. We will prove now that the kernel of $\mu$ is $\Fperpy$. Take $f = f(\bbeta, t) \in \ker \mu$. This means that
        \begin{equation*}
            f(\bbeta, t) \ \wcirc \ F(t \cdot \xx + \yy) = 0.
        \end{equation*}
        Let us expand $f(\bbeta, t) = \sum_{p=1}^n f_p(\bbeta) \cdot t^p$. Then
        \begin{equation*}
            f(\bbeta, t) \ \wcirc \ F(t \cdot \xx + \yy) = \sum_{p=1}^n (f_p(\bbeta) \ \wcirc \ F(t \cdot \xx + \yy) ) \cdot t^p = \sum_{p=1}^n (f_p(\bbeta) \ \wcirc \ F_{\yy})(t \cdot \xx + \yy) \cdot t^p.
        \end{equation*}
        by the chain rule. Suppose (for a contradiction) that there exist $p$ such that $f_p(\bbeta) \ \wcirc \ F_{\yy} \neq 0$. Take the smallest such $p$. Then the polynomial in variables $\xx, \yy$ near $t^p$ in $f(\bbeta, t) \ \wcirc \ F(t \cdot \xx + \yy)$ has to be $0$, so $(f_p(\bbeta) \ \wcirc \ F_{\yy})(t \cdot \xx + \yy)$ has all coefficients near $x^p y^j$ zero for $j=1, \dots, 6$. But this is only possible if $f_p(\bbeta) \ \wcirc \ F_{\yy} = 0$ which gives a contradiction. Thus all functions $f_p$ belong to $\Fperpy$ and so $f \in \Fperpy$ which finishes the proof.
    \end{proof}
    This lemma yields an immediate corollary.
    \begin{corollary}\label{cor:freeness}
        The algebra $\frac{\kk[t, \aalpha, \bbeta]}{F(t \cdot \xx + \yy)^{\perp}}$ is a free $\kk[t]$-module, because
        \begin{equation*}
            \frac{\kk[t, \aalpha, \bbeta]}{F(t \cdot \xx + \yy)^{\perp}} \simeq \frac{\kk[t, \bbeta]}{\Fperpy} \simeq \kk[t] \otimes_\kk \frac{\kk[\bbeta]}{\Fperpy}.
        \end{equation*}
    \end{corollary}
    Before stating the next algebraic preparation step, we recall some notion from commutative algebra.
    \begin{remark}\label{def:modulestructure}
        If $C$ is a finite dimensional $\kk[t]$-algebra, then the $\kk[t]$-module $\Hom_{\kk[t]}(C, \kk[t])$ has a $C$-module structure given by $(c \cdot \symbphi)(d) := \symbphi(c \cdot d)$ for $\symbphi \in \Hom_{\kk[t]}(C, \kk[t])$ and $c, d \in C$.
        Moreover, if $C = \frac{\kk[t, \aalpha]}{H^{\perp}}$ for some $H \in \kk[t, \xx]$, then there is a natural $C$-modules homomorphism
        \begin{equation*}
            \Phi: C \to \Hom_{\kk[t]}(C, \kk[t])
        \end{equation*}
        sending $c \in C$ to $\phi_c := \Phi(c)$ given by the formula
        \begin{equation*}
            \phi_c (d) = (ev_{c \ \wcirc \ H}(d))(\Bar{0}).
        \end{equation*}
        The evaluation at $\xx = \Bar{0}$ at the end makes sense as $(ev_{c \ \wcirc \ H}(d)) \in \kk[t, \xx]$ and the result lies indeed in $\kk[t]$.
    \end{remark}
    
    \begin{proposition}\label{prop:selfdual}
        Let us fix $H \in \kk[t, \xx]$ and consider the algebra $C := \frac{\kk[t, \aalpha]}{H^{\perp}}$. Suppose that $C$ is a free $\kk[t]$-module and for all $t_0 \in \kk$ the dimension of the $\kk$-linear space $\frac{\kk[\aalpha]}{H(t_0)^{\perp}}$ is constant.
        Then $\Phi : C \to \Hom_{\kk[t]}(C, \kk[t])$ is a $C$-module isomorphism.
    \end{proposition}
    Before the proof let us cite a lemma which is a translation of \cite[Proposition 2.12]{Jel18} to our context.
    \begin{lemma}\label{lem:selfdual}
        Assume that $\kk$ is algebraically closed. If $C = \frac{\kk[t, \aalpha]}{H^{\perp}}$ for $H \in \kk[t, \xx]$ is such that for all $t_0 \in \kk$ the dimension of the $\kk$-linear space $\frac{\kk[\aalpha]}{H(t_0)^{\perp}}$ is constant, then for all $t_0 \in \kk$ there is an isomorphism between $\frac{\kk[t, \aalpha]}{(H^{\perp}, (t-t_0))}$ and $\frac{\kk[\aalpha]}{H(t_0)^{\perp}}$ induced by the map
        \begin{equation*}
            \kk[t, \aalpha] \to \frac{\kk[\aalpha]}{H(t_0)^{\perp}},
        \end{equation*}
        which maps $f$ to $ev_{H(t_0)}(f)(t_0)$.
    \end{lemma}
    \begin{proof}[Proof of Proposition \ref{prop:selfdual}]
        First we prove that $\ker \Phi = 0$. This is because
        \begin{equation*}
            \Phi(c) = 0 \iff \im ev_{c \ \wcirc \ H} \subseteq (\xx) \iff c \ \wcirc \ H = 0 \iff c = 0 \in \frac{\kk[t, \aalpha]}{H^{\perp}}.
        \end{equation*}
        The middle equivalence holds, because if $c \ \wcirc \ H \neq 0$, then there exists some partial derivative such that after applying it to $c \ \wcirc \ H$ we get a non-zero element in $\kk[t]$.
        Now our goal is to show that $\Phi$ is onto. Suppose for a contradiction that it is not. Let $N := \Hom_{\kk[t]}(C, \kk[t])$. Then we get the following exact sequence of $C$-modules
        \[\begin{tikzcd}
        	0 & C & {N} & E & 0,
        	\arrow[from=1-1, to=1-2]
        	\arrow["\Phi", from=1-2, to=1-3]
        	\arrow[from=1-3, to=1-4]
        	\arrow[from=1-4, to=1-5]
        \end{tikzcd}\]
        where $E \neq 0$ is the quotient $C$-module. This is also a sequence of $\kk[t]$-modules and since $E \neq 0$ is finitely generated, there exists an ideal $\mm = (t - t_0) \subseteq \kk[t]$ (for some $t_0 \in \kk$) such that the quotient $E/(t- t_0) E$ is non-zero (take $\mm$ such that $E_{\mm} \neq 0$ and use the Nakayama's lemma). The ideal $\mm$ is of this form as $\kk$ is algebraically closed. We tensor everything by $\frac{\kk[t]}{(t-t_0)}$ and we get
        \[\begin{tikzcd}
        	& {C/(t-t_0)C} & {N/ (t-t_0)N} & {E/ (t-t_0)E} & 0.
        	\arrow[from=1-2, to=1-3]
        	\arrow[from=1-3, to=1-4]
        	\arrow[from=1-4, to=1-5]
        \end{tikzcd}\]
        To get the final contradiction, we will prove that ${C/(t-t_0)C} \to {N/ (t-t_0)N}$ is an isomorphism. First, by Lemma \ref{lem:selfdual} we get that $C/(t-t_0)C \simeq \frac{\kk[\aalpha]}{H(t_0)^{\perp}}$.
        On the other hand, $N/ (t-t_0)N \simeq \Hom_{\kk}(C/(t-t_0)C, \kk)$. Indeed, if we apply the functor $\Hom_{\kk[t]}(C, -)$ on a short exact sequence
        \[\begin{tikzcd}
        	0 & {(t-t_0)} & {\kk[t]} & {\frac{\kk[t]}{(t-t_0)}} & 0,
        	\arrow[from=1-1, to=1-2]
        	\arrow[from=1-2, to=1-3]
        	\arrow[from=1-3, to=1-4]
        	\arrow[from=1-4, to=1-5]
        \end{tikzcd}\]
        we get
        \[\begin{tikzcd}
        	0 & {\Hom_{\kk[t]}(C,(t-t_0))} & {\Hom_{\kk[t]}(C,\kk[t])} & {\Hom_{\kk[t]}\left(C,\frac{\kk[t]}{(t-t_0)}\right)} & 0,
        	\arrow[from=1-1, to=1-2]
        	\arrow[from=1-2, to=1-3]
        	\arrow[from=1-3, to=1-4]
        	\arrow[from=1-4, to=1-5]
        \end{tikzcd}\]
        where exactness follows from the fact that $C$ is a free $\kk[t]$-module. Now, the last sequence identifies with
        \[\begin{tikzcd}
        	0 & {(t-t_0) \Hom_{\kk[t]}(C,\kk[t])} & {\Hom_{\kk[t]}(C,\kk[t])} & {\Hom_{\kk}\left(C/(t-t_0)C,\frac{\kk[t]}{(t-t_0)}\right)} & 0,
        	\arrow[from=1-1, to=1-2]
        	\arrow[from=1-2, to=1-3]
        	\arrow[from=1-3, to=1-4]
        	\arrow[from=1-4, to=1-5]
        \end{tikzcd}\]
        so we get the desired isomorphism $N/ (t-t_0)N \simeq \Hom_{\kk}(C/(t-t_0)C, \kk)$. Now, note that by its definition, the homomorphism \[\frac{\kk[\aalpha]}{H(t_0)^{\perp}} \simeq {C/(t-t_0)C} \to {N/ (t-t_0)N} \simeq \Hom_{\kk}(C/(t-t_0)C, \kk) \simeq \Hom_{\kk}\left(\frac{\kk[\aalpha]}{H(t_0)^{\perp}}, \kk\right)\]
        is given by the same function $\Phi$ but coming from the Macaulay's duality without the variable $t$. The same argument as in the beginning of this proof shows that the obtained map 
        $\frac{\kk[\aalpha]}{H(t_0)^{\perp}} \to \Hom_{\kk}(\frac{\kk[\aalpha]}{H(t_0)^{\perp}}, \kk)$
        is injective. Its surjectivity follows from the equality of dimensions over $\kk$. But this contradicts the fact that we have an exact sequence
        \[\begin{tikzcd}
        	& {C/(t-t_0)C} & {N/ (t-t_0)N} & {E/ (t-t_0)E} & 0
        	\arrow[from=1-2, to=1-3]
        	\arrow[from=1-3, to=1-4]
        	\arrow[from=1-4, to=1-5]
        \end{tikzcd}\]
        with $E/ (t-t_0)E \neq 0$. This finishes the proof of surjectivity of $\Phi$.
    \end{proof}
    
    We come back to the analysis of the fractal family.
    Since $\Gamma(t) \ \wcirc \ (F_{\Bar{x}} \cdot F_{\Bar{y}}) = F(t \cdot \xx + \yy)$ by Lemma \ref{fact:Gamma}, we have $(F_{\xx} \cdot F_{\yy})^{\perp} \subseteq F(t \cdot \xx + \yy)^{\perp} \subseteq \kk[t, \aalpha, \bbeta]$. We get a short exact sequence
    \[\begin{tikzcd}
    	0 & L & {\frac{\kk[t, \aalpha, \bbeta]}{(F_{\xx} \cdot F_{\yy})^{\perp}}} & {\frac{\kk[t, \aalpha, \bbeta]}{F(t \cdot \xx + \yy)^{\perp}}} & 0,
    	\arrow[from=1-3, to=1-4]
    	\arrow[from=1-2, to=1-3]
    	\arrow[from=1-1, to=1-2]
    	\arrow[from=1-4, to=1-5]
    \end{tikzcd}\]
    where $L$ is the kernel of the natural projection. 
    This sequence of $\frac{\kk[t, \bbeta]}{\Fperpy}$-modules splits, as the right hand side module is a free $\frac{\kk[t, \bbeta]}{\Fperpy}$-module (isomorphic to $\frac{\kk[t, \bbeta]}{\Fperpy}$ by Lemma \ref{fact:free}). After applying the functor $\Hom_{\kk[t]}(-, \kk[t])$ we obtain a sequence of $\kk[t]$-modules
    \[\begin{tikzcd}
    	0 & {\Hom_{\kk[t]}(L, \kk[t])} & {\Hom_{\kk[t]}\left(\frac{\kk[t, \aalpha, \bbeta]}{(F_{\xx} \cdot F_{\yy})^{\perp}},\kk[t]\right)} & {\Hom_{\kk[t]}\left(\frac{\kk[t, \aalpha, \bbeta]}{F(t \cdot \xx + \yy)^{\perp}},\kk[t]\right)} & 0.
    	\arrow[from=1-4, to=1-3]
    	\arrow[from=1-3, to=1-2]
    	\arrow[from=1-2, to=1-1]
    	\arrow[from=1-5, to=1-4]
    \end{tikzcd}\]
    However, by Remark \ref{def:modulestructure} \[{\Hom_{\kk[t]}\left(\frac{\kk[t, \aalpha, \bbeta]}{(F_{\xx} \cdot F_{\yy})^{\perp}},\kk[t]\right)}\] has a $\frac{\kk[t, \aalpha, \bbeta]}{(F_{\xx} \cdot F_{\yy})^{\perp}}$-module structure and \[{\Hom_{\kk[t]}\left(\frac{\kk[t, \aalpha, \bbeta]}{F(t \cdot \xx + \yy)^{\perp}},\kk[t]\right)}\] has a $\frac{\kk[t, \aalpha, \bbeta]}{F(t \cdot \xx + \yy)^{\perp}}$-module structure. In particular both of these $\kk[t]$-modules have a structure of $\frac{\kk[t, \bbeta]}{\Fperpy}$-module. Thus the preceding sequence can be considered as a sequence of $\frac{\kk[t, \bbeta]}{\Fperpy}$-modules.
    If we fix a splitting on the first sequence
    \begin{equation*}
        \frac{\kk[t, \aalpha, \bbeta]}{(F_{\xx} \cdot F_{\yy})^{\perp}} \simeq L \oplus \frac{\kk[t,  \aalpha, \bbeta]}{F(t \cdot \xx + \yy)^{\perp}},
    \end{equation*}
    we get the induced splitting of $\frac{\kk[t, \bbeta]}{\Fperpy}$-modules
    \begin{equation}\label{eq:blablabla}
        \Hom_{\kk[t]}\left(\frac{\kk[t, \aalpha, \bbeta]}{(F_{\xx} \cdot F_{\yy})^{\perp}}, \kk[t]\right) \simeq \Hom_{\kk[t]}(L, \kk[t]) \oplus \Hom_{\kk[t]}\left(\frac{\kk[t,  \aalpha, \bbeta]}{F(t \cdot \xx + \yy)^{\perp}}, \kk[t]\right).
    \end{equation}
    Consider the $\kk$-algebra $\frac{\kk[t, \aalpha, \bbeta]}{(F_{\xx} \cdot F_{\yy})^{\perp}}$. It is a free $\kk[t]$-module, because 
    \[\frac{\kk[t, \aalpha, \bbeta]}{(F_{\xx} \cdot F_{\yy})^{\perp}} = \frac{\kk[t, \aalpha, \bbeta]}{(\Fperpx, \Fperpy)} \simeq \frac{\kk[\aalpha, \bbeta]}{(\Fperpx, \Fperpy)} \otimes_\kk \kk[t],\]
    where the first equality follows from Lemma \ref{fact:fracfamilyform}. Moreover, the polynomial $(F_{\xx} \cdot F_{\yy})$ does not depend on $t$, so the assumptions of Proposition \ref{prop:selfdual} are satisfied. Thus, the left hand side of equation \eqref{eq:blablabla} is isomorphic to $\frac{\kk[t, \aalpha, \bbeta]}{(F_{\xx} \cdot F_{\yy})^{\perp}}$ as a $\frac{\kk[t, \bbeta]}{\Fperpy}$-module, so it is a free $\frac{\kk[t, \bbeta]}{\Fperpy}$-module. This means that $\Hom_{\kk[t]}(L, \kk[t])$ is a projective $\frac{\kk[t, \bbeta]}{\Fperpy}$-module (as it is a direct factor in a free $\frac{\kk[t, \bbeta]}{\Fperpy}$-module). We will now prove that the \Mf is flat by showing it is isomorphic to $\Hom_{\kk[t]}(L, \kk[t])$ as a $\frac{\kk[t, \bbeta]}{\Fperpy}$-module.
    \begin{proposition}\label{prop:isotoff}
        $\Hom_{\kk[t]}(L, \kk[t])$ is isomorphic to the \Mf from Definition \ref{def:gammat}, as a $\frac{\kk[t, \bbeta]}{\Fperpy}$-module.
    \end{proposition}
    \begin{proof}
        Let $C := \frac{\kk[t, \aalpha, \bbeta]}{(F_{\xx} \cdot F_{\yy})^{\perp}}$, $C_0 := \frac{\kk[t, \aalpha, \bbeta]}{F(t \cdot \xx + \yy)^{\perp}}$ and $q: C \to C_0$ be the quotient map. From Proposition \ref{prop:selfdual} we know that $\Hom_{\kk[t]}(C, \kk[t]) \simeq C$, $\Hom_{\kk[t]}(C_0, \kk[t]) \simeq C_0$ by homomorphisms as in Remark \ref{def:modulestructure}. Thus in order to prove the proposition, it is enough to interpret the image of $\Hom_{\kk[t]}(C_0, \kk[t])$ inside $C \simeq \Hom_{\kk[t]}(C, \kk[t])$ as the ideal $(\Gamma(t))$. 
        
        Consider $c_0 \in C_0$ and let $\Bar{c}_0 \in C$ be a representative of $c_0$. We will calculate the image of $c_0$ coming from the morphism $C_0 \simeq \Hom_{\kk[t]}(C_0, \kk[t]) \to \Hom_{\kk[t]}(C, \kk[t]) \simeq C$. First, look at the diagram
        \[\begin{tikzcd}
        	0 & L & C & {C_0} & 0 \\
        	&&& \kk[t].
        	\arrow[from=1-1, to=1-2]
        	\arrow[from=1-2, to=1-3]
        	\arrow["q"', from=1-3, to=1-4]
        	\arrow[from=1-4, to=1-5]
        	\arrow["{\symbphi_{c_0}}", from=1-4, to=2-4]
        	\arrow["{\symbphi_{c_0} \circ q}"', from=1-3, to=2-4]
        \end{tikzcd}\]
        Here $\symbphi_{c_0}$ comes from Remark \ref{def:modulestructure}, so it is equal to $z \circ ev_{c_0 \ \wcirc \ F(t \cdot \xx + \yy)}$ where $z : \kk[t, \xx, \yy] \to \kk[t]$ is the map sending $\xx$ and $\yy$ to $0$. We claim that
        \begin{equation*}
            \symbphi_{c_0} \circ q = \symbphi_{\Bar{c}_0 \cdot \Gamma(t)}.
        \end{equation*}
        Indeed, that is because:
        \begin{equation*}
            \symbphi_{c_0} \circ q = z \circ ev_{c_0 \ \wcirc \ F(t \cdot \xx + \yy)} \circ q = z \circ ev_{\Bar{c}_0 \cdot \Gamma(t) \ \wcirc \ (F_{\xx} \cdot F_{\yy})} = \symbphi_{\Bar{c}_0 \cdot \Gamma(t)}.
        \end{equation*}
        So the induced homomorphism $C_0 \simeq \Hom_{\kk[t]}(C_0, \kk[t]) \to \Hom_{\kk[t]}(C, \kk[t]) \simeq C$ sends $c_0 \mapsto \Bar{c}_0 \cdot \Gamma(t)$ (and this doesn't depend on the choice of $\Bar{c}_0$). Hence the cokernel of this homomorphism is isomorphic to $C/(\Gamma(t))$ so to the \Mf by Lemma \ref{fact:fracfamilyform}.
    \end{proof}
    
    This yields the following.
    
    \begin{corollary}\label{cor:flatness_of_fractal_family}
        The \Mf is a flat $\frac{\kk[t, \bbeta]}{\Fperpy}$-module.
    \end{corollary}
    \begin{proof}
        From Proposition \ref{prop:isotoff} we know that \Mf is isomorphic to $\Hom_{\kk[t]}(L, \kk[t])$ which is a projective $\frac{\kk[t, \bbeta]}{\Fperpy}$-module by equation \eqref{eq:blablabla} and the fact that $\frac{\kk[t, \aalpha, \bbeta]}{(F_{\xx} \cdot F_{\yy})^{\perp}}$ satisfies the assumptions of Proposition \ref{prop:selfdual}. In particular, the \Mf is flat.
    \end{proof}
    
\subsection{Explicit description of the fiber $V$}\label{subsection:explicit}
    Let us recall the setup we are working in. We fix a general enough cubic $F$ (see Definition \ref{def:nice}) in six variables $\xx$, and $g, Q_1, \dots, Q_6 \in \kk[\aalpha]$ satisfying the assumptions from Notation \ref{notation:relative_Macaulay} (this is possible as $F$ is general enough). In particular, we will use Definition \ref{def:gammat}, so
    \begin{equation*}
            \Gamma(t) = g(\aalpha) + t \cdot \sum_{1 \leq i \leq 6} \beta_i \cdot Q_i(\aalpha) + t^2 \cdot \sum_{1 \leq i \leq 6} \alpha_i \cdot Q_i(\bbeta) + t^3 \cdot g(\bbeta),
    \end{equation*}
    and the fractal family is a $\frac{\kk[t, \bbeta]}{\Fperpy}$ - algebra $M = \frac{\kk[t, \aalpha, \bbeta]}{(\Fperpx, \Fperpy, \Gamma(t))}$ (depending on $F, g$ and $Q_i$'s). Also, $[I]$ denotes a point of $\cH = \mathrm{Hilb}^{13}_{\Affine^6/\kk}$ constructed from $F$ by taking $I = (\Fperp, g)$. We fix symbols $i, s, \pi$ for the maps
    \[\begin{tikzcd}
    	{\cB^{+}} & {\cB} \\
    	{\cB^{\G_m}}
    	\arrow["\pi", from=1-1, to=2-1]
    	\arrow["i"', from=1-1, to=1-2]
    	\arrow["s", bend left, from=2-1, to=1-1]
    \end{tikzcd}\]
    coming from the \BBname{} decomposition of the barycenter scheme $\cB \subseteq \cH$. We denote by $V$ the fiber $\pi^{-1}([I])$ of $\pi$ over $[I]$. In this subsection, in Theorem \ref{thm:completeV}, we prove that $V$ is isomorphic to $\Spec\left(\frac{\kk[t, \bbeta]}{\Fperpy}\right)$ and this isomorphism can be realised using the fractal family. The idea is to look at the morphism to $\cH$ defined by the fractal family and see what it does on tangent vectors. Recall that in Proposition \ref{prop:subspacew} we calculated the tangent space $T_{[I]}V$ by identifying it with a subspace $\nameY \subseteq T_{[I]}\cH = \Hom_{S}(I, S/I)$ spanned by homomorphisms
    \[ \symboly_i = \symbtan_i - \frac{1}{13} \partial_i, \]
    for $i=1, \dots, 6$, where $\symbtan_i$ maps $(\Fperp)_2$ to $0$ and $g$ to $Q_i$, and $\partial_i$ maps $f$ to $\frac{\partial f}{\partial \alpha_i}$ (see Example \ref{ex:additiveaction}).
    We equip $\Affine^1 \times \Spec\left(\frac{\kk[\bbeta]}{\Fperpy}\right)$ with the $\mathbb{G}_m$-action given by the grading of $\Spec\left(\frac{\kk[t,\bbeta]}{\Fperpy}\right)$, where all variables $t, \beta_1, \dots, \beta_6$ have degree one.
    Similarly, we fix the $\G_m$-action on $\Spec\left(\frac{\kk[\bbeta]}{\Fperpy}\right)$ coming from the induced action on the locus $t=0$. We also freely use actions on $\cH$ and $\cB$ introduced in Section \ref{subsection:groupsactiongonH}.
    We start with the following.

    \begin{lemma}
        The morphism 
        \[ v:\Affine^1 \times \Spec\left(\frac{\kk[\bbeta]}{\Fperpy}\right) \to \cH \]
        induced by the \Mf $M$ is $\Gmult$-equivariant.
    \end{lemma}
    \begin{proof}
        This is by Proposition \ref{prop:torusequivariance}.
    \end{proof}
    
    Let $\bar{v} : \Spec\left(\frac{\kk[\bbeta]}{\Fperpy}\right) \to \cH^+$ be the morphism induced from $v:\Affine^1 \times \Spec\left(\frac{\kk[\bbeta]}{\Fperpy}\right) \to \cH$, as in Definition \ref{def:xplusfunctorial}. We write $v|_{ \{1\} } : \Spec\left(\frac{\kk[\bbeta]}{\Fperpy}\right) \to \cH$ for the restriction of $v$ to $\{1\} \times \Spec\left(\frac{\kk[\bbeta]}{\Fperpy}\right)$. 
    
    By Proposition \ref{prop:hilbandbar} we have $\cH \simeq \Affine^n \times \cB$ and by Proposition \ref{prop:barisgminv} the projection $\namep : \cH \to \cB$ is $\G_m$-equivariant. Hence, the composition $u = \namep \circ v$
    \[\begin{tikzcd}
    	{\Affine^1 \times \Spec\left(\frac{\kk[\bbeta]}{\Fperpy}\right)} & \cH & \cB
    	\arrow["v", from=1-1, to=1-2]
    	\arrow["\namep", from=1-2, to=1-3]
    	\arrow["u", bend right, from=1-1, to=1-3]
    \end{tikzcd}\]
    is $\G_m$-equivariant. The morphism $u$ induces a morphism $\bar{u} : \Spec\left(\frac{\kk[\bbeta]}{\Fperpy}\right) \to \cB^+$ fitting in the diagram
    \[\begin{tikzcd}
    	{\Spec\left(\frac{\kk[\bbeta]}{\Fperpy}\right)} & {\cB^+} & {\cB^{\G_m}} \\
    	{\Affine^1 \times\Spec\left(\frac{\kk[\bbeta]}{\Fperpy}\right)} & \cB.
    	\arrow["i"', from=1-2, to=2-2]
    	\arrow["u"', from=2-1, to=2-2]
    	\arrow["{\bar{u}}", from=1-1, to=1-2]
    	\arrow["{\{ 1 \} \times id}"', hook, from=1-1, to=2-1]
    	\arrow["\pi", from=1-2, to=1-3]
    \end{tikzcd}\]
    From now on, the goal is to prove that $\bar{u}$ is an isomorphism onto $V \subseteq \cB^+$. We first need the following fact.
    \begin{lemma}\label{fact:factorisationV}
        The composition $\pi \circ \bar{u}$ is a constant morphism onto the point $[I]$.
    \end{lemma} 
    \begin{proof}
        Consider the composition $w := \pi_{\cH} \circ \bar{v}$. Note that the following diagram commutes
        \[\begin{tikzcd}
        	{\Spec\left(\frac{\kk[\bbeta]}{\Fperpy}\right)} & {\cH^+} & {\cB^+} \\
        	& {\cH^{\G_m}} & {\cB^{\G_m},}
        	\arrow["{\namep^+}", from=1-2, to=1-3]
        	\arrow["{\bar{v}}", from=1-1, to=1-2]
        	\arrow["{\pi_\cH}", from=1-2, to=2-2]
        	\arrow["\pi", from=1-3, to=2-3]
        	\arrow["{\namep^{\G_m}}"', from=2-2, to=2-3]
        	\arrow["{\bar{u}}", bend left, from=1-1, to=1-3]
        	\arrow["w"', from=1-1, to=2-2]
        \end{tikzcd}\]
        where $\namep^+$ and $\namep^{\G_m}$ are morphism induced by $\namep$ on \BBname{} decomposition functors and fix-point functors respectively. It is enough to prove that $w : \Spec\left(\frac{\kk[\bbeta]}{\Fperpy}\right) \to \cH^{\G_m}$ is a constant map onto $[I]$, because $\namep^{\G_m}([I]) = [I] \in \cB^{\G_m}$.
        
        By the functorial description from Section \ref{subsection:BBdecomp} of the map $\pi$, we see that $w = \pi_{\cH} \circ \bar{v} = v|_{ \{0\} \times \Spec\Big(\frac{\kk[\bbeta]}{\Fperpy}\Big) }$. As $v$ is given by the \Mf $M$ its restriction to zero is given by $M/(t)$ so by the algebra
        \[ M_0 := \frac{\kk[\aalpha, \bbeta]}{(\Fperpx, \Fperpy, \Gamma(0))}, \]
        but $\Gamma(0) = g(\aalpha)$. Hence $M_0$ is the product family
        \[ M_0 \simeq \frac{\kk[\aalpha]}{(\Fperpx, g(\aalpha))} \otimes_\kk \frac{\kk[\bbeta]}{(\Fperpy)}. \]
        As $I = (\Fperpx, g(\aalpha))$ this family corresponds to the constant morphism onto $[I]$. This finishes the proof.
    \end{proof}
    By Lemma \ref{fact:factorisationV} the morphism $\bar{u}:\Spec\left(\frac{\kk[\bbeta]}{\Fperpy}\right) \to \cB^+$ factors through $V \subseteq \cB^+$, because by its definition $V = \pi^{-1}([I])$. We will write $\bar{u}:\Spec\left(\frac{\kk[\bbeta]}{\Fperpy}\right) \to V$ for the factorisation. We would like to prove that $\bar{u}$ is an isomorphism. In order to do so, first we check that it is $\G_m$-equivariant.
    \begin{lemma}\label{fact:gmequibaru}
        The morphism $\bar{u}:\Spec\left(\frac{\kk[\bbeta]}{\Fperpy}\right) \to V$ is $\G_m$-equivariant with respect to the natural $\G_m$-action on $\Spec\left(\frac{\kk[\bbeta]}{\Fperpy}\right)$ and the usual $\G_m$-action on $\cH$.
    \end{lemma}
    \begin{proof}
        First we argue that the morphism $v|_{ \{1\} } : \Spec\left(\frac{\kk[\bbeta]}{\Fperpy}\right) \to \cH$ is $\G_m$-equivariant.
        The \Mf restricted to $t=1$ is given by the ideal $(\Fperpx, \Gamma(1)) \subseteq S_{\frac{\kk[\bbeta]}{\Fperpy}}$. The part $\Fperpx$ is homogeneous with respect to the standard grading on $S$ and $\Gamma(1)$ is homogeneous with respect to the grading with $\deg(\alpha_i) = \deg(\beta_i)=1$.
        Thus, by Proposition \ref{prop:torusequivariance} we get the $\G_m$-equivariance of $v|_{\{1\}} : \Spec\left(\frac{\kk[\bbeta]}{\Fperpy}\right) \to \cH$.
        
        Now we prove that this suffices for $\G_m$-equivariance of $\bar{u}$. Consider the following commutative diagram
        \[\begin{tikzcd}
        	{\Spec\left(\frac{\kk[\bbeta]}{\Fperpy}\right)} & {\cH^+} & {\cB^+} \\
        	& \cH.
        	\arrow["{\namep^+}", from=1-2, to=1-3]
        	\arrow["{\bar{v}}", from=1-1, to=1-2]
        	\arrow["{i_\cH}", from=1-2, to=2-2]
        	\arrow["{\bar{u}}", bend left, from=1-1, to=1-3]
        	\arrow["{v|_{ \{1\} }}"', from=1-1, to=2-2]
        \end{tikzcd}\]
        We know that $v|_{ \{1\} }, i_{\cH}, \namep^+$ are $\G_m$-equivariant. Thus, the $\G_m$-equivariance of $\bar{v}$ will finish the proof. Denote by $\mu : \G_m \times \Spec\left(\frac{\kk[\bbeta]}{\Fperpy}\right) \to \Spec\left(\frac{\kk[\bbeta]}{\Fperpy}\right)$ the $\G_m$-action on $\Spec\left(\frac{\kk[\bbeta]}{\Fperpy}\right)$ and $\nu, \nu^+$ the $\G_m$-action on $\cH, \cH^+$ respectively. We look at the following diagram
        \[\begin{tikzcd}
        	{\G_m \times \Spec\left(\frac{\kk[\bbeta]}{\Fperpy}\right)} & {\G_m \times \cH^+} & {\G_m \times\cH} \\
        	{\Spec\left(\frac{\kk[\bbeta]}{\Fperpy}\right)} & {\cH^+} & \cH.
        	\arrow["\mu", from=1-1, to=2-1]
        	\arrow["{id \times \bar{v}}", from=1-1, to=1-2]
        	\arrow["{\nu^+}"', from=1-2, to=2-2]
        	\arrow["{\bar{v}}", from=2-1, to=2-2]
        	\arrow["{id \times i_{\cH}}", from=1-2, to=1-3]
        	\arrow["\nu", from=1-3, to=2-3]
        	\arrow["{i_{\cH}}", from=2-2, to=2-3]
        	\arrow["{id \times v|_{ \{1\} }}", bend left, from=1-1, to=1-3]
        	\arrow["{v|_{ \{1\} }}", bend right, from=2-1, to=2-3]
        \end{tikzcd}\]
        We want to show that $\nu^+ \circ (id \times \bar{v}) = \bar{v} \circ \mu$. However, by Fact \ref{fact:mono} $i_{\cH}$ is a monomorphism of functors, so this is equivalent to the fact that \[i_{\cH} \circ \nu^+ \circ (id \times \bar{v}) = i_{\cH} \circ \bar{v} \circ \mu.\] By the commutativity of the diagram this translates to $\nu \circ (id \times v|_{ \{1\} }) = v|_{ \{1\} } \circ \mu$. But this is the $\G_m$-equivariance of $v|_{ \{1\} }$ so we are done.
    \end{proof}
    Now we investigate the behaviour of tangents under $v|_{ \{1\} }$.
    \begin{proposition}\label{prop:ending}
        Let $\{\beta_i^{\vee}\}_{1 \leq i \leq 6} \subseteq T_{ (\bbeta) } \Spec\left(\frac{\kk[\bbeta]}{\Fperpy}\right)$ be the dual basis to $\{\beta_1, \dots, \beta_6\} \subseteq (\bbeta)/(\bbeta)^2 = T_{ (\bbeta) }^{\vee} \Spec\left(\frac{\kk[\bbeta]}{\Fperpy}\right)$.
        The differential $d v|_{\{1\}} : T_{ (\bbeta) } \Spec\left(\frac{\kk[\bbeta]}{\Fperpy}\right) \to T_{[I]} \cH$ sends $\beta_i^{\vee}$ to $- \symbtan_i$ from Definition \ref{def:subspacew}.
    \end{proposition}
    \begin{proof}
        By definition $\beta_i^{\vee} \in T_{ (\bbeta) } \Spec\left(\frac{\kk[\bbeta]}{\Fperpy}\right)$ is the morphism $\Spec(\kk[\eps]) \to \Spec\left(\frac{\kk[\bbeta]}{\Fperpy}\right)$ that takes $\beta_j$ to $0$ for $j \neq i$ and $\beta_i$ to $\eps$, for $i=1, \dots, 6$. 
        Let us denote $M' = \frac{\kk[\aalpha, \bbeta]}{(\Fperpx, \Fperpy, \Gamma(1))}$.
        By composing with $v|_{\{1\}}$ we get a morphism $\Spec(\kk[\eps]) \to \cH$.
        This morphism is represented, in the sense of Definition \ref{def:Hilbert:first}, by the pullback of the family $\Spec(M')$ over $\Spec\left(\frac{\kk[\bbeta]}{\Fperpy}\right)$ to $\Spec(\kk[\eps])$.
        Thus, it is represented by the family
        \[\begin{tikzcd}
        	{\kk[\eps]} & {\frac{\kk[\aalpha, \eps]}{(\Fperpx, g(\aalpha) + \eps \cdot Q_i(\aalpha))}}.
        	\arrow[from=1-1, to=1-2]
        \end{tikzcd}\]
        By Example \ref{ex:tangentsw} this family corresponds to the tangent $- \symbtan_i = - \symbtan_{Q_i}$.
    \end{proof}
    As a corollary, we prove the main result of Subsection \ref{subsection:explicit}. For the convenience of the reader, we recall here the discussion from the beginning of Subsection \ref{sub:descriptionofv}. We argued that $V$ is affine, so $V = \Spec B$ for some graded $\kk$-algebra of finite type $B = \bigoplus_{n \geq 0} B_n$ with $B_0 = \kk$ and the maximal ideal $\mm := \bigoplus_{n \geq 1} B_n$. We fixed six homogeneous representatives of generators of $\frac{\mm}{\mm^2}$ called $\symboly_1^{\vee}, \dots, \symboly_6^{\vee} \in \mm$, forming the dual basis to $\symboly_1, \dots, \symboly_6 \in T_{[I]}V$. We also defined a graded surjective homomorphism $\kk[\symb_1, \dots, \symb_6] \to B$ sending $\symb_i$ to $\symboly_i^{\vee}$, whose kernel was denoted by $J$.
    \begin{theorem}\label{thm:completeV}
        The morphism $\bar{u}:\Spec\left(\frac{\kk[\bbeta]}{\Fperpy}\right) \to V$ is an isomorphism. Moreover:
        \[ J = (\Fperp)(\Bar{\symb}) := \{ g(\symb_1, \dots, \symb_6) \in \kk[\symb_1, \dots, \symb_6] : g(\alpha_1, \dots, \alpha_6) \in \Fperp \}. \]
    \end{theorem}
    \begin{proof}
        Consider the diagram
        \[\begin{tikzcd}
        	& {} \\
        	{\Spec\left(\frac{\kk[\bbeta]}{\Fperpy}\right)} & {\cH^+} & {\cB^+} \\
        	{\Affine^1 \times \Spec\left(\frac{\kk[\bbeta]}{\Fperpy}\right)} & \cH & \cB.
        	\arrow["v"', from=3-1, to=3-2]
        	\arrow["{i_{\cH}}", from=2-2, to=3-2]
        	\arrow["{\bar{v}}", from=2-1, to=2-2]
        	\arrow[hook, from=2-1, to=3-1]
        	\arrow["\namep", bend left, from=3-2, to=3-3]
        	\arrow["i"', from=2-3, to=3-3]
        	\arrow["{\namep^+}", from=2-2, to=2-3]
        	\arrow["u", bend right, from=3-1, to=3-3]
        	\arrow["{r^+}", bend left, from=2-3, to=2-2]
        	\arrow["r", from=3-3, to=3-2]
        	\arrow["{v|_{ \{1\} }}", from=2-1, to=3-2]
        	\arrow["{\bar{u}}"', bend left, from=2-1, to=2-3]
        \end{tikzcd}\]
        where $r$ is the embedding of $\cB$ in $\cH$ and morphisms $\namep^+, r^+$ come from $\namep$ and $r$ respectively, by Definition \ref{def:xplusfunctorial}. We are interested in the morphism
        \[ q := (r \circ i \circ \bar{u}) : \Spec\left(\frac{\kk[\bbeta]}{\Fperpy}\right) \to \cH. \] We calculate
        \begin{equation*}
            q = r \circ i \circ \bar{u} =  r \circ i \circ \namep^+ \circ \bar{v} = r \circ \namep \circ i_{\cH} \circ \bar{v} = r \circ \namep \circ v|_{ \{1\} }.
        \end{equation*}
        Thus, the differential $dq : T_{ (\bbeta) } \Spec\left(\frac{\kk[\bbeta]}{\Fperpy}\right) \to T_{[I]} \cH$ is the composition of $dv|_{ \{1\} } : T_{ (\bbeta) } \Spec\left(\frac{\kk[\bbeta]}{\Fperpy}\right) \to T_{[I]} \cH$ and the projection $d\namep:T_{[I]} \cH \to T_{[I]} \cB \subseteq T_{[I]} \cH$.
        By Proposition \ref{prop:ending} if $\beta_i^{\vee}$ is as in the statement of this proposition, then $dv|_{ \{1\} } (\beta_i^{\vee}) = - \symbtan_i$. Note that 
        \begin{equation*}
            - \symbtan_i = - \left(\symbtan_i - \frac{1}{13} \partial_i\right) - \frac{1}{13} \partial_i = - \symboly_i - \frac{1}{13} \partial_i
        \end{equation*}
        by Definition \ref{def:subspacey}. By Proposition \ref{prop:subspacew}, $\symboly_i \in T_{[I]} V \subseteq T_{[I]} \cB$, so the projection of $dv|_{ \{1\} } (\beta_i^{\vee})$ onto $T_{[I]} \cB$ is $- \symboly_i$, because $T_{[I]} \cH = T_{[I]} \cB \oplus \langle \partial_1, \dots, \partial_6 \rangle$. Hence
        \begin{equation*}
            dq(\beta_i^{\vee}) = - \symboly_i \textnormal{ for $i = 1, \dots, 6$.}
        \end{equation*}
        After identifying $V$ with the image of the closed embedding given by the composition \eqref{eq:compo} we get a factorisation
        \begin{equation}\label{eq:last}
            q : \Spec\left(\frac{\kk[\bbeta]}{\Fperpy}\right) \to V \subseteq \cH.
        \end{equation}
        Recall that at the beginning of Subsection \ref{sub:descriptionofv} we presented $V$ as $\Spec\left(\frac{\kk[\Bar{\symb}]}{J}\right)$ and $(\symb_i)_{1 \leq i \leq 6}$ was a basis of $T_{[I]}^{\vee} V$ dual to $(\symboly_i)_{1 \leq i \leq 6}$. Thus, by equation \eqref{eq:last}, we get that
        \begin{equation}\label{eq:lastlast}
            (dq)^{\vee}(\symb_i) = - \beta_i \textnormal{ for $i = 1, \dots, 6$.}
        \end{equation}
        Moreover, by Lemma \ref{fact:gmequibaru} $\bar{u}$ is $\G_m$-equivariant, so $q$ is $\G_m$-equivariant as well.
        In this way we see that $q : \Spec\left(\frac{\kk[\bbeta]}{\Fperpy}\right) \to \Spec\left(\frac{\kk[\Bar{\symb}]}{J}\right) \simeq V$ yields a graded homomorphism of algebras
        \[\begin{tikzcd}
        	{\frac{\kk[\bbeta]}{\Fperpy}} & {\frac{\kk[\Bar{\symb}]}{J}}
        	\arrow["q^{\#}", from=1-2, to=1-1]
        \end{tikzcd}\]
        and by equation \eqref{eq:lastlast} it is defined by $q^{\#}(\symb_i) = - \beta_i$ for $i=1, \dots, 6$. Thus if $f(\Bar{\symb}) \in J$ then $f(-\bbeta) \in \Fperpy$ and because $J$ and $\Fperpy$ are homogeneous it follows that $J \subseteq \Fperp(\Bar{\symb})$. From Proposition \ref{prop:boundj} we get that $\Fperp(\Bar{\symb}) \subseteq J$ so these two ideals are in fact equal. Hence, $q = r \circ i \circ \bar{u}$ is an isomorphism. The morphism $i$ is an isomorphism on some open neighbourhood of $[I]$ by Corollary \ref{cor:bplusisb} and $r$ is a closed embedding. Thus, $\bar{u}$ is an isomorphism onto $V \subseteq \cB^+$.
    \end{proof}
    
\section*{Acknowledgements}
The author would like to thank Joachim Jelisiejew for supervising this project and many helpful discussions. 
Moreover, many thanks go to all referees that gave comments to this text.
The author was supported by the NCN grant 2017/26/D/ST1/00913 while being a MSc student at the University of Warsaw.
This material is based upon work supported by the National Science Foundation under Grant No. DMS-2424441, and by the IAS School of Mathematics.


\begin{thebibliography}{{\v{S}}VdB18}

\bibitem[BB73]{ABB1} Andrzej Białynicki-Birula. Some theorems on actions of algebraic groups. Ann. of Math. (2) 98 (1973), 480–497.

\bibitem[BB76]{ABB2} Andrzej Białynicki-Birula. Some properties of the decompositions of algebraic varieties determined by actions of a torus. Bull. Acad. Polon. Sci. Sér. Sci. Math. Astronom. Phys. 24 (1976), no. 9, 667–674.

\bibitem[BCR19]{example1} Cristina Bertone, Francesca Cioffi and Margherita Roggero. Smoothable Gorenstein Points Via Marked Schemes and Double-generic Initial Ideals, Experimental Mathematics, 1-18, 2019, Taylor \& Francis.

\bibitem[Ber08]{notesbertin} Jose Bertin. The punctual Hilbert scheme: an introduction. \newblock \url{https://www-fourier.ujf-grenoble.fr/sites/ifmaquette.ujf-grenoble.fr/files/bertin_rev.pdf}

\bibitem[Bri77]{Bri77} Joël Briançon. 
Description de $\Hilb_n C\{{ x,y \}}$.
Invent. Math. 41 (1977), no. 1, 45–89.

\bibitem[Coh46]{Cohen} Irvin S. Cohen. On the structure and ideal theory of complete local rings. Trans. Amer. Math. Soc. 59 (1946), 54–106.

\bibitem[DJ+17]{example2} Theodosios Douvropoulos, Joachim Jelisiejew, Bernt Ivar Utstøl Nødland and Zach Teitler. The Hilbert scheme of 11 points in $\Bbb A^3$ is irreducible. Combinatorial algebraic geometry, 321–352, Fields Inst. Commun., 80, Fields Inst. Res. Math. Sci., Toronto, ON, 2017.

\bibitem[Eis05]{Eis1} David Eisenbud. The geometry of syzygies. A second course in commutative algebra and algebraic geometry. Graduate Texts in Mathematics, 229. Springer-Verlag, New York, 2005. xvi+243 pp. ISBN: 0-387-22215-4

\bibitem[ES87]{ES87} Geir Ellingsrud and Stein Arild Strømme. On the homology of the Hilbert scheme of points in the plane. Invent. Math. 87 (1987), no. 2, 343–352.

\bibitem[EI78]{EI78} Jacques Emsalem and Anthony Iarrobino.
Some zero-dimensional generic singularities; finite algebras having small tangent space.
Compositio Math. 36 (1978), no. 2, 145–188.

\bibitem[FGI+05]{FGI+05} Barbara Fantechi, Lothar Göttsche, Luc Illusie, Steven L. Kleiman, Nitin Nitsure and Angelo Vistoli. Fundamental algebraic geometry, volume 123 of Mathematical Surveys and Monographs. American Mathematical Society, Providence, RI, 2005. Grothendieck’s FGA explained.

\bibitem[Fog68]{Fogarty} John Fogarty. (1968). Algebraic Families on an Algebraic Surface. American Journal of Mathematics, 90(2), 511-521. doi:10.2307/2373541

\bibitem[Hai01]{Haiman} Mark Haiman.  Hilbert schemes, polygraphs and the Macdonald positivity conjecture. J. Amer. Math. Soc. 14 (2001), no. 4, 941–1006.

\bibitem[Har10]{DefHar} Robin Hartshorne.  Deformation theory. Graduate Texts in Mathematics, 257. Springer, New York, 2010. viii+234 pp. ISBN: 978-1-4419-1595-5

\bibitem[Iar72a]{Iar72a} Anthony Iarrobino.
Punctual Hilbert schemes.
Bull. Amer. Math. Soc. 78 (1972), 819–823.

\bibitem[Iar72b]{Iar72b} Anthony Iarrobino. Reducibility of the families of 0-dimensional schemes on a variety, Invent. Math. 15 (1972), 72–77.

\bibitem[Jel17]{Jel17} Joachim Jelisiejew (2017). Classifying local Artinian Gorenstein algebras. Collect. Math. 68 (2017) no. 1, 101-127.

\bibitem[Jel18]{Jel18} Joachim Jelisiejew. VSPs of cubic fourfolds and the Gorenstein locus of the Hilbert scheme of 14 points on $\Bbb{A}^6$. Linear Algebra Appl. 557 (2018), 265–286.

\bibitem[Jel19]{elemcomp} Joachim Jelisiejew. Elementary components of Hilbert schemes of points. J. Lond. Math. Soc. (2) 100 (2019), no. 1, 249–272.

\bibitem[Jel20]{Hilbpath} Joachim Jelisiejew. Pathologies on the Hilbert scheme of points. Invent. Math. 220 (2020), no. 2, 581–610.

\bibitem[Jel24a]{Jel24a} Joachim Jelisiejew. Open problems in deformations of Artinian algebras, Hilbert schemes and around, Contemporary Mathematics (AMS) volume Deformation of Artinian Algebras and Jordan Type, 2024.

\bibitem[Jel24b]{Jel24b} Joachim Jelisiejew. Generically nonreduced components of Hilbert schemes on fourfolds (for the memory volume of Gianfranco Casnati).
Rendiconti del seminario matematico di Torino, Vol. 82 No. 1 (2024).

\bibitem[JKS25]{jelisiejew2023behrend} Joachim Jelisiejew, Martijn Kool and Reinier F. Schmiermann. Behrend's function is not constant on $\mathrm{Hilb}^n(\mathbb{C}^3)$, Geometry \& Topology 29-8 (2025), 4469--4476.

\bibitem[JS19]{JS19} Joachim Jelisiejew and Łukasz Sienkiewicz. Białynicki-Birula decomposition for reductive groups. J. Math. Pures Appl. (9) 131 (2019), 290–325.

\bibitem[JŠ22]{JS21} Joachim Jelisiejew and Klemen Šivic (2021). Components and singularities of Quot schemes and varieties of commuting matrices, Journal für die reine und angewandte Mathematik (Crelle's Journal), vol. 2022, no. 788, 2022, pp. 129-187.

\bibitem[MNO+06]{maulik2006gromov} Davesh Maulik, Nikita Nekrasov, Andrei Okounkov and Rahul Pandharipande. Gromov-Witten theory and Donaldson-Thomas theory, I. Compositio Mathematica 142 (2006), no. 5, 1263-1285.

\bibitem[Osz26]{PiotrOszer2026} Piotr Oszer (2026). Deformations of the connected sum of Gorenstein algebras, \newblock \url{https://arxiv.org/abs/2601.03179}.

\bibitem[Ric24]{Ricolfi24} Andrea T. Ricolfi. A sign that used to annoy me, and still does. Journal of Geometry and Physics 195 (2024): 105032.
   
\bibitem[RS23]{Ramkumar_2023} Ritvik Ramkumar and Alessio Sammartano. On the parity conjecture for Hilbert schemes of points on threefolds, Ann. Sc. Norm. Super. Pisa Cl. Sci. (5) Vol. XXVI (2025), 1665-1676.

\bibitem[Sha90]{Sha90} Igor Rostislavovich Shafarevich. Deformations of commutative algebras of class 2, Algebra i Analiz, 2:6 (1990), 178–196; Leningrad Math. J., 2:6 (1991), 1335–1351.

\bibitem[Sha01]{Shafarevich} Igor Rostislavovich Shafarevich. Degeneration of semisimple algebras. Special issue dedicated to Alexei Ivanovich Kostrikin. Comm. Algebra 29 (2001), no. 9, 3943–3960.

\bibitem[Sta25]{Stacks} The Stacks project authors. The Stacks project. \newblock  \url{https://stacks.math.columbia.edu}

\bibitem[Vak17]{Vakil} Ravi Vakil, THE RISING SEA Foundations of Algebraic Geometry. \newblock  \url{http://math.stanford.edu/~vakil/216blog/FOAGnov1817public.pdf}

\end{thebibliography}
\end{document}